%------------------------------------------------------------------------------
% nak-new.tex  
%       Ringel: 
% Beginn: Oct-31-2020. Version Oct 12, 2021   
%------------------------------------------------------------------------------
\magnification=\magstep1   
\input amstex
\UseAMSsymbols
\input pictex 
%\hoffset=0truecm \voffset=0truecm 
%\nopagenumbers
\vsize=23truecm
\NoBlackBoxes
\parindent=18pt
  
   \font\rmk=cmr8      

\font\gross=cmbx10 scaled\magstep1 
%\hrule height 2pt \vskip 3pt \hrule \bigskip\bigskip    
\def\char{\operatorname{char}}

\def\mod{\operatorname{mod}}

\def\Hom{\operatorname{Hom}}

\def\Ext{\operatorname{Ext}}

\def\rad{\operatorname{rad}}

\def\Cok{\operatorname{Cok}}
\def\soc{\operatorname{soc}}

\def\top{\operatorname{top}}

\def\pd{\operatorname{pd}}
\def\id{\operatorname{id}}

\def\domdim{\operatorname{dom-dim}}

  \def\ss{\ssize }

\def\s{\hfill \square}

%%%%%%%%%%%%%%%%%%%%%%%%%%%%%%%%%%%%%%%%%%%%%%%%%%%%%%%%%%%%%%%%%% 
%%%%%%%%%%%%%%%%%%%%%%%%%%%%%%%%%%%%%%%%%%%%%%%%%%%%%%%%%%%%%%%%%%%
%%%%%%%%%%%%%%%%%%%%%%%%%%%%%%%%%%%%%%%%%%%%%%%%%%%%%%%%%%%%%%%%%%%%%
%%%%%%%%%%%%%%%%%%%%%%%%%%%%%%%%%%%%%%%%%%%%%%%%%%%%%%%%%%%%%%%%%%
\centerline{\gross Linear Nakayama algebras}
	\medskip
\centerline{\gross which are higher Auslander algebras.}
                     \bigskip
\centerline{Claus Michael Ringel}
                \bigskip\medskip
\noindent {\narrower Abstract:  \rmk An artin algebra $\ss A$ is said to be a higher Auslander
algebra 
provided the global dimension and the dominant dimension coincide.
We say that a linear Nakayama algebra is concave, provided its Kupisch
series first increases, then decreases. We are going to classify the concave Nakayama
algebras which are higher Auslander algebras. 
Let us stress that the classification strongly depends
on the parity of the global dimension of $\ss A$.
	\medskip
\noindent 
\rmk Key words. Nakayama algebra. Higher Auslander algebra. 
	\medskip
\noindent 
2020 Mathematics Subject classification. 
  Primary 16G10,   Secondary  16G20, 16G70, 16E05, 16E10, 16E65. 
\par}
	\medskip\bigskip
%%%%%%%%%%%%%%%%%%%%%%%%%%%%%%%%%%%%%%%%%%%%%%%%%%%%%%%%%%%%%%%%%%%%%
{\bf 0. Introduction.} 
	\medskip
Following Iyama [2], 
an artin algebra $A$ is said to be a {\it higher Auslander algebra}
provided the global dimension and the dominant dimension of $A$ are equal. 

In the present note, let $A$ be a 
Nakayama algebra with $n = n(A)$ simple modules and with global dimension $d = d(A).$
For simplicity, we assume that
$A$ is a $k$-algebra given by a quiver with relations, where $k$ is a field.
The modules $M$ to be considered are usually left $A$-modules with finite length $|M|$;
we denote by $PM$ and $IM$ a projective cover and an injective envelope of $M$, respectively,
and $\soc M,\ \rad M,$ and $\top M$ denote the socle, the radical and the top of $M$,
respectively. We write $\Omega M$ for the first syzygy module of $M$; this is the kernel
of a projective cover $PM \to M$ and $\pd M$ denotes the projective dimension of $M$;
for the zero module, we use the convention $\pd 0 = -\infty$. 
An indecomposable module $M$ is said to be {\it even}, or {\it odd},
if $\pd M$ is even, or odd, respectively. We write $\Sigma M$ for the first suspension
module of $M$: it is the cokernel of an injective envelope $M \to IM$.
	\medskip
The maximal length of indecomposable modules will be called the {\it height} $h(A)$ 
of $A$, and the indecomposable modules of length $h(A)$ are called {\it summits} (they have to
be both projective and injective). 
If $A$ is a linear Nakayama algebra, we denote by $\omega_A$ its 
simple injective module.

It is well-known that if $A$ has finite global
dimension, then $d(A)\le 2n(A)-2.$ Of course, if $A$ is a linear Nakayama algebra, then
$A$ has global dimension $d(A) < n(A)$. 
A quite surprising recent 
paper [5] by Madsen, Marczinzik and Zaimi shows that 
for any numbers $n\le d\le 2n-2,$ there is a unique (necessarily cyclic) 
Nakayama algebra $A$ with $n = n(A)$ and $d = d(A)$ which is
a higher Auslander algebra. 
An explanation of this fact 
has been given by Sen [7] using the powerful $\epsilon$-reduction which he
has introduced in previous papers. Sen's reduction shows:
in order to classify the Nakayama algebras which are higher 
Auslander algebras, it remains to consider the linear ones.

The purpose of the
present note is to analyze the linear Nakayama algebras which are higher
Auslander algebras of global dimension $d$. 
It turns out that it is decisive to look at the parity of $d$. 
We are going to provide a complete system 
of invariants for the concave Nakayama algebras
which are higher Auslander algebras (a linear Nakayama algebra is said to be 
{\it concave} provided
its Kupisch series is first increasing, then decreasing); some of them have been
exhibited already in [8] and [1]. 
If $A$ is a concave algebra, there is a unique indecomposable  module $P$
of length $h(A)$
such that $\rad P$ is projective. We say that $P$ is the {\it first summit} of $A$.
	\medskip
If $A$ is a Nakayama algebra and $M$ is an odd module, with
composition series $0 = M_0 \subset M_1 \subset \cdots M_m = M,$ let
$\char M = (\pd M_1/M_0,\dots, \pd M_m/M_{m-1}).$ We call $\char M$ 
the {\it characteristic} of $M$. 
In addition, we look also at the zero module and define $\char 0$ to be the empty
sequence.
	\medskip 
A semisimple algebra has global dimension zero or $-\infty$ and dominant
dimension $\infty$, therefore {\it a higher Auslander algebra is never semisimple.} 
Thus, if $A$ is a higher Auslander Nakayama algebra and $P$ is a summit, then
$|P| \ge 2,$ thus $\soc P \subseteq \rad P.$
	\medskip
{\bf Theorem 1.} {\it Let $A$ be a concave Nakayama algebra with first summit $P$. 
If $A$ is a higher Auslander algebra of odd
global dimension $d$,  let $Z(A) = \rad P/\soc P.$ 
Then $Z(A)$ is
an odd module or the zero module 
and $\char Z(A)$ is a decreasing sequence of odd
numbers which are bounded by $d$. 

In this way, we obtain a bijection between the isomorphism classes of the 
concave higher Auslander Nakayama algebras $A$ of odd
global dimension $d$ and the decreasing sequences of odd numbers bounded by $d$,
by sending $A$ to $\char Z(A)$.}
	\medskip
In particular, Theorem asserts that for 
any decreasing sequence $d \ge c_1\ge c_2 \ge \cdots \ge c_u$
of odd numbers, there is a (unique) concave higher Auslander algebra 
$H = H_d(c_1,\dots,c_u)$ of
global dimension $d$ with $\char Z(H) = (c_1,\dots,c_u).$ The algebra
$H = H_d(\emptyset)$ is the linear Nakayama algebra $H$ with $h(H) = 2$ and $n(H) = d+1.$ 
Here is a sketch of $\mod H_d(c_1,\dots,c_u)$ in case $u\ge 1$:
$$
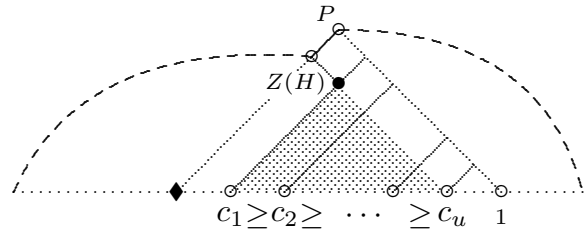
{\beginpicture
    \setcoordinatesystem units <.3cm,.3cm>
%%%%%%%%%%%%%%%%%%%%%%%%%%%%%%%%%%%%%%%%%%%%%%%%%%
%\put{$\mod H_d(c_1,\dots,c_u)$} at -4 6
\multiput{} at -1 0  15  4  /
\plot 5 5  6 6  /
\setdots <.3mm>
\plot 2 0  7 5  /
\plot 4 0  8 4  /
\plot 8 0  10 2 /
\plot 10 0  11 1 /
\plot 5 5  6 4 /
\setdashes <1mm>
\setsolid
\multiput{$\circ$} at  5 5  6 6 /
\multiput{$\bullet$} at  6 4  /
\setdots <1mm>
\plot -6 0  15 0 /
\setdashes <1mm>
\setquadratic
\plot -6 0  -2 4.1  5 5 / 
\plot 6 6  12 4.5  15 0 /
%\put{$P$} at 5.6 6.6 
\put{$\ss P$} at 5.5 6.6  
\put{$\ss Z(H)$} at 4.4 4  
\setlinear
\setshadegrid span <.4mm>
\vshade 2 0 0 <z,z,,> 6 0 4 <z,z,,> 10 0 0 /
\put{$c_1$} at 2 -0.9
\put{$c_2$} at 4 -0.9
\put{$c_{u}$} at 10.2 -0.9
\put{$\ss 1$} at 12 -0.9
\multiput{$\ge$} at 3 -0.9  5 -.9  9 -.9 /
\put{$\cdots$} at 7 -.9
\setdots <.5mm>
\plot 0 0  3.5 3.5 /
\plot 4.5 4.5  5 5 /
\plot 6 6  12 0 /

%\plot 6 4  10 0 /
\put{$\ss\blacklozenge$} at 0 0
\multiput{$\circ$} at 2 0  4 0  8 0  10 0  12 0 /
\endpicture}
$$
The bullet is $Z(H) = \rad P/\soc P$; it is an odd module 
with $\char Z(H) = (c_1,\dots,c_u).$
The shaded part are the subfactors of $Z(H).$
We have added below any odd composition factor of $P$ its projective dimension
(since $P$ is the first summit of $H$, the module $\rad P$ is projective, thus
$\pd \top P = 1$). The socle of $P$ 
has even projective dimension, namely $c_1-1$, and is marked by a 
black lozenge $\ss\blacklozenge$. 
	\bigskip

{\bf Theorem 1$'$.} {\it Let $A$ be a concave Nakayama algebra with first summit $P$. 
If $A$ is a higher Auslander algebra of even
global dimension $d$,  let $Z'(A) = \rad^2 P.$ 
Then $Z'(A)$ is
an odd module or the zero module and $\char Z'(A)$ is a strictly increasing sequence of odd
numbers bounded by $d$. 

In this way, we obtain a bijection between the isomorphism classes of the 
concave higher Auslander Nakayama algebras $A$ of even
global dimension $d$ and the strictly increasing sequences of odd numbers bounded by $d$,
by sending $A$ to $\char Z'(A)$.}
	\medskip
As in the case of $d$ being odd, Theorem 1$'$ asserts that for $d$ even and 
any strictly increasing sequnece $c_1\ge c_2 \ge \cdots \ge c_u$
of odd numbers bounded by $d$, there is a (unique) concave higher Auslander algebra 
$H = H_d(c_1,\dots,c_u)$ of
global dimension $d$ now with $\char Z'(H) = (c_1,\dots,c_u).$
The algebra
$H = H_d(\emptyset)$ is the linear Nakayama algebra $H$ with $h(H) = 2$ and $n(H) = d+1.$ 
Here is a sketch of $\mod H_d(c_1,\dots,c_u)$ in case $u \ge 1$:
$$
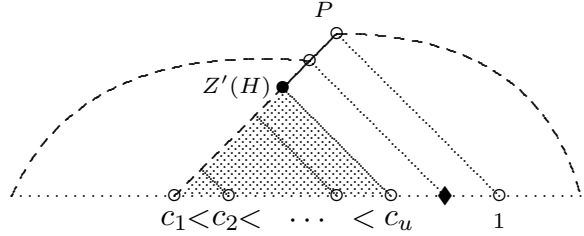
{\beginpicture
    \setcoordinatesystem units <.3cm,.3cm>
%%%%%%%%%%%%%%%%%%%%%%%%%%%%%%%%%%%%%%%%%%%%%%%%%%
\put{$\ss P$} at 5.5 6.9  
\put{$\ss Z'(H)$} at 2.3 4.05  
\multiput{} at -1 0  15  4  /
\plot 4 4  6 6  /
\setdots <.3mm>
\plot 4 4  8 0 /
\plot 3 3  6 0 /
\plot 1 1  2 0 /

\setdashes <1mm>
\plot 0 0  4 4 /
\setsolid 
\multiput{$\circ$} at   6 6  5 5 /
\multiput{$\bullet$} at 4 4  /
\setdots <1mm>
\plot -6 0  15 0 /
\setdashes <1mm>
\setquadratic
\plot -6 0  -2 4.1  5 5 / 
\plot 6 6  12 4.5  15 0 /
\setlinear
\setshadegrid span <.4mm>
\vshade 0 0 0 <z,z,,> 4 0 4 <z,z,,> 8 0 0 /
\put{$c_1$} at 0 -0.9
\put{$c_2$} at 1.9 -0.9
\put{$c_{u}$} at 8.3 -0.9
\put{$\ss 1$} at 12 -0.9
\multiput{$<$} at .9 -0.85  2.8 -.85  7 -.85 /
\put{$\cdots$} at 5 -.9
\setdots <.5mm>
\plot 5 5  10 0 /

\plot 6 6  12 0 /
\put{$\ss\blacklozenge$} at 10 0
\multiput{$\circ$} at 0 0  2 0  6 0  8 0  12 0 /
\endpicture}
$$
The bullet is $Z'(H) = \rad^2 P$; it is an odd module with 
$\char Z'(H) = (c_1,\dots,c_u).$
The shaded part are the subfactors of
$Z'(H).$
Again, we have added below any odd composition factor of $P$ its projective dimension
(since $P$ is the first summit of $H$, the radical $\rad P$ of $P$ is projective, thus
$\pd \top P = 1$). 
The composition factor 
$\top \rad P$ has even projective dimension, namely $c_u+1$, and has been 
marked by a black lozenge $\ss\blacklozenge$. 
	\bigskip
The algebras $H_d(c_1,\dots,c_u)$ can be constructed directly quite easily, 
as the following Theorems 2 and 2$'$ show. 
In order to describe these constructions, we need some further considerations.
	\medskip
If $A$ is a Nakayama algebra of finite global dimension, we attach to {\bf every} indecomposable
module $M$ (not only the odd ones), its characteristic sequence $\char M$, see 1.3.
If $M$ has length $m$, then $\char M = (z_1,\dots,z_m)$ is a sequence of $m$
non-negative numbers, with all but
at most one of the numbers $z_i$ being odd (the numbers $z_i$ exhibit the projective dimension of $M$ itself or of composition factors of $M$).

Starting with a sequence $\bold z = (z_1,\dots,z_m)$ of non-negative numbers,
with one of the numbers zero, whereas the remaining numbers are odd, 
there is a concave Nakayama algebra $A(\bold z)$ 
of height $m$ with $\char P(\omega_A) = \bold z,$ called the ascent
algebra of $\bold z$, see Proposition 2.1.

For any linear Nakayama algebra $A$, and $d$ a positive integer, 
the partial $d$-closure $C_d(A)$ of $A$ will be defined in section 3.
	\bigskip
{\bf Theorem 2.} {\it Let $d$ be odd. 
Let $c_1\ge c_2 \ge \cdots \ge c_u$ be odd numbers bounded by $d$.
Let $H = H_d(c_1,\dots,c_u)$ be the partial $d$-closure of $A(0,c_1,\dots,c_u,1)$.
Then $H$
is a concave higher Auslander algebra of global dimension $d$ with $\char Z(H) = (c_1,\dots,c_u).$}
	\medskip
{\bf Theorem 2$'$.} {\it Let $d$ be even.  
Let $c_1 < c_2 < \cdots < c_u$ be odd numbers bounded by $d$.
Let $H = H_d(c_1,\dots,c_u)$ be the partial $d$-closure of $A(c_1,\dots,c_u,0,1)$.
Then $H$
is a concave higher Auslander algebra of global dimension $d$ with $\char Z'(H) = (c_1,\dots,c_u).$}
	\bigskip
Until now, we have considered only 
the first summit of the algebras $H_d(c_1,\dots,c_u).$ 
Given a concave Nakayama algebra, there is a unique indecomposable module $Q$ of length $h(A)$ 
such that $Q/\soc Q$ is injective. We say that $Q$ is the {\it last summit} of $A$
and call $Q/\soc Q$ the {\it principal cliff module} of $A$.
 
	\medskip
{\bf Theorem 3.} 
{\it Let $d$ be odd.
Let $c_1 \ge c_2 \ge \cdots \ge c_u$ be odd numbers bounded by $d$.
Let $H = H_d(c_1,\dots,c_u)$. Let $P$ be the first summit of $H$ and $Q$ the last summit.

If $u = 0$, then $H$ has $d$ summits, 
$\char P = (0,1)$, and $\char Q = (0,d).$

If $u \ge 1$, let $t = \frac12(d-c_1).$ 
Then $t\ge 0$, the algebra $H$ has $(u+2)t+1$ summits, and }
$$
 \char P = (0,c_1,\dots,c_u,1), \quad
 \char Q = (0,c_1\!+\!2t,\dots,c_{u}\!+\!2t,1\!+\!2t).
$$
	\medskip
{\bf Theorem 3$'$.} 
{\it Let $d$ be even.
Let $c_1 < c_2 < \cdots < c_u$ be odd numbers bounded by $d$.
Let $H = H_d(c_1,\dots,c_u).$ Let $P$ be the first summit of $H$ and $Q$ the last summit.

If $u = 0$, then $H$ has $d$ summits, $\char P = (0,1),$ and $\char Q = (d\!-\!1,0).$

If $u \ge 1$, let $t = \frac12(d-c_u-1).$ 
Then $t\ge 0$,  $H$ has $(u+2)t+u$ summits, and}
$$
 \char P = (c_1,\dots,c_u,0,1),\quad
 \char Q = (d\!-\!1,0,1\!+\!2t,c_1\!+\!2t+2,\cdots,c_{u-1}\!+\!2t+2).
$$
  	\bigskip 
%%%%%%%%%%%%%%%%%%%%%%%%%%%%%%%%%%%%%%%%%%%%%%%%%%%
{\bf Outline of the paper.}
	\smallskip
If $M_0,\dots,M_m$ are indecomposable modules with $\Hom(M_{i-1},M_i) \neq 0$ for
$1\le i \le m$, we say that $M_0$ is a {\it predecessor} of $M_m$ and $M_m$ a 
{\it successor} of $M_0.$ As usual, we denote by $\tau$ the Auslander-Reiten translation
(with $\tau P = 0$ for $P$ projective). If $M$ is an indecomposable module, let
$\Cal F(M)$ be the set of the (isomorphism classes of) non-zero subfactors of $M$.
A {\it Serre subcategory} of $\mod A$ is the full subcategory given by the extension closure of a 
set of simple modules. The Serre subcategory {\it generated by a set of modules} is the Serre subcategory given by the composition factors of these modules.
If $A$ is a Nakayama algebra, the modules with a fixed socle are said to form a {\it ray};
those with a fixed top are said to form a {\it coray.}
	\medskip
In section 1, we introduce and analyze the characteristic sequence $\char M$ of 
an indecomposable module $M$. This seems to be of independent 
interest. As we will see in 1.6 (see also 4.8),
this sequence determines the projective dimension of 
all subfactors of $M$. 
		\medskip
	
A linear Nakayama algebra $A$ is called {\it ascending} provided for 
indecomposable projective modules  $P, P'$ with $\Hom(P,P') \neq 0$, we have
$|P| \le |P'|,$ or, equivalently, provided $P|(\tau S)| \le |PS|$ for any simple
module $S$. Note that a 
linear Nakayama algebra $A$ is ascending iff $A$ is concave and $P(\omega_A)$
is a summit.
In section 2, we construct for any projective characteristic sequence 
$\bold z$ its ascent algebra
$A = A(\bold z)$; it is an ascending Nakayama algebra $A$ with $\char P(\omega_A) = \bold z$.
In section 3, we construct the partial $d$-closure of a Nakayama algebra.

Section 5 is devoted to Nakayama algebras which are higher Auslander algebras.
We show that if $A$ is a higher Auslander Nakayama algebra of global dimension $d$ and
$S$ is a simple module, then $S$ is torsionless or else $\pd S = \pd IS = d.$

Section 7 contains the proof of Theorems 1, 2, and 3;
section 8 the proof of Theorems 1$'$, 2$'$, and 3$'.$
In both sections, we provide various characterizations of the concave
higher Auslander algebras $A$ of odd or even global dimension, respectively, with
reference to the existence of suitable modules $M$
of length $h(A)-1$. 
The modules $M$ which we single out are, on the one hand, the radical of the first
summit (and certain $\tau^{-}$-shifts), and, on the other hand, the principal cliff
module (and certain $\tau$-shifts).
	\medskip
{\bf Left rotations, right rotations.} 
The main tool used in the paper are the left rotations (see section 2) 
and the right rotations (see setion 3), both are defined for characteristic sequences.
If $A$ is a concave Nakayama algebra and $P$ a summit, left rotations $\lambda$ 
are used in order to reconstruct all predecessors of $P$: here, using $\lambda$, 
we obtain from $PS$, where $S$ is a simple module, the module $P(\tau S)$.
The right rotations are used
in order to construct successors of a given indecomposable module. 
If $I$ is an indecomposable module which is both projective and injective, then
the right rotation $\rho$ sends $\rad I$ to $I/\soc I$. 

Let $A$ be concave and $P$ a summit. Whereas the predecessors of $P$ can be reconstructed
from $P$ using $\lambda$, it is not possible in general, 
to reconstruct the successors of $\char P$. However, this is possible in case the 
algebra $A$ is $d$-bound as defined in section 6 
(using now the right rotation $\rho$). 
All algebras which we are interested in, are $d$-bound: either by
construction (the algebras $H_d(c_1,\dots,c_u)$), or by assumption (the higher Auslander
algebras).
	\medskip
{\bf Memory quivers.}
Throughout the paper, we work with ``memory functions'' on translation quivers.
Special memory quivers will be discussed in section 4. The use of
characteristic sequences in section 1 is already a kind of shadow.

Let $\Gamma$ be a translation quiver. 
A function $\mu\:\Gamma_0 \to \Bbb N\cup\{\infty\}$ is called a {\it memory function}
provided there is an artin algebra $C$ with a Serre subcategory $\Cal C \subseteq \mod C$ 
such that $\Gamma$ is the Auslander-Reiten quiver of $\Cal C$ and 
we have $\mu(X) = \pd {}_CX$ for any indecomposable object $X\in \Cal C.$ 
The pair $(\Gamma,\mu)$ will be called a {\it memory translation quiver} or just
a {\it memory quiver.} In this paper, we only will deal with memory functions with values
in $\Bbb N.$
	\medskip
{\bf The parameter $d$.} 
Note that there are no higher Auslander algebras of global dimension $0$ (since
semisimple algebras have infinite dominant dimension). 
It is easy to see that a connected higher Auslander algebra of global dimension at most $1$ is
the path algebra of a quiver of type $\Bbb A_n$ (without relations).
Thus, we usually may assume that $d\ge 2.$ 
Since for a Nakayama algebra $A$ the injective envelope of an indecomposable projective
module is an indecomposable projective module, the dominant dimension of $A$ is at least 
$1$. 

Iyama, who has introduced the higher Auslander algebras, called 
a higher Auslander algebra $A$ of global dimension $d$ a $(d\!-\!1)$-Auslander algebra,
see [3]. 
Of course, this index shift is well-thought since in this way the Auslander algebras
themselves are just the $1$-Auslander algebras.
(There is the similar, and of course related, problem concerning the 
calibration of the representation dimension: a shift by 1 has the advantage that then
the value $1$ (instead of $2$) is attached to the non-semisimple representation-finite algebras.)

However, we have refrained from following this convention (as also several other
authors do), since in our discussion of the higher Auslander algebras of global dimension $d$,
we want to stress that it is the number $d$ (and not $d-1$) and its parity, which are decisive.

	\bigskip\medskip 
%%%%%%%%%%%%%%%%%%%%%%%%%%%%%%%%%%%%%%%%%%%%%%%%%%%%%%%%%%%%%%%%%%%%%
{\bf 1. The characteristic sequence of an indecomposable module.} 
	\medskip
For any indecomposable module $M$, we are going to introduce 
its characteristic sequence $\char M$;
it is a sequence $(z_1,\dots,z_m)$ of non-negative integers with at most one even entry,
where $m = |M|$. In addition, let $\char 0$ be the empty sequence.  
The relevance of $\char M$ 
relies on Madsen's maximum principle 1.1 and 
the subfactor formula 1.4 which shows that $\char M$ determines the 
projective dimension of all subfactors of $M$.
	\medskip
{\bf 1.1. Proposition (Madsen's maximum principle).} {\it Let $A$ be a Nakayama algebra. 
Let $M$ be an indecomposable module.
\item{\rm(1)} At most one of the composition factors of $M$ is even.
\item{\rm(2)} If all composition factors $F_iM$ of $M$ are odd, then $\pd M = \max_i \pd F_iM;$ 
   thus
   $M$ is odd. Conversely, if $M$ is odd, then all its composition factors are odd. 
\item{\rm(3)} If at least one composition factor of $M$ is even, then $M$ is even.
\item{\rm(4)} If $S$ is simple, with $\tau^iS$ being odd for $0 \le i \le m$, then
   $|PS| \ge m+2.$\par}
	\medskip
In case $M$ is even, we should stress that
the projective dimension of the odd composition factors of $M$ do not determine the
projective dimension of $M$, see 1.5.
For completeness, we recall from [4] that {\it the existence of an even simple
module implies that $A$ has finite global dimension;} this is the reason that (2) implies (3).
	\medskip
Proof of proposition. For (2) and (3), see [4], or also [6].
	\smallskip
Thus, let us look at (1).
Assume that $M$ has a submodule $U$ with $\top U$ even and that 
$M/U$ has only odd composition factors. 
Since $S = U/\rad U$ is even, $\Omega S = \rad PS$ is odd, thus, according to (3),
all composition factors of $\rad PS$ are
odd. 
We have $U = PS/\Omega U,$ and $\rad U = \rad PS/\Omega U$. 
Since all composition factors of $\rad PS$ are odd, 
all composition factors of $\rad U$ are odd. We see that $U/\rad U$ is the only
even composition factor of $M$.
	\smallskip
Proof of (4).  Assume that $|PS| \le m+1$. The composition factors of $PS$ are
the modules $\tau^iS$ with $0 \le i \le |PS|-1,$ thus they are odd. According to the maximum
principle, it follows that $PS$ is odd. But $\pd PS = 0.$
$\s$
	\medskip
%%%%%%%%%%%%%%%%%%%%%%%%%%%%%%%%%%%%%%%%%%%%%%%%%%%%+
{\bf 1.2. Further consequences.} Let us add some consequences of 1.1  which one 
should keep in
mind. Recall that an indecomposable 
module $M$ is said to be {\it multiplicity-free} provided any simple
module occurs at most once as a composition factor of $M$. 
	\medskip
{\it Let $A$ be a Nakayama algebra.}

\item{(1)} {\it If $S$ is simple and not projective, at most one of $S$, $\tau S$ can be even.}
\item{(2)} {\it If $A\neq 0$, then not all simple modules are odd.}
\item{(3)} {\it An indecomposable odd module is multiplicity-free.}
\item{(4)} {\it If $A$ has finite global dimension and rank $n$, 
  then the height of $A$ is bounded by $2n-1$.} 
	\medskip
Proofs. (1) If $S, \tau S$ both would be even, 
the middle term of the corresponding Auslander-Reiten
sequence would be an indecomposable module with two even composition factors, in contrast ti 
1.1.(1).

(2) If $A \neq 0$, there is an indecomposable projective module $P$. Since
$\pd P = 0$, we see that $P$ is even. According to 1.1.(2), not all composition
factors of $M$ can be odd. 

(3) If $M$ is indecomposable and not multiplicity-free, then all simple modules occur as
composition factors. But if $M$ is odd, the even simple modules do not occur as composition
factors.

(4) If an indecomposable module $M$ has length at least $2n$, then every simple module 
occurs in $M$ with multiplicity at least 2. According to (2), there is an even simple module $S$,
and (1) asserts that $S$ occurs at most once as a composition factor of $M$. This contradiction
shows that any indecomposable module has length at most $2n-1$. 
$\s$
	\bigskip
	
{\bf 1.3. Characteristic sequences.} 
A sequence $(z_1,\dots,z_m)$ of non-negative numbers will be said to be a
{\it characteristic sequence} provided at most one of the entries is even.
A characteristic sequence $(z_1,\dots,z_m)$ with one entry being zero is said to be {\it projective.}

	\medskip
{\bf 1.4. The characteristic of an indecomposable module.}
Assume that $A$ has finite global dimension. 
To every indecomposable module we attach a sequence of
numbers, called its characteristic. On a first reading, the choice and the
ordering of the numbers may seem to be rather curious, but we hope to convince
the reader that the definition is natural. The characteristic  $\char M$
of an
indecomposable module of length $m$ will be a sequence of $m$
non-negative numbers, with all but at most one being odd (and all numbers
being positive iff $M$ is not projective). 

Let $M$ be an indecomposable module of length $m$,
say with composition series 
$$
 0 = M_0 \subset M_1 \subset \cdots \subset M_m = M,
$$
thus $F_i(M) = M_i/M_{i-1}$ (with $1\le i \le m$) are the composition factors of $M$.
	
We define $\char M$ as a sequence of $n$ non-negative integers $(\char M)_i$ with $1\le i \le m$
$$
 (\char M)_i = \left\{\matrix \pd F_i(M) & \text{if} & F_i(M)\ \text{is odd}, \cr
                        \pd M     &  \text{if} & F_i(M)\ \text{is even}.   \cr
                      \endmatrix   \right.     
$$
In addition, we look also at the zero module and define $\char 0 = (\emptyset)$ 
(the empty sequence). The sequence $\char M$ will be called the {\it characteristic of $M$}.
	\medskip 
{\bf Lemma.} 
{\it Let $A$ be a Nakayama algebra of finite global dimension.

\item{\rm(1)} The characteristic of an indecomposable module is a characteristic
sequence.

\item{\rm(2)} An indecomposable module is projective iff its characteristic 
is projective. 

\item{\rm(3)} An indecomposable module is odd iff its characteristic has only 
odd entries.} 
	\medskip
Proof. (1) At most one of the coefficients of $\char M$ is even, see 1.1 (1).

(2) If $(\char M)_i$ is even, then by definition $(\char M)_i = \pd M$. 
If all coefficients of $\char M$ are odd, then the maximum principle for odd modules
asserts that $\pd M = \max_i\char M_i$. 
As a consequence, if $M$ is non-zero and projective, not all coefficients of
$\char M$ can be odd, since otherwise $\pd M = 0$, but $\max_i(\char M_i) > 0$. Thus
$(\char M)_i$ is even, for some $i$, and then $(\char M)_i = \pd M = 0$. 
Conversely, if $(\char M)_i = 0$, for some $i$, 
then $\pd M = (\char M)_i = 0$, thus $M$ is projective (and
non-zero). 

(3) See 1.1.(2). 
$\s$
	\bigskip 
{\bf 1.5.} Given a projective characteristic sequence $(z_0,\dots,z_m)$ of
length $m+1$, let 
$$
 Y(z_0,\dots,z_m) = \left\{
     \matrix (z_1,\dots,z_m) & \text{\it if} & z_0 = 0, \cr
             (z_1,\dots,z_{v-1},z_0+1,z_{v+1},\dots,z_m) &\text{\it if} & z_v = 0, & v\ge 1.
     \endmatrix \right.
$$
This is a non-projective characteristic sequence of length $m$.

Given a non-projective characteristic sequence $(z_1, \dots, z_m)$ of length $m$, let
$$
 P(z_1,\dots,z_m) = \left\{
     \matrix (0,z_1,\dots,z_m) & \quad\ \text{\it  if all $z_i$ are odd,}  \cr
             (z_v\!-\!1,z_1,\dots,z_{v-1},0,z_{v+1},\dots,z_m) &\text{\it if $z_v$ is even.}
     \endmatrix \right.
$$
This is a projective characteristic sequence of length $m+1$. 
	\medskip
{\it The maps $Y(-)$ and $P(-)$ provide bijections between the set of projective characteristic
sequences and the set of non-projective characteristic sequences which are inverse 
to each other.}
$\s$
	\medskip
{\bf Lemma.} 
{\it Let $A$ be a Nakayama algebra of finite global dimension.
Let $P$ be indecomposable projective and $Y = P/\soc P$. If $\char P =
(z_0,\dots,z_m)$, then $\char Y = Y(z_0,\dots,z_m).$}
	\medskip
Proof. First, let $z_0 = 0.$ 
Then the numbers $z_i$ with $1\le i \le m$ are odd, and are the projective
dimension of the composition factors of $Y$. 

Second, let $z_v = 0$ for some $v \ge 1$.
Then $P/\soc P$ has an even composition factor, thus $\soc P$ is an odd composition factor,
therefore the first entry of $\char P$ is $z_0 = \pd \soc P.$ 
Since $\Omega Y = \soc P$, we have $\pd Y = z_0+1.$ It follows that  for $i\neq v$, we have
$(\char Y)_i = \pd F_iY = \pd F_{i+1}P = z_i$ and
$(\char Y)_v = \pd Y = z_0+1.$ 
$\s$
	\bigskip
%%%%%%%%%%%%%%%%%%%%%%%%%%%%%%%%%%%%%%%%%%%%%%%%
{\bf 1.6. The subfactor formula.} In 1.1, we have mentioned Madsen's maximum principle: 
{\it The projective dimension of an odd
module is the maximum of the projective dimension of its composition factors.}
Here is a corresponding formula for even modules.
	\medskip
{\bf Proposition.} {\it Let $A$ be a linear Nakayama algebra. 
Let $M$ be an indecomposable module. Let
$M_1 \subset M_2 \subseteq M$ be submodules and assume that $M_2/M_1$ is even.
Then $M_1$ and $M/M_2$ are odd or zero, and} 
$$
 \pd M_2/M_1 = \max\{\pd M_1+1,\pd M, \pd M/M_2-1\}.
$$
	\medskip
Proof. According to 1.1, $M$ has precisely one even composition factor, and this has
to be a subfactor of $M_2/M_1.$ In particular, we see that $M_1$ and $M/M_2$ are odd
or zero. 
	\smallskip
Case 1: $M_2 = M.$ If, in addition, $M$ is projective, then
the assertion is clear: if $M_1 = 0$, then $M_2/M_1 = M$ is a non-zero projective module
(thus both the left hand side and the right hand side are equal to $0$), 
whereas if $M_1 \neq 0,$ then $\Omega (M/M_1) = M_1$ shows that $\pd M/M_1 = \pd M_1 + 1.$

Thus, let $M_2 = M$ and $M$ non-projective (in particular, $M_1 \neq 0$). 
Let $PM$ be a projective cover of $M$. 
$$
{\beginpicture
    \setcoordinatesystem units <0.2cm,.2cm>
%%%%%%%%%%%%%%%%%%%%%%%%%%%%%%%%%%%%%%%%%%%%%%%%%%
\multiput{} at 0 0  20 13 /
\plot 0 4  4 0  5 1  6 0  10 4  14 0  15 1  16 0  20 4  10 14  0 4 /
\plot 5 9  10 4  15 9 / 
\plot -4 0  0 4 /
\plot 20 4  24 0  /
\setdots <1mm>
\plot -4 0  24 0 /
\multiput{$\ss\bullet$} at /
\put{$\Omega M$} at -1.9 4.4
\put{$\Omega(M/M_1)$} at 1 9.7
\put{$PM$} at 11.9 14.9
\put{$M$} at 16.5 9.5
\put{$M_1$} at 12 4 
\put{$M/M_1$} at 23 4.2 
\endpicture}
$$
Then $\pd M/M_1 = \pd \Omega(M/M_1)+1$
and there is an exact sequence $0 \to \Omega M \to \Omega(M/M_1) \to M_1 \to 0$.
Since both modules $\Omega M$ and $M_1$ are odd, 
$\pd \Omega(M/M_1) = \max\{\pd M_1,\pd \Omega M\}
=  \max\{\pd M_1,\pd M-1\}$, therefore $\pd M/M_1 = \pd \Omega(M/M_1)+1 =
\max\{\pd M_1+1,\pd M\}.$
	\smallskip 
Case 2: $M_2$ is a proper submodule of $M$. 
We first may assume that $S = \top M$
is the simple injective module (we just delete the simple modules $\tau^{-i}S$ with $i\ge 1$),
and consider the one-point extension $B = A[M/M_1]$, say with extension vertex 
$\omega$ (which now is the simple injective module), 
thus $\rad P(\omega) = M/M_1$, and therefore $\Omega(\omega) = M/M_1.$ 
Of course, $B$ is again a linear Nakayama algebra with $\tau(\omega) = S.$ Since $M/M_1$ has the even factor module $V = M_2/M_1$,
we see that $M/M_1$ is even, thus $\Omega(\omega) = M/M_1$ shows that $\omega$ is odd. 
   
$$
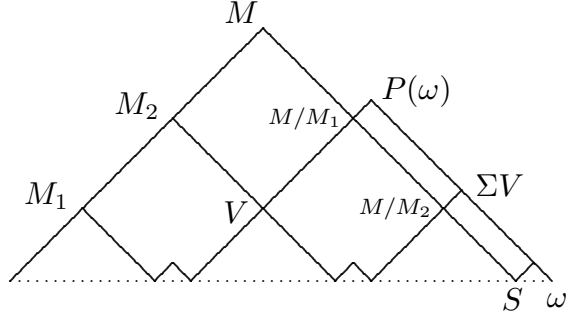
{\beginpicture
    \setcoordinatesystem units <0.2cm,.2cm>
%%%%%%%%%%%%%%%%%%%%%%%%%%%%%%%%%%%%%%%%%%%%%%%%%%
\multiput{} at 0 0  20 13 /
\plot -4 0  0 4  4 0  5 1  6 0  10 4  14 0  15 1  16 0  20 4  10 14  0 4 /
\plot 5 9  10 4  15 9 / 
\plot 15 9  16 10  26 0 /
\plot 25 1  24 0  20 4  21 5 /
\setdots <1mm>
\plot -4 0  26 0 /
\multiput{$\ss\bullet$} at /
\put{$M_1$} at -2 4.7
\put{$\ss M/M_1$} at 12.3 9
\put{$\ss M/M_2$} at 17.3 4
\put{$M_2$} at 3 9.7
\put{$M$} at 8.8 14.9
\put{$P(\omega)$} at 18.5 10.5
\put{$V$} at 8.5 3.8
\put{$\Sigma V$} at 23.1 5.3 
\put{$S$} at 23.8 -1 
\put{$\omega$\strut} at 26.3 -1 
\endpicture}
$$
The module $V$ is a submodule of the projective-injective module
$P(\omega)$ and $\Sigma V = P(\omega)/V$ is an extension
of $M/M_2$ by $\omega$ 
and therefore an odd module with
$$
 \pd \Sigma V = \max\{\pd M/M_2, \pd \omega\}.
$$
Since $P(\omega)$ is projective-injective, we have $\Omega\Sigma V = V$, 
thus 
$$
 \pd V = -1+\pd \Sigma V = \max\{-1+\pd M/M_2, -1+\pd \omega\}.
$$
Now, $-1+\pd \omega = \pd \Omega(\omega) = \pd M/M_1$. Thus, 
it remains to observe that
$$
  \pd M/M_1 = \max\{\pd M_1+1,\pd M\},
$$
but this we know already, due to case (1).
$\s$
	\medskip
{\bf Corollary.} 
{\it Let $A$ be a Nakayama algebra of finite global dimension.
Let $M$ be an indecomposable module and $X$ a non-zero 
subfactor of $M$. Then $\pd X$ is determined by $\char M$.}
$\s$
	\medskip
{\bf Remark.}
Let $X$ be a subfactor of the indecomposable module $M$. 
It is cumbersome, but not difficult, to write down an explicit formula
for $\pd X$ in terms of the coefficents of $\char M.$ On the other hand,
there is a very efficient algorithm for calculating $\pd X$ 
in terms of the coefficents of $\char M$, usind piles, see
4.8.
	\bigskip 
%%%%%%%%%%%%%%%%%%%%%%%%%%%%%%%%%%%%%%%%%%%%%%%%%%%%%%%%%%%%%%%%%%%%%%
{\bf 1.7. Remarks. (1)} 
{\it The sequence of numbers $\pd F_i(M)$ 
does {\bf not} determine $\pd M$, thus not $\char M$.}
Examples: Here are three modules $M$ of length $2$ (they are encircled), all
with $\pd F_1(M) = 3,\ \pd F_2(M) = 4.$ But we have $\pd M$ equal to $0,2,4,$
respectively:
$$
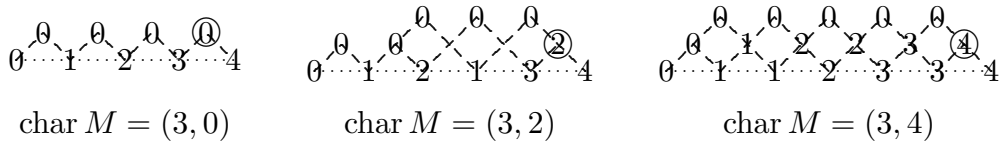
{\beginpicture
    \setcoordinatesystem units <.3cm,.3cm>
%%%%%%%%%%%%%%%%%%%%%%%%%%%%%%%%%%%%%%%%%%%%%%%%%%
\put{\beginpicture
\multiput{} at 0 0  8 1 /
\setdots <1mm>
\plot 0 0  8 0 /
\setdashes <1mm>
\plot 0 0  1 1  2 0  3 1  4 0  5 1  6 0  7 1  8 0 /
\multiput{$0$} at 0 0  1 1  3 1  5 1  7 1 /
\multiput{$1$} at 2 0 /
\multiput{$2$} at 4 0 /
\multiput{$3$} at 6 0 /
\multiput{$4$} at 8 0 /
\put{$\bigcirc$} at 7 1
\endpicture} at 0 0 
%%%%%%%%%%%%%%%%%%%%%%%%%%%%%%%%%%%%%%%%%%%%%%%%%%
\put{\beginpicture
\multiput{} at 0 0  10 2 /
\setdots <1mm>
\plot 0 0  10 0 /
\setdashes <1mm>
\plot 0 0  1 1  2 0  4 2  6 0  8 2  10 0 /
\plot 3 1  4 0  6 2  8 0  9 1 /
\multiput{$0$} at 0 0  1 1  3 1  4 2  6 2  8 2 /
\multiput{$1$} at 2 0  6 0 /
\multiput{$2$} at 4 0  9 1 /
\multiput{$3$} at 8 0 /
\multiput{$4$} at 10 0 /
\put{$\bigcirc$} at 9  1
\endpicture} at 12 0 
%%%%%%%%%%%%%%%%%%%%%%%%%%%%%%%%%%%%%%%%%%%%%%%%%%
\put{\beginpicture
\multiput{} at 0 0  12 2 /
\setdots <1mm>
\plot 0 0  12 0 /
\setdashes <1mm>
\plot 0 0  2 2  4 0  6 2  8 0  10 2  12 0 /
\plot 1 1  2 0  4 2  6 0  8 2  10 0  11 1    /
\multiput{$0$} at 0 0  1 1  2 2  4 2  6 2  8 2  10 2  /
\multiput{$1$} at 2 0  3 1  4 0 /
\multiput{$2$} at 5 1  6 0  7 1 /
\multiput{$3$} at 8 0  9 1 10 0 /
\multiput{$4$} at  11 1  12 0  /
\put{$\bigcirc$} at 11 1
\endpicture} at 26 0 
\put{$\char M = (3,0)$} at 0 -3
\put{$\char M = (3,2)$} at 12 -3
\put{$\char M = (3,4)$} at 26 -3
\endpicture}
$$

	\medskip 
{\bf (2)} {\it Also, $\char M$ is not determined by $\pd M$ and 
the sequence of the odd values $\pd F_i(M)$} 
(we have to know in addition the index $i$ with $\pd F_i(M)$ being
even):   
$$
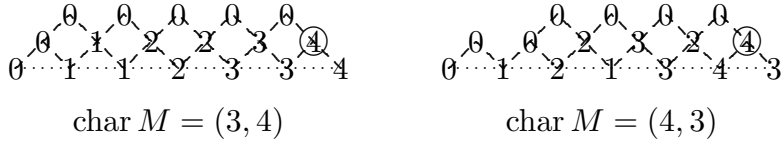
{\beginpicture
    \setcoordinatesystem units <.3cm,.3cm>
%%%%%%%%%%%%%%%%%%%%%%%%%%%%%%%%%%%%%%%%%%%%%%%%%%
%%%%%%%%%%%%%%%%%%%%%%%%%%%%%%%%%%%%%%%%%%%%%%%%%%
\put{\beginpicture
\multiput{} at 0 0  12 2 /
\setdots <1mm>
\plot 0 0  12 0 /
\setdashes <1mm>
\plot 0 0  2 2  4 0  6 2  8 0  10 2  12 0 /
\plot 1 1  2 0  4 2  6 0  8 2  10 0  11 1    /
\multiput{$0$} at 0 0  1 1  2 2  4 2  6 2  8 2  10 2  /
\multiput{$1$} at 2 0  3 1  4 0 /
\multiput{$2$} at 5 1  6 0  7 1 /
\multiput{$3$} at 8 0  9 1 10 0 /
\multiput{$4$} at  11 1  12 0  /
\put{$\bigcirc$} at 11 1
\endpicture} at 0 0 
%%%%%%%%%%%%%%%%%%%%%%%%%%%%%%%%%%%%%%%%%%%%%%%%%%

\put{\beginpicture
\multiput{} at 0 0  12 2 /
\setdots <1mm>
\plot 0 0  12 0 /
\setdashes <1mm>
\plot 0 0  1 1  2 0  4 2  6 0  8 2  10 0  11 1  /
\plot 3 1  4 0  6 2  8 0  10 2  12 0 /
\multiput{$0$} at 0 0  1 1  3 1  4 2  6 2  8 2  10 2 /
\multiput{$1$} at 2 0  6 0  /
\multiput{$2$} at 4 0  5 1   9 1 /
\multiput{$3$} at 7 1  8 0  12 0  /
\multiput{$4$} at 10 0  11 1 /
\put{$\bigcirc$} at 11 1
\endpicture} at 16 0 
%%%%%%%%%%%%%%%%%%%%%%%%%%%%%%%%%%%%%%%%%%
\put{$\char M = (3,4)$} at 0 -3
\put{$\char M = (4,3)$} at 16 -3
\endpicture}
$$
	\medskip
{\bf (3)} 
{\it The numbers $\pd\ \rad^i P$ do not determine $\char P$.} Here are two projective
modules $P,P'$ with $\char P \neq \char P',$ but $\pd \rad^i P = \pd\rad^i P'$ for all $i$.
$$
{\beginpicture
    \setcoordinatesystem units <0.3cm,.3cm>
%%%%%%%%%%%%%%%%%%%%%%%%%%%%%%%%%%%%%%%%%%%%%%%%%%
\put{\beginpicture
\multiput{} at -6 0  4 2 /
\setdots <.5mm>
\plot -6 0  -4 2  -2 0  0 2  2 0  3 1  /
\plot -5 1  -4 0  -2 2  0 0  2 2  4 0 /
\setdots <1mm>
\plot -6 0  4 0 /
\multiput{$\ss 0$} at -6 0  -5 1  -4 2  -2 2  0 2  2 2 /
\multiput{$\ss 1$} at -4 0  -3 1  -2 0 / 
\multiput{$\ss 2$} at -1 1  0 0  1 1 /
\multiput{$\ss 3$} at  2 0  3 1  4 0 /
\put{$\bigcirc$} at 2 2 
\endpicture} at 0 0
%%%%%%%%%%%%%%%%%%%%%%%%%%%%%%%%%%%%%%%%%%%%%%%%%%
\put{\beginpicture
\multiput{} at -4 0  4 2 /
\setdots <.5mm>
\plot -4 0  -3 1  -2 0  0 2  2 0  3 1 /
\plot -1 1  0 0  2 2  4 0 /
\setdots <1mm>
\plot -4 0  4 0 /
\multiput{$\ss 0$} at -4 0  -3 1  0 2  2 2 /
\multiput{$\ss 1$} at -2 0  -1 1  2 0 / 
\multiput{$\ss 2$} at 0 0  1 1 /
\multiput{$\ss 3$} at   3 1  4 0 /
\put{$\bigcirc$} at 2 2 
\endpicture} at 14 0
\endpicture}
$$
	\bigskip
{\bf 1.8. Warning.} We will show in 4.9 the following: If $A$ is a linear Nakayama algebra
and $M$ is an indecomposable module with $\char M = (z_1,\dots,z_m),$ then
$\pd X \le 1+\max z_i$ for all subfactors $X$ of $M$. Let us stress that the stronger
inequality $\pd X \le \max z_i$ holds true for odd modules, but 
is in general false, see for example the first two
examples in 1.7 (1).
	\bigskip\medskip
%%%%%%%%%%%%%%%%%%%%%%%%%%%%%%%%%%%%%%%%%%%%%%%%%%%%%%%%%%%%%%%%%%%%%
{\bf 2. The ascent algebra of a projective characteristic sequence.} 
	\medskip
Recall that a linear Nakayama algebra $A$ is said to be ascending provided for 
indecomposable projective modules  $P, P'$ with $\Hom(P,P') \neq 0$, we have
$|P| \le |P'|.$     
	\medskip
{\bf 2.1. The left rotation $\lambda$.} Assume that $m\ge 1.$ 
We define the {\it left rotation} $\lambda$ of projective 
characteristic sequences
as follows: Let $\bold z = (z_0,\dots,z_m)$
be a projective characteristic sequence of length $m+1$. Let
$$
 \lambda(z_0,\dots,z_m) = \left\{ 
   \matrix
        (0,z_0,\dots,z_{m-1}) &\text{if} & z_m = 0, \cr
        (z_0,\dots,z_{m-1}) &\text{if} & z_m = 1, \cr
        (z_m-2,z_1,\dots,z_{m-1}) &\text{if} & z_m \ge 2.\cr
   \endmatrix
   \right.
$$
	\medskip
{\bf 2.2. Proposition.} {\it Let $A, A'$ be ascending Nakayama algebras with
$\char P(\omega_A) = \char P(\omega_{A'})$. Then $A$ and $A'$ are isomorphic.}
	\medskip
Proof, by induction on the rank $n$ of $A$. If the rank of $A$ is $n=1$, then $A = k$
and $\omega_A = P(\omega_A) = k$, thus $\char P(\omega_A) = (0)$. If also $A'$
is a linear algebra with $\char P(\omega_{A'}) = (0),$ then $P(\omega_{A'})$ has length 1.
Thus $\omega_{A'}$ is projective. Since $A'$ is connected, we must have $A' = k$, thus
$A$ and $A'$ are isomorphic.

In general, in order to show that two linear Nakayama algebras $A,A'$ of rank $\ge 2$
are isomorphic, it is convenient to write 
$A = B[M],\ A'=B'[M']$ where $B,B'$ are linear Nakayama algebras, such that
$M$ is an indecomposable $B$-module with $\top M = \omega_B$, and 
$M'$ is an indecomposable $B'$-module with $\top M' = \omega_{B'}$. 
If we can show that $B,B'$ are isomorphic, and that $|M| = |M'|$, 
then $A, A'$ are isomorphic.

Thus, let us assume that $A,A'$ are ascending Nakayama algebras such that
$A$ has rank $n\ge 2$ and $\char P(\omega_A) = \char P(\omega_{A'}) = \bold z
= (z_0,\dots,z_m)$. In particular, the $A$-module $P(\omega_A)$ and 
the $A'$-module $P(\omega_{A'})$ both have length $m+1$, thus the $A$-module
$M = \rad P(\omega_A)$ and the $A'$-module $M' = \rad P(\omega_{A'})$
both have length $m$. 
The Serre subcategory in $\mod A$ generated by the 
simple modules $\tau^i\omega_A$ with $i\ge 1$ is of the form $\mod B$,
where $B$ is an ascending Nakayama algebra of rank $n-1.$ Of course, $M$ is
a $B$-module with $\top M = \omega_B$ and $A = B[M].$
Similarly, the Serre subcategory in $\mod A'$ generated by the 
simple modules $\tau^i\omega_{A'}$ with $i\ge 1$ is of the form $\mod B'$
for some ascending Nakayama algebra $B'$.
The $B'$-module $M'$ has $\top M' = \omega_{B'}$ and $A' = B'[M'].$

Thus, it remains to show that $B$ and $B'$ are isomorphic. 
We will show below that
$$
 \char P(\omega_B) = \lambda\, \bold z, \tag{$*$}
$$
and similarly, $\char P(\omega_{B'}) = \lambda\,\bold z.$
Since the rank of $B$ is $n-1$, and $\char P(\omega_B) = \char P(\omega_{B'})$
we know by induction that 
$B$ is isomorphic to $B'$, as we want to show. 
	\medskip
In order to verify $(*)$, 
three cases have to be distinguished.

First case: $\pd \omega_A = 1,$
thus $\bold z = (z_0,\dots,z_{m-1},1)$ and therefore
$\lambda\,\bold z = (z_0,\dots,z_{m-1})$. In this case $\rad P(\omega_A)
= P(\omega_B)$ and $\char P(\omega_B) = (z_0,\dots,z_{m-1}),$ thus
$\char P(\omega_B) = \lambda\,\bold z,$ as required.

If $\pd \omega_A \ge 2,$ 
then $\rad P(\omega_A)$ is a proper factor
module of $P(\omega_B).$ Since $A$ is ascending, we have
$|P(\omega_B)| \le |P(\omega_A)|,$ and therefore
$\rad P(\omega_A) = P(\omega_B)/\soc P(\omega_B).$ 

We consider now the second case: $\pd \omega_A \ge 2,$ and 
$\rad P(\omega_A)$ is odd.
In this case $\soc P(\omega_B)$ has to be even (since an indecomposable projective module
has an even composition factor), thus
$\char P(\omega_B) = (0,z_0,\dots,z_{m-1})$ and by definition
$\lambda\,\bold z = (0,z_0,\dots,z_{m-1}).$ Thus, $(*)$ holds true also in this case.

Third case: $\pd \omega_A \ge 2,$ and $\rad P(\omega_A)$ is even.
Since $\pd \omega_A \ge 2,$ we see that $\rad P(\omega_A)$ is not projective.
Since $\rad P(\omega_A)$ is even. $\omega_A$ is odd, thus $z_m = \pd \omega_A \ge 3$
and $\soc P(\omega_B) = \Omega^2\omega_A$ has projective dimension $z_m-2$.
It follows that $\char P(\omega_B) = (z_m-2,z_1,\dots,z_{m-1})$,
but this is just $\lambda\,\bold z.$
$\s$
	\medskip
{\bf 2.3. Proposition.} {\it 
For any projective characteristic sequence $\bold z = (z_0,\dots,z_m)$,
there is a (necessarily unique) ascending algebra $A = A(\bold z)$ with
$\char P(\omega_A) = \bold z.$} The algebra 
$A(\bold z)$ will be called the {\it ascent algebra} of $\bold z$.
	
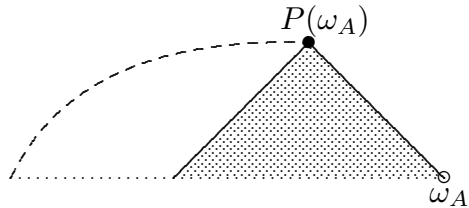
\medskip
$$
{\beginpicture
    \setcoordinatesystem units <.3cm,.3cm>
%%%%%%%%%%%%%%%%%%%%%%%%%%%%%%%%%%%%%%%%%%%%%%%%%%
\multiput{} at -1 0  15  4  /
\plot 0 0  5 5  10 0 /
\setsolid 
\multiput{$\bullet$} at   5 5  /
\setdots <1mm>
\plot -6 0  10 0 /
\setdashes <1mm>
\setquadratic
\plot -6 0  -2 3.9  5 5 /
\put{$P(\omega_A)$} at 5.6 5.8 
\multiput{$\circ$} at   10 0   /
 
\put{$\omega_A$} at 10.2 -.7  
\setlinear
\setshadegrid span <.4mm>
\vshade 0 0 0 <z,z,,> 5 0 5 <z,z,,> 10 0 0 /
\setdots <.5mm>
\endpicture}
$$
	\medskip
Proof of Proposition. If $\bold z = (z_0,\dots,z_m)$ is a projective 
characteristic sequence with $z_v = 0$, let $\epsilon\bold z = \sum_i z_i + v,$
this is a non-negative integer. We use induction on $\epsilon\bold z.$ 
The smallest possible value for $\epsilon\bold z$
is $0$, it occurs for $\bold z = (0)$, and $A = k$ satisfies $\char P(\omega_A) = (0).$
We assume now that  $\epsilon\bold z,$ thus $m\ge 1.$ 	\smallskip
Claim: {\it If $m\ge 1$, then
$\epsilon\lambda(z_0,\dots,z_m) < \epsilon(z_0,\dots,z_m).$}

Namely, if $z_m = 0,$ then $\epsilon\lambda(z_0,\dots,z_m) =  \sum_i z_i + 0 <
\sum_iz_i+m = \epsilon(z_0,\dots, z_m).$
If $z_m = 1,$ and $z_v$ is even, then 
$\epsilon(z_0,\dots,z_m) = \sum_i z_i + v$, whereas
$\epsilon(z_0,\dots,z_{m-1}) = \sum_{i=0}^{m-1}z_i + v = \sum_i z_i + v-1$.
Finally, if $z_m \ge 2,$ and $z_v$ is even, then the even entry of $\lambda(z_0,\dots,z_m)$ 
has index $v+1$, thus again 
$\epsilon\lambda(z_0,\dots,z_m) = 
\epsilon(z_m-2,z_0,\dots,z_{m-1}) = \sum_i z_i -2 +v+1 = \sum_i z_i + v - 1.$
This completes the proof of the claim. 
	\smallskip

By induction, there is the ascent algebra $A' = A(\lambda(z_0,\dots,z_m))$ with
last summit $P'$ and $\omega' = \top P' = \omega_{A'}$. We have $\char P' = 
\lambda(z_0,\dots,z_m).$
We are going to construct $A$ as a one-point extension of $B$, namely either
$A = B[P']$ or $A = B[P'/\soc P']$. We denote the extension vertex by $\omega = \omega_A$, thus
$P = P(\omega)$ has radical equal to $P'$ or to $P'/\soc P'$, respectively. 
In both cases, with $B$ also $A$ is ascending
(in case $A = B[P']$, we have $|P'| < |P|$; in case $A = B[P'/\soc P']$, 
we have $|P'| = |P|$). Note that for any $B$-module $X$, we have
$\pd{}_BX = \pd{}_AX$ and the indecomposable $A$-modules which are not $B$-modules
are factor modules of $P$. 

First, let $z_m = 1,$ thus $\lambda(z_0,\dots,z_m) = (z_0,\dots,z_{m-1}).$ 
In this case, let $A = A'[P']$.
Now $\rad P = P'$ shows that the composition factors of $P$ going upwards
are those of $P'$ followed by $\omega = P/P'$. Since $\pd \omega = 1$, we see that
$\char P = (z_0,\dots,z_{m-1},1)$, as we want to show.

Next, let $z_m \neq 1.$ Let $M = P'/\soc P'$ and $A = B[M].$ 
The composition factors of $P$ going upwards
are those of $M$ followed by $\omega = P/M$. 
In order to calculate $\char P$, we have to distinguish two cases. 

If $z_m = 0$, then the remaining entries $z_i$ with $0\le i \le m-1$ are odd
and $\char P' = (0,z_0,\dots,z_{m-1})$. 
Then $\char M = (z_0,\dots,z_{m-1})$. The composition factors of $P$ are
those of $M$ and in addition $\omega$. Since the composition factors of $M$ are
odd, it follows that $\omega$ has to be even, therefore $\char P = (z_0,\dots,z_{m-1},0)
= (z_0,\dots,z_m)$.

Now assume that $z_m \ge 2$ (and then $z_m$ has to be odd, since we deal with a
projective characteristic sequence).  
We have $\Omega^2\omega = \Omega M = \soc P'$, thus $\pd \omega = \pd (\soc P')+2.$
But $\pd(\soc P') = z_m-2$, thus $\pd \omega = z_m.$ 
Since $\pd(\soc P')$ is odd, $M$ has an even composition factor. Thus
$(z_0,\dots,z_m)$ is obtained from the sequence of numbers  
$\pd F_iM$, replacing the even number by $0$.
In order to determine $\char P$, we have to 
take the sequence of numbers $\pd F_iP$, and replace the even number by $0$. 
Since $\pd \omega$ is odd, we have to take the sequence of numbers $\pd F_iM$, replace
the even number by $0$ and add at the end $\pd \omega.$ 
It follows that $\char P = (z_0,\dots,z_m),$ as we want to show.$\s$
	\medskip 
{\bf Remark.} The proof shows:
	\medskip
{\it Let $A$ be a linear Nakayama algebra. 
If $S$ is a simple module, and $|P(\tau S)| \le |PS|,$ then}
$$
    \lambda \char PS = \char P(\tau S). 
$$
\vglue-.5cm
$\s$

{\bf 2.4. Corollary.} {\it If $A$ is an ascending Nakayama algebra, $A = A(\char(P(\omega_A)).$}
	\medskip
Proof. For any projective characteristic sequence $\bold z$, we have 
$\char P(\omega_{A(\bold z)}) = \bold z.$ 
Thus, 2.2 asserts that $A$ is isomorphic to $A(\char P(\omega_A)).$
$\s$
	\bigskip 

{\bf 2.5.} {\it For $h\ge 1,$ the maps
$$
  A \mapsto \char P(\omega_A).
$$
and $\bold z \mapsto A(\bold z)$ provide inverse bijections between
\item{$\bullet$} the ascending algebras $A$ of height $h$, and
\item{$\bullet$} the projective characteristic sequences $\bold z$ of length $h$.\par}
	\medskip
Proof. If we start with a characteristic sequence $\bold z$ and form $A = A(\bold z)$,
then, by construction $\char P(\omega_{A(\bold z)}) = \bold z.$ On the other hand,
Corollary 2.3 asserts that $A = A(\char P(\omega_A)).$
$\s$
	\medskip
Let us describe the cases $h = 1$ and $h = 2$ in detail:
	\smallskip
The only ascending algebra of height $1$ is $k$. Similarly, there is just one projective
characteristic sequence of length 1, namely $(0)$.

The ascending algebras of height $2$ are the radical-square-zero Nakayama 
algebras $A$ of type $\Bbb A_n$ with $n\ge 2.$ The module $P(\omega_A)$ has two composition factors, namely its socle with projective dimension $n-2$ and its top
with projective dimension $n-1$, thus $\char P(\omega_A) = (0,n-1)$ in case
$n$ is even, and $\char P(\omega_A) = (n-2,0)$ in case $n$ is odd.
	\bigskip\bigskip
%%%%%%%%%%%%%%%%%%%%%%%%%%%%%%%%%%%%%%%%%%%%%%%%%%%%%%%%%%%%%%%%%%%%%
{\bf 3. The partial $d$-closure of a linear Nakayama algebra.} 
	\medskip
{\bf 3.1.} Let $A$ be a linear Nakayama algebra and $d$ a positive natural number.
A simple module $S$ will be called {\it $d$-closed} provided $S$ is torsionless or
there is a module $M$ with socle $S$ and $\pd M \ge d.$

The algebra $A$ will be called {\it $d$-closed} provided all simple modules
are $d$-closed. Also, $A$ will be called {\it almost $d$-closed} 
provided any simple module which is not $d$-closed
is a composition factor of $P(\omega_A).$ And $A$ will be called 
{\it partially $d$-closed} provided all
composition factors of $P(\omega_A)$ are $d$-closed.
Thus, {\it $A$ is $d$-closed iff $A$ is both almost $d$-closed and partially $d$-closed.}
		\bigskip
{\bf 3.2.} Let $A$ be a linear Nakayama algebra. 
The algebra $B$ is said to be an
{\it extension of} $A$ provided $\mod A$ is a full subcategory of $\mod B$ which is 
closed under submodules, factor modules and 
projective covers. Note that in this case $\mod A$ is also closed under extensions, thus it
is a Serre subcategory of $\mod B$; also, 
the simple $B$-modules which are not $A$-modules are of the form $\tau^{-t}\omega_A$
with $t\ge 1.$
	
An extension $B$ of $A$ is said to be {\it descending} provided 
$$
 |P(\omega_A)| \ge |P(\tau_B^-\omega_A)| \ge |P(\tau_B^{-2}\omega_A)| \ge \cdots,
$$
thus provided $\pd \tau_B^{-i}\omega_A \neq 1$ for all $i\ge 1.$
	\bigskip
{\bf 3.3.} Let $A$ be a linear Nakayama algebra and $d$ a positive natural number.
The {\it $d$-cliff module} of $A$ is the module $Y_A = P(\omega_A)/U$,
where $U$ is the maximal submodule of $P(\omega_A)$ whose composition factors are
$d$-closed. 
Obviously, {\it the algebra $A$ is partially $d$-closed iff the $d$-cliff module of $A$ is zero.}

If the $d$-cliff module $Y_A$ is non-zero, let
$E_d(A) = A[Y_A]$ be the one-point extension using $Y_A$ and call it the
{\it $d$-cliff extension} of $A$. If $A$ is partially $d$-closed, then we write $E_d(A) = A$.
If we iterate this procedure, we get a sequence
$A, E_d(A), E_d^2(A),\dots$, with all algebras being descending extensions of $A$. 
	\medskip
{\bf Proposition.} {\it Let $A$ be a linear Nakayama algebra and $d$ a positive natural number.
There is $t \ge 0$ such that $E_d^t(A)$ is $d$-closed.}
If $E_d^t(A)$ is $d$-closed, we write 
$C_d(A) = E_d^t(A)$ and call it the {\it partial $d$-closure} of $A$.
	\medskip
Proof. Recall that for any indecomposable module $M$, $\Cal F(M)$ denotes the set of 
non-zero subfactors. Let $||M||_d$ be the sum of 
the values $\pd N$, with $N \in \Cal F(M)$ and $\pd N \le d.$

(a) {\it Assume that the $d$-cliff module $Y$ of $A$ is non-zero. 
Let $\omega_E$ be the extension vertex of 
$E = E_d(A)$, then}
$$
 ||P(\omega_E)/\soc P(\omega_E)||_d = ||Y||_d + |Y|.
$$
	\medskip
Proof. Let $P = P(\omega_E)$, thus $\rad P = Y$. Let $M = Y/\soc Y$ and $Z = P/\soc P.$ 
Let $\Cal Y$ be the set of non-zero submodules
of $Y$ and $\Cal Z$ the set of non-zero factor modules of $Z$.
There is the bijection $\Omega\:\Cal Z \to \Cal Y$ and we have $\pd \Omega N + 1 = \pd N$,
for $N \in \Cal Z$.
By assumption, all $Y\in \Cal Y$ have projective dimension at most $d-1$,
thus all $N\in \Cal Z$ have projective dimension at most $d$.
We see that $\Cal F(Y)$ is the disjoint union of the sets $\Cal Y$
and $\Cal F(M),$ and that $\Cal F(Z)$
is the disjoint union of the sets $\Cal Z$ and $\Cal F(M).$ 
$$
{\beginpicture
    \setcoordinatesystem units <.3cm,.3cm>
\multiput{} at 0 0  10 5 /
\setdots <1mm> \plot 0 0  10 0 /
\setsolid
\plot 0 0  5 5  10 0 /
\setdots <0.3mm>
\plot 1 1  2 0  6 4 /
\plot 2.5 5.5  8 0  9 1 /
\put{$\ss Y$} at 3 4 
\put{$\ss M$} at 5.9 3 
\put{$\ss P$} at 5.5 5.7
\put{$\ss Z$} at 7 4.1 
\multiput{$\bullet$} at 4 4  6 4 /
\put{$\Cal Y\strut$} at 1.3 2.5
\put{$\ss\Cal F(M)$} at 5 1 
\put{$\Cal Z\strut$} at 8.7 2.5
%\put{$\ss 1$} at 10 -.5
\setdashes <1mm>
\setquadratic
\plot -8 0  -4 4.7  2.5 5.5 /
\multiput{$\circ$} at 8 0  10 0  5 5  5 3 /
\put{$\ss \omega_A$} at 8 -.7
\put{$\ss \omega_E$} at 10 -.7
\endpicture}
$$
Since $\pd N = 1+ \pd \Omega  N$ for $N \in \Cal Z$, 
we see that
$$
 \sum\nolimits_{N\in \Cal Z} \pd N = 
 |Y| + \sum\nolimits_{X\in \Cal Y} \pd X,
$$
therefore $||Z||_d = |Y| + ||Y||_d.$
$\s$
	\medskip
(b) {\it If $Y$ is the $d$-cliff module, then $||Y||_d \le d|Y|(|Y|+1).$}
	\medskip
Proof. The set $\Cal F(Y)$ has cardinality $\frac12 |Y|(|Y|+1)$, and any summand
of $||Y||_d = \sum_{N} \pd N$ is bounded by $d$. $\s$
	\medskip
Consider now a sequence of $d$-cliff extensions
$A,\ E_d(A),\ \dots, \ E_d^t(A),$ say $E^i_d(A) = E_d^{i-1}(A)[Y_{i-1}]$ for
$1\le i \le t$, where $Y_{i-1} \neq 0$ is the $d$-cliff module of $E_d^{i-1}(A)$.

Let $\omega_0 = \omega_A$ and let $\omega_i = \tau_B^{-i}\omega_0$ be the extension vertex
for the extension $E_d^{i-1}(A) \subset E_d^i(A)$. Then we have $|P(\omega_i)| =
|Y_{i-1}|+1$ and $|P(\omega_0)| \ge |P(\omega_1)| \ge |P(\omega_2)| \ge \cdots.$
Assume that we have equalities
$$
 |P(\omega_0)| = |P(\omega_1)| = \cdots = |P(\omega_s)|, 
$$
thus $|Y_i| = |Y_0|$ for $1\le i < s.$ Applying (a) several times, we have
$$
 ||Y_{s-1}||_d = ||Y_0||_d + (s-1)|Y_0|.
$$
According to (b), we see that $s$ is bounded. It follows that there is $s > 0$
such that $|Y_0| > |Y_s|.$ By induction, it follows that $t$ has to be bounded.
Thus there is $t > 0$ such that the $d$-cliff module of $E_d^t(A)$ is zero. 
$\s$
	\bigskip
{\bf 3.4. Proposition.} {\it Let $A$ be a non-zero linear Nakayama algebra 
of global dimension at most $d$. Let $B$ be a descending extension of $A$. 
If $B$ is $d$-closed, then $B = C_d(A)$ and the global dimension of $B$ is $d$.}
	\medskip
Proof, by induction on the number $m$ of simple $B$-modules which are not $A$-modules.
If $m = 0$, then $A = B$ is $d$-closed, thus, by definition,  $C_d(A) = A$.

Assume now that $m\ge 1.$ Let $T = \tau_B^-\omega_A$ and $M = \rad P_B(T).$ We claim that
$M$ is the $d$-cliff module of $A$. 
First of all, $M$ is a factor module of $P(\omega_A),$ say $M = P(\omega_A)/U$, where $U$
is a submodule of $P(\omega_A).$ 
Since $B$ is a descending 
extension of $A$, we have $|P(\omega_A)| \ge |P_B(T)|,$ therefore $M$ is a proper
factor module of $P(\omega_A)$. In particular, no composition factor of $M$ is
torsionless (as an $A$-module).

We claim that the composition factors $S$ of $U$ are $d$-closed as $A$-modules. 
First of all, there is the socle $S = \soc U = \soc P(\omega_A)$. Of course,
this module is torsionless as an $A$-module, thus a $d$-closed simple $A$-module. 
Let $S$ be different from $\soc P(\omega_A).$ Then 
$S$ is not torsionless as a $B$-module. Since $S$ is $d$-closed as a $B$-module,
there is a $B$-module $M$ with $\soc M = S$ and $\pd_B M = d.$ Clearly, 
$M$ has to be an $A$-module, and $\pd_AM = \pd_BM = d$. This shows that $S$ is
$d$-closed as an $A$-module. 

On the other hand, $S = \soc M$ is not $d$-closed as an $A$-module. Namely, otherwise there
would exist an $A$-module $N$ with socle $S$ and $\pd N = d.$ But then $N$ is a submodule
of $P_B(T)$ and therefore $\pd_B P_B(T)/N = 1+\pd N = d+1.$ But this contradicts the assumption
that the global dimension of $B$ is at most $d$.

Altogether, we see that 
$U$ is the maximal submodule of $P(\omega_A)$ whose composition factors are $d$-closed.
This shows that $M = P(\omega_A)/U$ is the $d$-cliff module of $A$.
 
Let $\mod A'$ be the Serre subcategory of $\mod B$ generated by 
$\mod A$ and the simple module $T$.
Then $A' = A[M] = E_d(A)$ is a linear Nakayama algebra and $B$ is a deseending
extension of $A$. The number of simple $B$-modules which are not $A'$-modules is $m-1$.
Thus, by induction $B = C_d(A') = C_d(E_d(A)) = C_d(A).$ 
$\s$
	\bigskip
{\bf Corollary.} {\it Let $B$ be a concave Nakayama algebra which is $d$-closed and has
global dimension at most $d$.
Let $P$ be a summit of $B$. Then $B = C_dA(\char P).$}
	\medskip
Proof. Let $T = \top P$ and $\mod A$ the Serre subcategory of $\mod B$ generated by the
simple modules $\tau^iT$ with $i\ge 0.$ Since $P$ is a summit of $B$, the algebra $A$
is ascending and the algebra $B$ is a descending extension of $A$.
By construction, $P = P(\omega_A)$. 
According to 2.2, $A = A(\char P).$ Thus $B = C_d(A) = C_dA(\char P)).$
$\s$
	\bigskip
{\bf 3.5.} Let $A$ be a linear Nakayama algebra and $d$ a positive natural number.
Here are some properties of $C_d(A).$ 
	\smallskip
(1) {\it If $X$ is an indecomposable $C_d(A)$-module, but not an $A$-module, then
$\pd X \le d.$}
	\medskip 
(2) {\it If $C = C_d(A)$ is a proper extension of $A$, then $\pd \omega_C = d.$}
	\medskip
(3) {\it If $A$ is a non-zero algebra, then the global dimension of $C_d(A)$ 
is the maximum of $d$ and of the global dimension
of $A$.\par}
	\medskip
(4) {\it If $A$ is concave, also $C_d(A)$ is concave.}
	\medskip
Proof. (1), (2) and (4) follow immediately from the definitions involved. (3) is a direct consequence
of (1) and (2). $\s$
	\bigskip 
{\bf 3.6. Proposition.} {\it 
Let $A$ be a linear Nakayama algebra and $d$ a positive natural number.
If $A$ is almost $d$-closed, then $E_d(A)$ is almost $d$-closed.}
	\medskip
Proof. Let $A$ be almost $d$-closed and 
$Y$ the $d$-cliff module of $A$. Let $E = E_d(A) = A[Y]$ with extension vertex $\omega_E$.
Let $P = P(\omega_A)$ with socle $S$ and let $P' = P(\omega_E)$ with socle $S'$.
We have to show that all simple modules $S''$ which are predecessors of $S'$
are $d$-closed in $E$.  
The simple modules which are predecessors of $S$ are $d$-closed in $A$, thus in $E$, 
since $A$ is almost $d$-closed. If $S'' = S'$, then $S''$ is torsionless in $E$, since
it is the socle of $P'$. Thus, it remains to assume that $S''$ is a proper predecessor 
of $S'''$ and a proper successor of $S$. Then the maximality of $Y$ implies that
$IS''$ has a submodule of projective dimension at least $d$, thus $S''$ is
$d$-closed. 
$\s$
	\medskip
{\bf Corollary.} {\it Let $A$ be a linear Nakayama algebra and $d$ a positive natural number.
If $A$ is almost $d$-closed, then $C_d(A)$ is $d$-closed.} $\s$
	\bigskip
{\bf 3.7.} Finally, let us provide another characterization of the ascending algebras which
shows that ascending algebras are almost $d$-closed.
	\medskip
{\bf Lemma.} {\it A linear Nakayama algebra $A$ is ascending iff 
any simple module is torsionless or a composition factor of $P(\omega_A).$}
	\medskip
Proof. First, assume that $A$ is ascending. Let $S$ be a simple module which is not torsionless. 
Then $IS$ is not projective, thus $|PIS| > |IS|$. 
Assume that $S$ is also not a composition factor of $P(\omega_A).$ Then $\top IS \neq \omega_A$,
thus there is a simple module $T$ with $\tau T = \top IS.$ Since $IS$ is injective, 
we must have $|IS| \ge |PT|$. Therefore $|P(\tau T)| = |PIS| >|IS| \ge |PT|$, but this contradicts our
assumption that $A$ is  ascending.

Conversely, let $A$ be a linear Nakayama algebra and assume that $A$ is is not
ascending, thus there is a simple module $T$ with $|P(\tau T)| > |PT|$.
We show that there is a simple module $S$ 
which is not torsionless and not a composition factor of $P(\omega_A).$
Let $X$ be the factor module of $P(\tau T)$ with $|X| = |PT|$. Since 
$|P(\tau T)| > |PT|$, we see that $X$ is a proper non-zero factor module of 
an indecomposable projective module, thus not projective. On the other hand,
it follows from $|X| = |PT|,$ that $X$ is injective. Let $S = \soc X.$ Then $X = IS$.
Since $IS$ is not projective, $S$ is not torsionless. Since $\top IS = \tau T$
for some simple module $T$, we see that $\top IS \neq \omega_A$, and therefore $S$
is not a composition factor of $P(\omega_A).$ 
$\s$
	\medskip
{\bf Corollary.} {\it Let $A$ be an ascending Nakayama algebra and $d$ a positive natural number.
Then $A$ is almost $d$-closed, thus
$C_d(A)$ is $d$-closed.}
	\medskip
Proof. Let $A$ be an ascending Nakayama algebra. Then all the simple modules which are
not composition factors of $P(\omega_A)$ are torsionless, thus $d$-closed. Therefore
$A$ is almost $d$-closed. Corollary 3.5 shows that then $C_d(A)$ is $d$-closed.
$\s$
	\bigskip\bigskip

%%%%%%%%%%%%%%%%%%%%%%%%%%%%%%%%%%%%%%%%%%%%%%%%%%%%%%%%%%%%%%%%%% 
{\bf 4. The right rotation $\rho$ (and memory piles).}
	\medskip
In this section $m\ge 1$ will be a fixed positive integer.
	\medskip
{\bf 4.1. The rotation $\rho$.} Let $\Cal Z_m$ be the set of $m$-tupels of integers with at most one even entry; in particular, $\Cal Z_1 = \Bbb Z.$
Note that the elements of $\Cal Z_m$ with non-negative entries are just the characteristic
sequences.
	\medskip
We define $\rho = \rho_m\: \Cal Z_m \to \Cal Z_m$ (and call it the {\it right rotation} 
for $\Cal Z_m$) as follows:
for $\bold z = (z_1,\dots,z_m)\in \Cal Z_m$, let
$$
 \rho\,\bold z = \left\{
   \matrix 
     (z_2,\dots,z_m,z_1\!+\!1), & \text{if $z_2,\dots,z_m$ are odd,}  \cr
     (z_2,\dots,z_{v-1},z_1\!+\!1,
          z_{v+1},\dots,z_m,z_v\!+\!1), &\ \ \text{if $v>1$ and $z_v$ is even.} \cr
   \endmatrix
 \right.
$$
(For $m=1$, $\rho(z) = z+1$ is just the addition by $1$.) Note that 
if $z_1$ is even, then $\rho\,\bold z$ has only odd entries; otherwise, $\rho\,\bold z$
has an even entry (of course, just one).
	\medskip
Clearly, $\rho$ is bijective (thus invertible), and it is easy to write down the
corresponding formula for $\rho^{-1}.$ Let $\bold z = (z_1,\dots,z_m) \in \Cal Z_m.$
$$
 \rho^{-1}\bold z = \left\{
   \matrix 
     (z_m\!-\!1,z_1,\dots,z_{m-1}), &\quad \text{if $z_1,\dots,z_{m-1}$ are odd,}  \cr
     (z_v\!-\!1,z_1,\dots,z_{v-1},z_m\!+\!1,
          z_{v+1},\dots,z_m),   &\ \ \text{if  $v>1$ and $z_v$ is even.} \cr
   \endmatrix
 \right.
$$
	\medskip
{\bf Lemma.} {\it Let $\bold z \in \Cal Z_m$ be a characteristic sequence. 
Then $\rho^{-1}(\bold z)$ has a negative entry iff 
$\bold z$ is a projective characteristic sequence.} 
	\medskip
Proof. If $\bold z$ is projective, then $z_i = 0$ for some $i$, thus
the first entry of $\rho^-\bold z$ is $-1$ (and all other entries are non-negative).
If $\bold z$ is not projective, then all entries of $\bold z$ are positive, thus
all entries of $\rho^-\bold z$ are non-negative.
$\s$
	\bigskip
{\bf 4.2. Memory piles.}
A {\it pile} $\Gamma$ of height $h\ge 2$ is the Auslander-Reiten quiver of a linear
Nakayama algebra $A$ with Kupisch series the form $(1,2,\dots,h\!-\!1,h,\dots,h)$.
A pile has a unique projective vertex $R$ of length $h-1$, it is called its {\it radical}.
Similarly, it has a unique injective vertex $Y$ of length $h-1$, it is called its {\it cliff.}
If $\Gamma$ has $t$ summits, we may label the simple modules of $A$ by $S_1,\dots, S_{t+h-1},$
with $\tau S_i = S_{i-1}$ for $i\ge 2$; then $IS_{t+1}$ is the cliff module and 
the $A$-modules with socle $S_{t+1}$ form the {\it cliff ray}.

A {\it memory pile} is a pile with a memory function $\mu$ such that $\mu(x) = 0$
for all summits $x$ (actually, we only need to assume that $\mu(x) = 0$ for 
the first summit $x$).
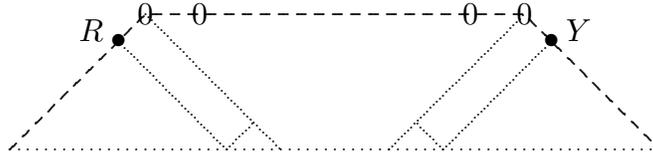
$$
{\beginpicture
    \setcoordinatesystem units <0.3cm,.3cm>
%%%%%%%%%%%%%%%%%%%%%%%%%%%%%%%%%%%%%%%%%%%%%%%%%%
\multiput{} at 0 0  24 5 /
\setdots <1mm>
\plot 0 0  24 0 /
\setdashes <1mm>  
\plot 0 0  5 5  19 5  24 0 /
\multiput{$0$} at  5 5  7 5  17 5  19 5 /
\put{$R$} at 3 4.4 
\put{$Y$} at 21 4.4 
\setdots <.5mm>
\plot 20 4  16 0  15 1 /
\plot 14 0  19 5 /
\plot 5 5  10 0 /
\plot 4 4  8 0  9 1 /
\multiput{$\bullet$} at 4 4  20 4 /
\endpicture}
$$
	\medskip
Serre subcategories which are piles play an 
important role in the sequel, in particular the summit pile and the descent piles mentioned in 6.5.
	\medskip
Dealing with a Serre subcategory of a module category 
may lead to confusion. Whenever it seems necessary, we add a corresponding
subscript. If $\Gamma$ is a pile, and $x$ is a non-projective vertex in $\Gamma$, we may
write $\Omega_\Gamma x$ for the first syzygy of $x$ in $\Gamma$, and $\pi_\Gamma x$ for the
Auslander-Reiten translate of $x$ in $\Gamma.$ Note that for a memory pile $(\Gamma,\mu)$, 
and $x$ a non-projective vertex of $\Gamma,$ 
we have $\mu(\Omega_\Gamma M) = -1+\mu(M).$
Also, if $(\Gamma,\mu)$ is a memory quiver, 
and $x$ is a vertex of $\Gamma$, we may write $\char_\mu x$ for 
the characteristic sequence of $x$ with respect to $\mu$.
	\bigskip
{\bf 4.3. Lemma (The module theoretic interpretation of $\rho$).}
{\it If $(\Gamma,\mu)$ is a memory pile with a single summit, with 
radical $R$ and cliff $Y,$ then
$\rho(\char R) = \char Y.$}
	
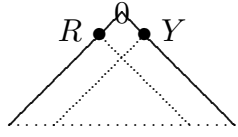
\medskip
$$
{\beginpicture
    \setcoordinatesystem units <0.25cm,.25cm>
%%%%%%%%%%%%%%%%%%%%%%%%%%%%%%%%%%%%%%%%%%%%%%%%%%
\put{\beginpicture
\multiput{} at 0 0  10 5 /
\setdots <1mm>
\plot 0 0  10 0 /
\setsolid
\plot 0 0  5 5  10 0 /
\setdots <.5mm>
\plot 2 0  6 4 /
\plot 4 4  8 0 /
\put{$R$} at 2.7 4.2
\put{$Y$} at 7.3 4.2
\multiput{$\bullet$} at 4 4  6 4 /
%\multiput{$\ss\blacksquare$} at   6 4 /
%\put{$\diamond$} at 0 0 
\multiput{$0$} at 5 5 /
\endpicture} at 0 0
\endpicture} 
$$

Proof. Let $\char R = (z_1,\dots,z_m).$
	\smallskip
First, let $z_2,\dots,z_m$ be odd. Then $\char Y = (z_2,\dots,z_m,x)$
for some $x$. 
If $z_1$ is odd, then $z_1 = \pd \soc R$. Now
$Y = IR/\soc R$, thus $\pd Y = 1+\pd \soc R = 1+z_1$. This shows that $Y$ is even,
therefore $\top Y$ has to be even, and $x = \pd Y = z_1+1.$
If $z_1$ is even, then $z_1 = \pd R$. Now $R = \Omega \top Y$, thus $\pd \top Y = 
1+ \pd \top R$ is odd, therefore $x = \pd \top R = 1+z_1.$
	\smallskip
Second, assume that $z_v$ is even, where $v > 1.$ 
Then $z_v = \pd R$. Since $\Omega \top Y = R$, we see that the last entry of 
$\char Y$ is $z_v+1$. We have $F_{i}(Y) = F_{i+1}(R)$ for $1\le i < m$.
In particular, $F_{v-1}(Y) = F_{v}(R)$ is even, therefore the entry of 
$\rho\,z $ with index $v-1$ is $\pd Y = \pd F_1(R)+1 = z_1+1$. The remaining
entries are obtained by the obvious index shift from entries of $\bold z$. 
$\s$ 
	\medskip
By induction on $t$, we obtain:
	\medskip
{\bf Corollary.} {\it Let $(\Gamma,\mu)$ be the memory pile with radical $R$
and $t$ summits. Let $Y$ be its cliff. Then $\char Y = \rho^t\char R.$}
$\s$
	\bigskip
%%%%%%%%%%%%%%%%%%%%%%%%%%%%%%%%%%%%%%%%%%%%%%%%%%%%%%%%%%%%%%%%%%
{\bf 4.4. Shift-Lemma.} For all characteristic sequences $\bold z$, we have
$$
 \rho^{m+1}(\bold z) = \bold z + (2,\dots,2).
$$
	\medskip
Proof. Ler $\Gamma$ be the pile of height $m+1$ with $m+1$ summits. Let $R$ be the radical
of $\Gamma$. We assume that $(\Gamma,\mu)$ is a memory pile with $\char_\mu R = \bold
z = (z_1,\dots,z_m)$. 

Let $Y$ be the cliff of $\Gamma$. Let $\char_\mu Y = (z'_1,\dots,z'_{h-1}).$
If $z_i$ is odd, then $z_i = \pd F_iR$ and $\pd F_iY = 2+\pd F_iR$ is also odd,
thus $z'_i = \pd F_iY = 2+z_i.$
If $z_i$ is even, then $\pd F_iR$ is even and $z_i = \pd R$.
In this case, also $\pd F_iY = 2+\pd F_iR$ is even and
$z'_i = \pd Y = 2+\pd R = 2+z_i.$ Altogether, we see that $z'_i = 2+z_i$
for all $i$, thus $\char Y = \char R + (2,2,\dots,2).$
$\s$
	\medskip
{\bf Corollary.} {\it Let $\bold z = (z_1,\dots,z_{h-1})$ and $t\ge 1$. Then 
$\rho^{th}(\bold z) = \bold z + t(2,\dots,2).$} $\s$
	\medskip 
Here is the pile $\Gamma$ of height $h$ with $th$ summits and 
with radical $R$. The black squares mark the vertices
$Y_i = \tau_\Gamma^{-ih}R$, the small circles are the vertices
$S_{1+ih} = \tau_\Gamma^{-ih}(\soc R).$

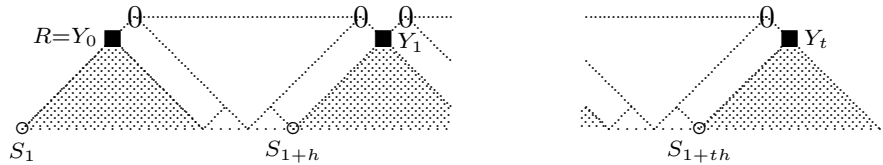
$$
{\beginpicture
    \setcoordinatesystem units <0.25cm,.25cm>
%%%%%%%%%%%%%%%%%%%%%%%%%%%%%%%%%%%%%%%%%%%%%%%%%%
\multiput{} at 0 0  35 5 /
\setdots <1mm>
\plot 0 0  19 0 /
\plot 26 0  38 0 /
\setsolid
\setdots <.5mm>
\plot 0 0  5 5  10 0  11 1  12 0  17 5  19 3 /
\plot 25 3  28 0  29 1  30 0  34 4  33 5 /
\plot 11 1  15 5  16 4 /
\plot 29 1  33 5 /
\plot 25 1  26 0  27 1 /

\plot 34 4  38 0 /
\plot 5 5   15 5 /
\plot 17 5  19 5 /
\plot 25 5  33 5 /
\multiput{$0$} at 5 5  15 5  17 5  33 5 /
\multiput{$\ss\blacksquare$} at   4 4  16 4 34 4 /
\put{$\ss R = Y_0$} at 1.9 4.1
\put{$\ss Y_1$} at 17.2 3.85
\put{$\ss Y_t$} at 35.2 4
\put{$\ss S_1$} at 0 -1
\put{$\ss S_{1+h}$} at 12 -1
\put{$\ss S_{1+th}$} at 30 -1
\multiput{$\circ$} at 0 0  12 0  30 0 /
\setdots <.5mm>
\plot 4 4  8 0  9 1 /
\plot 16 4  19 1 /
%\plot 34 4  35 5  40 0 /
%\plot 38 0  39 1 /
\setshadegrid span <.4mm>
\hshade 0 0 8 <z,z,,> 4 4 4   /
\hshade 0 12 19 <z,z,,>  1 13 19 <z,z,,>  4 16 16 /

\hshade 0 24.7 26 <z,z,,>  1 24.7 25 /
\hshade 0 30 38 <z,z,,> 4 34 34 /
\endpicture}
$$
	\bigskip
{\bf 4.5. The piles used in section 6.}
We will consider piles of height $h = m+1$ with radical $R$, cliff $Y$, and $s$ summits.
In (1) and (2), 
we start with $\char R$ and calculate $\char Y = \rho^s\char R$.
In (1$'$) and (2$'$), 
we start with $\char Y$ and calculate $\char R = \rho^{-s}\char Y$.
	\medskip
{\it Let $\Gamma$ be a pile with radical $R$, cliff $Y$ and $s$ summits.
	\smallskip
\item{\rm(1)} If $\char R = (0,c_2,\dots,c_m)$, $t \ge 0$ and $s = (m+1)t+1$, 
    then, of course, the numbers $c_i$ are odd. We have
    $\char Y = (c_2+2t,\dots,c_m+2t,1+2t)$.
	\smallskip
\item{\rm(1$'$)} Let $\char Y = (c_1',c_2',\dots,c_m')$, with all numbers $c_i'$ odd, let
    $0 \le t < \frac12 c'_i$ for all $c'_i$ and let $s = (m+1)t+1$.
    Then $\char R = (c_m'-2t-1,c_1'-2t,\dots,c_{m-1}'-2t)$.
	\smallskip
\item{\rm(2)} If $\char R = (c_1,\dots,c_m)$, $t\ge 0$ and $s = ht$, 
    then $\char Y = (c_1+2t,\dots,c_m+2t)$.
	\smallskip
\item{\rm(2$'$)}
    Let $\char Y = (c_1',\dots,c_m')$, with all numbers $c'_i$ being odd. Let
$0 \le t < \frac12 c'_i$ for all $i$ and let $s = (m+1)t$. Then $\char R = 
(c_m'-2t-1,c_1'-2t,\dots,c_m'-2t)$.\par}
	\smallskip
Proof: Assertion 
(2) is just the shift lemma 4.3. In order to establish (1), we first use the
definition of $\rho$ which asserts that 
$\rho(0,c_2,\dots,c_m) = (c_2,\dots,c_m,1)$. Then, according to (2), we
have $\rho^{ht}(c_2,\dots,c_m,1) = (c_2+2t,\dots,c_m+2t,1+2t)$. Altogether, we get
$\rho^{ht+1}(0,c_2,\dots,c_m) = (c_2+2t,\dots,c_m+2t,1+2t)$.
	\medskip
(2$'$) is again the shift lemma 4.3, but now formulated for $\rho^{-1}$. It remains
to prove (1$'$). According to (2$'$), $\rho^{-(m+1)t}(c_1',\dots,c_m') =
(c_1'-2t,\dots,c_m'-2t)$. We have to apply $\rho^{-1}$ once more. 
We have $\rho^{-1}(c_1'-2t,\dots,c_m'-2t) =
(c_m'-2t-1,c_1'-2t,\dots,c_{m-1}'-2t).$ Altogether we see that 
$\char R = (c_m'-2t-1,c_1'-2t,\dots,c_{m-1}'-2t)$.
$\s$
	\medskip
We will use these assertions in section 6 for a given concave algebra. 
The assertions (1) and (1$'$) concern the ``summit pile'', whereas
the assertions (2) and (2$'$) concern the ``descent piles'', see 6.5.

	\bigskip
{\bf 4.6. Lemma.} {\it Let $c_1,\dots,c_m$ be odd natural numbers. Then}
$$
  \rho^{m}(c_1,\dots,c_{m}) = (c_m+1,c_1+2,\dots,c_{m-1}+2).
$$ 
	\medskip
Proof. Let $R$ be the radical, $Y$ the cliff of $\Gamma$.
Let $S_i = \tau_\Gamma^{-i+1}\soc R$ for $1\le i \le 2m.$
Let $\rho^m(c_1,\dots,c_m) = (c'_1,\dots,c'_m)$. 
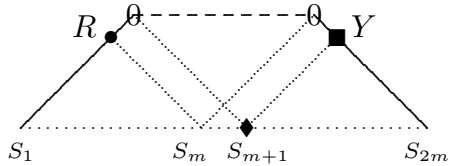
$$
{\beginpicture
    \setcoordinatesystem units <0.25cm,.25cm>
%%%%%%%%%%%%%%%%%%%%%%%%%%%%%%%%%%%%%%%%%%%%%%%%%%
\put{\beginpicture
\multiput{} at 0 0  18 5 /
\setdots <1mm>
\plot 0 0  18 0 /
\setsolid
\plot 0 0  5 5 /
\plot 13 5  18 0 /
\setdashes <1mm>
\plot 5 5  13 5 /
\setdots <.5mm>
\plot 4 4  8 0  13 5 /
\plot 5 5  10 0  14 4 /
\put{$R$} at 2.8 4.5
\put{$Y$} at 15.2 4.5
\multiput{$\bullet$} at 4 4  /
\multiput{$\ss\blacksquare$} at   14 4 /
\put{$\ss\blacklozenge$} at 10 0 
\multiput{$0$} at  5 5  13 5 /
\put{$\ss S_1$\strut} at 0 -1
\put{$\ss S_m$\strut} at 7.5 -1
\put{$\ss S_{m+1}$\strut} at 10.5 -1
\put{$\ss S_{2m}$\strut} at 18 -1
\endpicture} at 20  0
\endpicture} 
$$
We have $\Omega S_{m+1} = R$, and $\pd R = c.$  
Thus $S_{m+1}$ is even. It follows that $c'_1 = \pd Y.$
Since $\Omega Y = S_m$, we have $\pd Y = 1+c_m.$
According to 4.3, we have $\rho^{m+1}(c_1,\dots,c_m) =
(c_1+2,\dots,c_m+2)$, therefore $c'_i = c_{i-1}+2$ for 
$2\le i \le m.$
$\s$
	\bigskip
{\bf 4.7. Lemma.} {\it Let $m \ge 2.$ Let $(z_1,\dots,z_m)$ be a characteristic sequence. 
\item{\rm(a)} Let $z_m$ be even.
Then
$$
 \rho^{m-1}(z_1,\dots,z_m) = (z_{m-1}+1,z_{m}+1,z_1+2,\dots,z_{m-2}+2).
$$ 

\item{\rm(b)} Let $(z_1,\dots,z_m)$ be non-projective with $z_1$ even. Then}
$$
 \rho^{-m+1}(z_1,\dots,z_m) = (z_3-2,\dots,z_m-2,z_1-1,z_2-1).
$$ 
	\medskip
$$
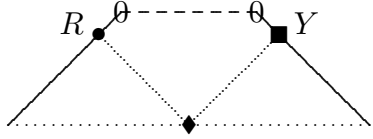
{\beginpicture
    \setcoordinatesystem units <0.25cm,.25cm>
%%%%%%%%%%%%%%%%%%%%%%%%%%%%%%%%%%%%%%%%%%%%%%%%%%
\put{\beginpicture
\multiput{} at 0 0  16 5 /
\setdots <1mm>
\plot 0 0  16 0 /
\setsolid
\plot 0 0  5 5 /
\plot 11 5  16 0 /
\setdashes <1mm>
\plot 5 5  11 5 /
\setdots <.5mm>
\plot 4 4  8 0  12 4 /
\put{$R$} at 2.8 4.5
\put{$Y$} at 13.2 4.5
\multiput{$\bullet$} at 4 4  /
\multiput{$\ss\blacksquare$} at   12 4 /
\put{$\ss\blacklozenge$} at 8 0 
\multiput{$0$} at 5 5  11 5 /
\endpicture} at 20  0
\endpicture} 
$$
	\medskip
Proof. We consider the pile $(\Gamma)$ with $m-1$ summits and radical $R$.
We denote the simple objects by $S_1,\dots,S_{2m-1}$ with $\tau S_i = S_{i-1}$
for $2\le i \le 2m-1.$ We have $\Omega S_{m+1} = R,$ $\Omega Y = S_{m-1}$ 
and $\Omega^2 S_{i+m-1} = S_i$ for $m+2 \le i \le 2m-1.$ 

Since $S_m$ is even, we have $\char R = (z_1,\dots,z_m) = (\pd S_1,\dots,\pd S_{m-1},\pd R)$
and $\char Y = (\pd Y,\pd S_{m+1},\dots,\pd S_{2m-1}) = (z'_1,\dots,z'_m).$ 

Since $\Omega S_{m+1} = R,$ we have $z'_2 = \pd S_{m+1} = \pd R + 1 = z_m+1.$
Since $\Omega Y = S_{m-1}$, we have $z'_1 = \pd Y = \pd S_{m-1}+1 = z_{m-1}+1.$
Finally, for $3 \le i \le m$, we have
$z'_i = \pd S_{i+m-1} = z_{i-2}+2,$ since $\Omega^2 S_{i+m-1} = S_{i-2}$.
This yields (a).

The characteristic sequences $(z'_1,\dots,z'_n)$ which we obtain in this way are
non-projective, and $z'_1 = z_m+1$ is even. Of course, since $\rho$ is intervible, 
$\rho^{-m+1}(z'_1,\dots,z'_m) = (z_1,\dots,z_m)$, and we have
$z_{i} = z'_{i+2}-2,$ for $1\le i \le m-2$ and $z_{m-1} = z'_1-1$, 
$z_m = z'_2-1.$ This yields the assertion (b).
$\s$
	\bigskip
{\bf 4.8. An algorithm for determining the projective dimension of subfactors of an indecomposable module.} Let $A$ be a linear Nakayama algebra and 
$M$ an indecomposable module with
$\char M = (z_1,\dots,z_m).$ We have seen in 1.6 that $\char M$ determines the projective
dimension of any subfactor of $M$, but we have refrained from writing down an
explicit formula. 

Here we outline an effective
algorithm in order to obtain the projective dimension of all subfactors of $M$.
There is no problem in case $M$ is odd, since one just uses the
maximum principle. Thus, we can assume that $M$ is {\bf even}. 
So, let $\bold z = (z_1,\dots,z_m)$ and assume that some $z_v$ is even. 
Let $\Gamma$ be the pile of height $m+1$ with $v$ 
summits. We want to determine a memory function $\mu$ on $\Gamma$,
with $\mu$ being obtained from $\bold z$
by applying $\rho, \rho^2,\cdots, \rho^v$. Let $Y$ be the cliff of $\Gamma$,
thus 
$$
  \char Y = \rho^v\bold z = 
(z_{v+1},\dots,z_m,z_v\!+\!1,z_1\!+\!2,\dots,z_{v-1}\!+\!2).
$$
In particular, {\it all coefficients of $\rho^v\bold z$ are odd,} 
thus the maximum principle yields immediately the values $\mu(x)$ for all 
subfactors $x$ of $Y$. We now use downward induction,
in order to calculate $\mu(x)$ for the vertices $x$ in the ray with index $i$, where
$1\le i \le v.$
	\medskip
{\bf Example.} Let $\bold z = (5,1,4,1)$. Here $v = 3$ (since $z_3 = 4$ is even). 
Thus, we have to deal with the pile of height 5 with 3 summits.
$$
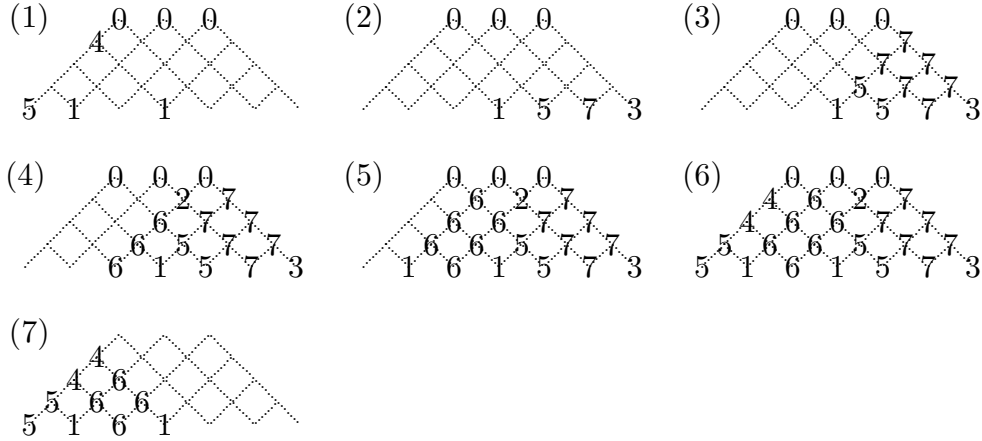
{\beginpicture
    \setcoordinatesystem units <0.25cm,.25cm>
%%%%%%%%%%%%%%%%%%%%%%%%%%%%%%%%%%%%%%%%%%%%%%%%%%
\put{\beginpicture
\multiput{} at 0 0  12 4 /
\put{(1)} at 0 4 
\setdots <.5mm>
\plot 0 0  4 4  8 0  10 2 /
\plot 1 1  2 0  6 4  10 0   11 1 /
\plot 2 2  4 0  8 4  12 0 /
\plot 3 3  6 0  9 3 /
\multiput{$0$} at 4 4  6 4  8 4  /
\multiput{$1$} at 2 0  6 0 /
\multiput{$3$} at /
\multiput{$4$} at 3 3 /
\multiput{$5$} at 0 0 /
\multiput{$6$} at /
\multiput{$7$} at /
\endpicture} at 0 0
%%%%%%%%%%%%%%%%%%%%%%%%%%%%%%%%%%%%%%%%%%%%%%%%%%
\put{\beginpicture
\multiput{} at 0 0  12 4 /
\put{(2)} at 0 4 
\setdots <.5mm>
\plot 0 0  4 4  8 0  10 2 /
\plot 1 1  2 0  6 4  10 0   11 1 /
\plot 2 2  4 0  8 4  12 0 /
\plot 3 3  6 0  9 3 /
\multiput{$0$} at 4 4  6 4  8 4  /
\multiput{$1$} at 6 0 /
\multiput{$3$} at 12 0 /
\multiput{$4$} at /
\multiput{$5$} at 8 0 /
\multiput{$6$} at /
\multiput{$7$} at 10 0 /
\endpicture} at 15 0
%%%%%%%%%%%%%%%%%%%%%%%%%%%%%%%%%%%%%%%%%%%%%%%%%%
\put{\beginpicture
\multiput{} at 0 0  12 4 /
\put{(3)} at 0 4 
\setdots <.5mm>
\plot 0 0  4 4  8 0  10 2 /
\plot 1 1  2 0  6 4  10 0   11 1 /
\plot 2 2  4 0  8 4  12 0 /
\plot 3 3  6 0  9 3 /
\multiput{$0$} at 4 4  6 4  8 4  /
\multiput{$1$} at 6 0 /
\multiput{$2$} at /
\multiput{$3$} at 12 0 /
\multiput{$4$} at /
\multiput{$5$} at 7 1  8 0 /
\multiput{$6$} at /
\multiput{$7$} at 8 2  9 1  9 3  10 0  10 2  11 1 /
\endpicture} at 30 0
%%%%%%%%%%%%%%%%%%%%%%%%%%%%%%%%%%%%%%%%%%%%%%%%%%
\put{\beginpicture
\multiput{} at 0 0  12 4 /
\put{(4)} at 0 4 
\setdots <.5mm>
\plot 0 0  4 4  8 0  10 2 /
\plot 1 1  2 0  6 4  10 0   11 1 /
\plot 2 2  4 0  8 4  12 0 /
\plot 3 3  6 0  9 3 /
\multiput{$0$} at 4 4  6 4  8 4  /
\multiput{$1$} at  6 0 /
\multiput{$2$} at 7 3 /
\multiput{$3$} at 12 0 /
\multiput{$4$} at  /
\multiput{$5$} at 7 1  8 0 /
\multiput{$6$} at 4 0  5 1  6 2 /
\multiput{$7$} at 8 2  9 1  9 3  10 0  10 2  11 1 /
\endpicture} at 0 -7
%%%%%%%%%%%%%%%%%%%%%%%%%%%%%%%%%%%%%%%%%%%%%%%%%%
\put{\beginpicture
\multiput{} at 0 0  12 4 /
\put{(5)} at 0 4 
\setdots <.5mm>
\plot 0 0  4 4  8 0  10 2 /
\plot 1 1  2 0  6 4  10 0   11 1 /
\plot 2 2  4 0  8 4  12 0 /
\plot 3 3  6 0  9 3 /
\multiput{$0$} at 4 4  6 4  8 4  /
\multiput{$1$} at 2 0  6 0 /
\multiput{$2$} at 7 3 /
\multiput{$3$} at 12 0 /
\multiput{$4$} at  /
\multiput{$5$} at 7 1  8 0 /
\multiput{$6$} at 3 1  4 0  4 2  5 1  5 3  6 2 /
\multiput{$7$} at 8 2  9 1  9 3  10 0  10 2  11 1 /
\endpicture} at 15 -7
%%%%%%%%%%%%%%%%%%%%%%%%%%%%%%%%%%%%%%%%%%%%%%%%%%
\put{\beginpicture
\multiput{} at 0 0  12 4 /
\put{(6)} at 0 4 
\setdots <.5mm>
\plot 0 0  4 4  8 0  10 2 /
\plot 1 1  2 0  6 4  10 0   11 1 /
\plot 2 2  4 0  8 4  12 0 /
\plot 3 3  6 0  9 3 /
\multiput{$0$} at 4 4  6 4  8 4  /
\multiput{$1$} at 2 0  6 0 /
\multiput{$2$} at 7 3 /
\multiput{$3$} at 12 0 /
\multiput{$4$} at 2 2  3 3 /
\multiput{$5$} at 0 0  1 1  7 1  8 0 /
\multiput{$6$} at 3 1  4 0  4 2  5 1  5 3  6 2 /
\multiput{$7$} at 8 2  9 1  9 3  10 0  10 2  11 1 /
\endpicture} at 30 -7
%%%%%%%%%%%%%%%%%%%%%%%%%%%%%%%%%%%%%%%%%%%%%%%%%%
\put{\beginpicture
\multiput{} at 0 0  12 4 /
\put{(7)} at 0 4 
\setdots <.5mm>
\plot 0 0  4 4  8 0  10 2 /
\plot 1 1  2 0  6 4  10 0   11 1 /
\plot 2 2  4 0  8 4  12 0 /
\plot 3 3  6 0  9 3 /
\multiput{$1$} at 2 0  6 0 /
\multiput{$4$} at 2 2  3 3 /
\multiput{$5$} at 0 0  1 1  /
\multiput{$6$} at 3 1  4 0  4 2  5 1   /
\endpicture} at 0 -14
\endpicture}
$$ 
We have labelled the pictures. Picture (1) shows the given data
(namely, the entries of $\bold z$ are the projective dimension of $M$ itself and of
3 of the 4 compositon factors of $M$, where $M$ is an indecomposable module with
$\char M = \bold z$),  whereas the last picture (7) shows the
projective dimension of all non-zero subfactors of $M$.
The formula for $\rho^3\bold z$ yields
the numbers exhibited in (2) (then we may forget the remaining numbers mentioned in (1);
they will systematically be recovered).
We use the maximum principle in order to obtain (3). Next, we calculate $\mu(x)$
for the vertices in the third ray, this yields (4). Then we deal with the second ray,
this yields (5), and finally with the first ray, this yields (6). 
The information which we were aiming at is presented in picture (7). 
	\medskip
How to get from picture (1) to picture (2)? We have mentioned that one may use
the formula for $\rho^v$. However, the following procedure seems to be more
efficient: Note that we present in (1) the characteristic 
of an {\bf even} module $M$ of length $m$. Let us label the vertices at the
lower boundary by $S_1, S_2,\dots$ going from left to right.
In (1) there is given $\pd M$ as well as $\pd S_i$ 
for $m-1$ composition factors $S_i = F_i(M)$ of $M$,
namely $\pd S_i$ with $1\le i < v$ as well as  $v < i \le m$.
Since $\Omega S_{m+1} = M$, we obtain $\pd S_{m+1} = 1+\pd M$ (in our example, this yields
the entry $5= 1+\pd M = \pd S_5$ in (2)). For $1\le i < v$, we obtain
$\pd S_{m+i}$ by the formula
$\pd S_{m+1+i} = 2+\pd S_i,$ since $\Omega^2 S_{m+1+i} = S_i$ (in our example, this yields
the numbers $7= 2+\pd S_1 = \pd S_6 $,  and $3= 2+\pd S_2 = \pd S_7 $).
In this way, we know $\pd S_i$ for all the indices $i$ with $v < i < m+1+v$,
thus for $m$ consecutive vertices at the lower boundary. All these numbers
$\pd S_i$ are odd, thus we can use the maximum principle
in order to get picture (3).
	\medskip 
{\bf Warning.} On a first sight, this procedure can be interpreted as being 
purely formal! 
On the other hand, the procedure also has a module
theoretic interpreation. But we should stress the following: In case we start with an 
$A$-module $M$ and consider ``a pile $\Gamma$ with radical $M$'', 
say with $s\ge 1$ summits, 
this pile $\Gamma$ usually cannot be 
realized by $A$-modules, but we have to use $B$-modules for some related linear Nakayama algebra $B$. Here is the recipe for obtaining such an algebra $B$. Let $T$ be the top of
$M$. Let $\mod A_0$ be the Serre subcategory of $\mod A$ generated by the
simple modules $\tau^tT$, with $t\ge 0$\,;
this yields a linear Nakayama algebra $A_0$ and we may consider $M = M_0$ as an
$A_0$-module (note that $M$, considered as an $A_0$-module, is injective). 
Second, we construct inductively 
$s$ one-point extensions as follows: we start with $A_1 = A_0[M_0]$,
then we form $A_2 = A_1[M_1]$, where $M_1 = \tau^-_{A_1}M_0,$ and so on, finally let
$B = A_s.$ Note that the operator $\Omega$ mentioned in the previous paragraph
is the syzygy functor $\Omega = \Omega_B$ in the category $\mod B.$
	\bigskip
{\bf 4.9. Proposition.} {\it Let $A$ be a linear Nakayama algebra and $M$ an indecomposable
module with $\char M = (z_1,\dots,z_m).$ If $X$ is a subfactor of $M$, then
$\pd X \le 1+\max z_i.$}
	\medskip 
Proof. If $M$ is odd, then we even have $\pd X \le \max z_i.$ Thus, let us assume that
$M$ is even. Say let $F_v(M)$ be even for some $1\le v \le m.$. Then we have
$$
 \rho^v(z_1,\dots,z_m) = (z_{v+1},\dots,z_m,z_v+1,z_{1}+2,\dots,z_{v-1}+2),
$$
and all these numbers are odd, thus $\pd Y \le 2+\max z_i$, where $Y$ is the cliff of the
pile with radical $R$ and $v$ summits. On the other hand, $F_v(M) = \Omega_\Gamma Y,$
therefore $\pd F_v(M) = -1+\pd Y \le 1+\max z_i.$

If $X$ is a subfactor of $M$, then the set of composition factors of $X$ is
a subset of the set of composition factors of $M$, and therefore
$\pd X \le \max_i F_i(M) \le 1+\max z_i.$
$\s$
	\bigskip\bigskip
%%%%%%%%%%%%%%%%%%%%%%%%%%%%%%%%%%%%%%%%%%%%%%%%%%%%%%%%%%%%%%%%%%%%%
{\bf 5.  Higher Auslander algebras.} 
	\medskip
{\bf 5.1. Proposition.} {\it Let $A$ be a connected Nakayama algebra of finite global
dimension, but not semisimple.  Let $d\ge 1.$ 
The following conditions are equivalent.
\item{\rm (i)} $A$ is a higher Auslander algebra of global dimension $d$.
\item{\rm (ii)} Any indecomposable module $M$ is torsionless or satisfies $\pd M = d.$
\item{\rm (iii)} Any indecomposable module $M$ is either torsionless or else satisfies 
    $\pd M = d$ (and not both).\par}
	\medskip
The proof will be given in 5.6. 
	\medskip
{\bf Corollary.} {\it A higher Auslander Nakayama 
algebra of global dimension $d$ is $d$-closed.}
	\bigskip

{\bf 5.2. Lemma.} {\it Let $A$ be an artin algebra with $\domdim A \ge d\ge 2$ and $S$
a simple module. Then $S$ is torsionless or $\Ext^i(S,A) = 0$ for $1\le i < d.$}
	\medskip
Proof. We calculate $\Ext^i(S,A)$ using a minimal injective coresolution of ${}_AA$,
say
$$
 0 \to {}_AA \to I_0 \to \cdots \to I_{d-1} \to I_d \to \cdots.
$$
Since the dominant dimension of $A$ is at least $d$, the modules $I_0,\dots,I_{d-1}$
are projective.
Now let $1\le i < d.$ The group
$\Ext^i(S,A)\neq 0$ is the homology of the complex $\Hom(S,I_\bullet)$
at the position $i$:
$$
   \Hom(S,I_{i-1}) \to \Hom(S,I_i) \to \Hom(S,I_{i+1}).
$$
Assume that $\Ext^i(S,A) \neq 0$. Then, we must have $\Hom(S,I_i) \neq 0.$ 
Since $1\le i < d$, the module $I_i$ is projective, thus $S$ is torsionless.
$\s$
	\medskip
{\bf 5.3. Lemma.} {\it Let $A$ be an artin algebra which is a higher Auslander algebra
of global dimension $d$.  If $P$
is indecomposable projective, $\id P \in\{0,d\}$. If $I$
is indecomposable injective, $\pd I \in\{0,d\}$.}
	\medskip
Proof. We show the first assertion, the second follows by duality. Let $P$ be indecomposable
projective. Assume that $1\le e = \id P < d$. We take
a minimal injective coresolution of $P$,
say
$$
 0 \to P \to I_0 \to \cdots \to I_{e-1} \to I_e \to 0.
$$
Since the dominant dimension of $P$ is $d$, all the modules $I_0,\dots,I_e$
are projective. In particular, $I_e$ is projective, thus the map $I_{e-1} \to I_e$ is
split epi. But this is impossible. 
$\s$
	\medskip
{\bf 5.4. Lemma (Madsen).} {\it Let $A$ be a connected Nakayama algebra. Let $S$ be a simple
module with $\pd S = \pd IS$. Then $A$ has finite global dimension and 
$\pd M = \pd S$ for all modules $M$ with $\soc M = S$.}
	\medskip
Proof. Let $M$ be a module with
$\soc M = S$, thus $IM = IS.$  
According to [6], Theorem 3.6, $\pd S$ is finite. Thus $A$
has finite global dimension (see for example [6], 3.7 (6)). 
In particular, $\pd M$ is finite, thus even or odd.
Now we use [6] 3.2 (a) and (b). If $\pd M$ is odd, then $\pd S \le \pd M$ and $\pd S$
is also odd. Thus $\pd IS$ is odd, and therefore $M\subseteq IS$ shows that 
$\pd M \le \pd IS$. 
Similarly, if $\pd M$ is even, then $\pd M \ge \pd IS$ and $\pd IS$ is even. Thus
$\pd S$ is even, and therefore $S \subseteq M$ implies that $\pd S \ge \pd M.$
$\s$
	\medskip
{\bf 5.5. Lemma.} {\it Let $A$ be an artin algebra of dominant dimension at least $1.$
Let $d\ge 1$. Assume that any module
of projective dimension smaller than $d$ is torsionless. Then the dominant
dimension of $A$ is at least $d$.}
	\medskip
Proof. If $d = 1$, nothing has to be shown. Thus, we may assume that $d\ge 2.$
Since the dominant dimension of $A$ is at least $1$, the injective envelope
of any projective module $P$ is projective, thus $\Sigma^1 P$ has projective dimension
at most $1$, thus $\Sigma^1P$ is torsionless and $\Sigma^2P$ has projective
dimension at most $2$. By induction, $\Sigma^{d-1}P$ has projective dimension at
most $d-1$, thus is torsionless, thus its injective envelope is projective
(a torsionless module $N$ is a submodule of a
projective module $Q$, thus the injective envelope of $N$ 
is a submodule, and therefore a direct summand,
of the injective envelope of $Q$).
In this way, we see that the first $d$ terms of a
minimal injective coresolution of $P$ are projective. 
$\s$
	\bigskip
{\bf 5.6. Proof of Proposition 5.1.} 

(i) implies (iii).
We assume that $A$ is a Nakayama algebra 
which is a higher Auslander algebra of global dimension $d$.
Let $M$ be an indecomposable module which is not torsionless. We have
to show that $\pd M = d.$ Since $M$ is not torsionless, it cannot be
projective, thus $1 \le \pd M \le d.$ 
If $d = 1$, then $\pd M = d$. Thus, we can assume that $d\ge 2.$

Let $S = \soc M$ and $I = IM = IS.$
Since $M$ is not torsionless, also $S$ is not torsionless.
Let $e = \pd S$. Since $S$ is not torsionless, it cannot be projective, 
thus $e\ge 1.$ 
Since the global dimension of $A$ is $d$, we have $e \le d.$ 
Now $\Omega^eS$ is projective and 
a minimal projective resolution of $S$ yields a non-zero element in
$\Ext^e(S,\Omega^eS).$ This shows that $\Ext^e(S,{}_AA) \neq 0.$ 
Lemma 5.2 implies that $e = d.$ 
According to 5.3, $\pd I = d$. Thus 5.4 asserts that $\pd M = d.$
	\smallskip
It remains to show
that an indecomposable module $M$ with $\pd M = d$ cannot be torsionless.
Thus, assume that $M$ is indecomposable, torsionless, and not projective. 
There is an embedding 
$u\:M \to P$ with $P$ projective such that the image of $u$ is contained in the
radical of $P$. But this implies that $P$ is a projective cover of $\Cok(u)$, thus
$\pd \Cok(u) = \pd M + 1 > d$, whereas the global dimension of $A$ is $d$. 
	\medskip
Obviously, (iii) implies (ii).
	\medskip
(ii) implies (i). Let us assume that every indecomposable module $M$ is torsionless or satisfies $\pd M = d.$ By assumption, the algebra $A$ has finite global dimension.
Let us show that the global dimension of $A$ cannot be greater than $d$. Namely, otherwise there exist an indecomposable module $M_0$ with $\pd M_0 = e > d.$ 
The condition (ii) asserts that $M_0$ is torsionless,
say a submodule of the projective module $P$. Then $M_1 = P/M_0$ 
$\pd M_1 = \pd M +1= e+1$. Inductively, we obtain modules $M_i$ with $\pd M_i = e+i$
for all $i\ge 0.$ Since by assumption $A$ has finite global dimension, we obtain
a contradiction. 

Not all modules are torsionless: otherwise the injective modules are torsionless,
thus $A$ is selfinjective. But a selfinjective algebra of finite global dimension is
semi-simple, contrary to our assumption. Thus, there exists a module with projective
dimension $d$, and therefore 
the global dimension of $A$ is equal to $d$.

Finally, we have to determine the dominant dimension of $A$. Since $A$ is a Nakayama algebra,
the dominant dimension of $A$ is at least $1$. According to 5.5, the dominant dimension 
of $A$ is at least $d$. Let us assume that the dominant dimension of $A$ is greater or
equal to $d+1$. Then there is an exact sequence
$$
 0 \to {}_AA \to I_0 \to I_1 \to \cdots \to I_d \to Z \to 0
$$
with all modules $I_i$ injectie and projective. Since we know already that the global
dimension of $A$ is $d$, we see that the image of the map $I_0\to I_1$ is projective. 
But then the embedding ${}_AA \to I_0$ splits and $A$ is self-injective. Since $A$ has
finite global dimension, it follows that $A$ is semisimple, a contradiction.
$\s$
	\medskip
{\bf 5.7. Remark. The ray property.} The modules which are not torsionless form rays:
{\it If $A$ is a higher Auslander algebra of global dimension $d$,
then the modules $M$ in a fixed ray are either all torsionless
or else have $\pd M = d.$}
	\bigskip\bigskip 
%%%%%%%%%%%%%%%%%%%%%%%%%%%%%%%%%%%%%%%%%%%%%%%%%%%%%%%%%%%%%%%%%%%%%%
{\bf 6. $\pd$-controlled  modules, $d$-piles and $d$-bound algebras.}
	\medskip
{\bf 6.1.}
An indecomposable module $M$ with $\pd Z \le \pd M$ for all subfactors
$Z$ of $M$ will be said to be {\it $\pd$-controlled .} 
According to the maximum principle for odd modules, {\it all odd modules
are $\pd$-controlled.} Even modules are usually not $\pd$-controlled: for example,
an indecomposable projective module $P$ is $\pd$-controlled only in case
$P$ is simple.  
In sections 7 and 8, we are going to work a lot with 
$\pd$-controlled  modules of a fixed projective dimension $d$. In case $d$ is even
(section 8), we will introduce the ``plus-strictly-increasing'' modules: these are special
$\pd$-controlled  modules, see Lemma 8.3(c).
	\medskip
{\bf Lemma.} {\it Let $M$ be an indecomposable module. If
there is a submodule $U$ with $\pd U \ge \pd Z$ for all subfactors $Z$
of $N$, then $M$ is $\pd$-controlled .}
	\medskip
Proof of Lemma. Since all odd modules are $\pd$-controlled , we can assume that $M$ is even.
Let $U$ be a submodule of $M$ with $d = \pd U \ge \pd Z$ for all 
subfactors $Z$ of $M$. 

We claim that $d$ is even. Assume for the contrary that $d$ is odd.
If $U \subset U' \subseteq M$ and $U'$ is odd, then there is the
inequality $\pd U \le \pd U'$ (since $U'$ is odd), but also the inequality
$\pd U \ge \pd U'$. Therefore $\pd U = \pd U'$.
Thus, we can assume that $U$ is not a proper submodule 
of an odd submodule of $M$. Therefore $U = \rad V$ for some even
submodule $V$ of $M$. Let $T = V/U$; this is a simple module
which has to be even. 
Consider the pile $\Gamma$ with radical $M$ and one summit. There is 
the vertex $U'$ such that $\Omega_\Gamma U' = U$, thus
$\pd U' = 1+d$. The socle of $U'$ is equal to $T$.
Since $T \subseteq U'$ are even modules, we have
$\pd T \ge \pd U'$. Therefore
$\pd T \ge \pd U' = 1+d > d$. 
On the other hand, $T$ is a subfactor of $M$, thus $\pd T \le \pd U = d$.
This contradiction shows that $U$ cannot be odd.

If $U = M$, nothing has to be shown. 
Thus, assume that $U$ is a proper submodule of $M$. There is a submodule $V$ of $M$ such that
$U = \rad V$. Again, consider the pile $\Gamma$
with radical $M$ and one summit. Let $U'$ be the vertex of $\Gamma$
with $\Omega_\Gamma U' = U$. Thus $\pd U' = d+1$ is odd.
Let $Z = \rad U'$. This is subfactor of $M$, and therefore
$\pd Z \le d.$ Since $\pd U' = d+1$ and $\pd\rad U' = d,$
the maximum principle asserts that $\pd U'/\rad U' = d+1.$
But we have $\Omega_\Gamma (\pd U'/\rad U') = M$. Thus, we see that
$\pd M = d.$ 
$\s$
	\medskip
{\bf 6.2.} 
{\bf Lemma.} {\it Let $(\Gamma,\mu)$ be a memory pile with cliff $Y$. Let $e$ be the
maximum of $\mu(Y')$, where $Y'$ is a subfactor of $Y$.
If $X$ is not a subfactor of $Y$, then $\mu(X) < e.$}
	\medskip
Proof. Assume that $X$ is not a subfactor of $Y$. If $X$ is
a summit, then $\mu(X) = 0$. Otherwise $X = \Omega X'$ for some
vertex $X'$ of the pile, which again is not a summit. By induction, we see that
$X = \Omega^iY'$ for a subfactor $Y'$ of $Y$ and $i \ge 1.$ 
Therefore $\mu(X) = \mu(Y')-i < \mu(Y') \le e.$ 
$\s$
	\medskip
{\bf Definition of a $d$-pile.} 
A memory pile $(\Gamma,\mu)$ is said to be a {\it $d$-pile} provided 
the cliff module $Y$ is $\pd$-controlled  and $\pd Y = d$ (thus, provided $d = \mu(Y) \ge
\mu(Z)$ for all subfactors $Z$ of $Y$, or, equivalently, provided
the maximum of the values $\mu(U)$, where $U$ is a submodule of the cliff $Y$, is $d$,
and second, $\mu(Z) \le d$ for all non-zero subfactors $Z$ of $Y$).
	\medskip
{\bf Corollary.} {\it Let $(\Gamma,\mu)$ be a $d$-pile, with cliff module $Y$.
If $X$ is a vertex of $\Gamma$ and not a subfactor of $Y$, then $\mu(X) < d.$}
$\s$
	\bigskip 
{\bf 6.3. The radical of a $d$-pile.} We are going to characterize the redicals
of the $d$-piles.
	\medskip
{\bf  Lemma.}
{\it Let $A$ be a linear Nakayama algebra and $R$ an indecomposable module.
\item{\rm(a)} The module $R$ is the radical of a $d$-pile iff
$\pd U < d$ for all submodules $U$ of $R$ and $\pd X \le d$ for all subfactors $X$ of $R$.
\item{\rm(b)} If $R$ is the radical of a $d$-pile $(\Gamma,\mu),$ then $(\Gamma,\mu)$
is uniquely determined by $\char R.$ \par}
	\medskip
Proof. First, let $(\Gamma,\mu)$ be a $d$-pile and $R$ its radical. Since $\mu(x) \le d$
for all vertices of $\Gamma,$ we have $\pd X \le d$ for all subfactors $X$ of $R$.
Also, if $U$ is a submodule of $R$, then there is $z\in \Gamma$ such that
$U = \Omega_\Gamma z$ and therefore $\pd U = -1+\mu(z) < d.$

Conversely, let us assume that $R$ is an indecomposable module such that
$\pd U < d$ for all submodules $U$ of $R$ and that $\pd X \le d$ for all subfactors $X$ of $R$.
We want to show that $R$ is the radical of a $d$-pile. We may assume that $\top R = \omega_A$.
We consider the one-point extension $A[R]$, say with extension vertex $S$. Then
$\rad PS = R$. Since all submodules of $R$ have projective dimension at most $d-1$,
the factor modules of $PS$ have projective dimension at most $d$. It follows that
all subfactors of $PS$ have projective dimension at most $d$. 

Let $\bold z = \char PS$ and $C = C_d(A(\bold z)).$ We have $|PS| = m+1.$
Let $s$ be minimal such that $|P(\tau^{-s}S)| \le m,$ and let $Q = P(\tau^{-s+1}S).$
Note that all the projective modules $P(\tau^{-i}S)$ with $0\le i < s$ have length $m+1$.
The Serre subcategory generated by these modules $P(\tau^{-i}S)$ with $0\le i < s$
is a pile $\Gamma$ with $s$ summits. 

We claim that $\Gamma$ is a $d$-pile. 
Since all subfactors of $PS$ have projective dimension at most $d$,
the definition of $E_d$ shows that the same is true for all subfactors of $P(\tau^{-i}S)$
with $i\ge 0.$ Also, let $Y$ be the cliff of $\Gamma$, thus $Y = Q/\soc Q$. 
Then there has
to be a submodule $U$ of $Y$ of projective dimension $d$, since $C$ is $d$-closed and
$Y$ is not torsionless. 

Finally, we have to show that $\Gamma$ is the unique $d$-pile with radical $R$. 
We show this by induction on the number of summits $s$ of $\Gamma$. If $s = 1,$
then $\Gamma$ is the pile with radical $R$ and one summit. Now assume that
$s \ge 2.$ Then $R$ determines uniquely the module $R' = \pi_\Gamma^{-}R$, and we have
$\char R' = \rho(\char R)$. Since $s\ge 2$, we see that $R'$ is the radical of a
$d$-pile $\Gamma'$ with $s-1$ summits, and by induction $\Gamma'$ is uniquely
determined by $R'$. 
$\s$
	\bigskip
{\bf 6.4. $d$-bound algebras.}
Let $A$ be a linear Nakayama algebra.
The algebra $A$ will be said to be {\it 
$d$-bound,} provided the global dimension of $A$ is at most $d$
and $A$ is $d$-closed. Any ray of a $d$-bound algebra contains either a 
projective module or else a module with projective dimension $d$. In particular, {\it if $A$
is $d$-bound, and $I$ is in indecomposable module which is injective and not projective,
then $I$ is $\pd$-controlled  with $\pd I = d.$}
	\bigskip
{\bf Proposition.} {\it Let $A$ be a $d$-bound linear Nakayama algebra. Let $S$ be a simple
module and $s\ge 1$ such that 
$$
 |PS| = |P(\tau^-S)| = \cdots = |P(\tau^{-s+1}S)| > |P(\tau^{-s}S|.
$$
Let $\Cal C$ be the Serre subcategegory of $\mod A$ generated by the modules
$P(\tau^{-i}S)$ with $0\le i < s$. Let $\Gamma$ be the Auslander-Reiten quiver of $\Cal C$.
Then $(\Gamma,\pd_A)$ is a $d$-pile.}
$$
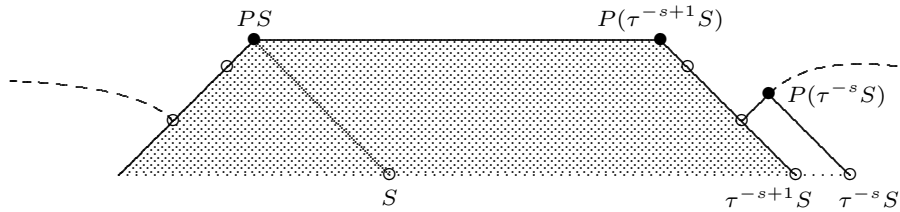
{\beginpicture
    \setcoordinatesystem units <.3cm,.3cm>
%%%%%%%%%%%%%%%%%%%%%%%%%%%%%%%%%%%%%%%%%%%%%%%%%%
\multiput{} at 0 0  27 5  /
\setdots <1mm>
\plot 0 0  27 0 /
\setdots <.3mm>
\plot 5 5  10 0 /
\setsolid
\plot 0 0  5 5  20 5  25 0 /
\plot 23 2  24 3  27 0 /

\setdashes <1mm>
\setquadratic
%\plot -2 5.5  2 5   4 4 /
\plot -4 3.5  0 3   2 2 /
\plot 24 3  26 4  29 4 /
\put{$\ss PS$} at 5 5.8
\put{$\ss S$} at 10 -.8
\put{$\ss P(\tau^{-s+1}S)$} at 20 5.8
\put{$\ss \tau^{-s+1}S$} at 24 -.8
\put{$\ss P(\tau^{-s}S)$} at 26.5 3
\put{$\ss \tau^{-s}S$} at 27.7 -.8
\multiput{$\bullet$} at 5 5  20 5  24 3 /
\multiput{$\circ$} at 2 2   4 4 10 0  25 0  27 0  21 4  23 2  /
\setlinear
\setshadegrid span <.4mm>
\hshade 0 0 25 <z,z,,> 5 5 20   /
\endpicture}
$$
	\medskip
Proof: Of course, $\Gamma$ is a pile. Since $A$ has global dimension at most $d$,
we have $\mu(x) \le d$ for all vertices $x$ of $\Gamma$. Let $Y$ be the cliff of $\Gamma$, thus
$Y = P(\tau^{-s+1}S)/\soc P(\tau^{-s+1}S).$ Then $Y$ is not torsionless. Since $A$ is
$d$-closed, there is a submodule $U$ of $Y$ with $\pd U = d.$
$\s$
	\bigskip
{\bf 6.5.} 
Dealing with a Nakayama algebra $A$, Serre subcategories which are piles play an 
important role. In the present paüer, we are mainly interested in
concave algebras. Such an algebra has a summit pile and usually several descent piles,
defined as follows. 
	\medskip
{\bf The summit pile.}
Let $A$ be concave with first summit $P$ and last summit $Q$.
The non-zero subfactors of the indecomposable injective modules 
which are successors of $P$ and predecessors of $Q$
form a pile, the {\it summit pile} of $A$.
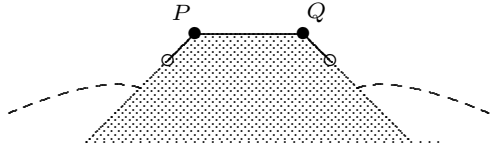
$$
{\beginpicture
    \setcoordinatesystem units <.3cm,.3cm>
\multiput{} at 0 0  12 4 /
\plot 3 3   4 4  8 4  9 3  /
\multiput{$\bullet$} at 4 4  8 4 /
\multiput{$\circ$} at 3 3  9 3 /
\put{$\ss P$} at 3.5 4.8 
\put{$\ss Q$} at 8.5 4.8
\setdots <.5mm>
\plot 0 0  4 4  8 4  12 0 /
\setdots <1mm>
\plot 0 0  13 0 /
\setdashes <1mm>
\setquadratic
%\plot 3 3  1 3  -2 2 /
%\plot 9 3  11 3  14 2 /
\plot 2 2   0 2  -3 1 /
\plot 10 2  12 2  15 1 /
\setlinear
\setshadegrid span <.4mm>
\hshade 0 0 12 <z,z,,> 4 4 8   /
\endpicture}
$$
Clearly, {\it the cliff module of the summit pile is the principle cliff module. }
	\medskip
{\bf The descent piles.} Let $A$ be concave. Let $0\le i \le h-2$. It is easy to see that
there is a unique indecomposable injective module $I_i$
of length $h-i-1$ which is not projective. The module $I_0$ is the principal
cliff module, and we have $I_{h-2} = \omega_A.$ 

Assume that $i\ge 1$ and that
$I_{i}$ is not a factor module of $I_{i-1}$. Then the
non-zero subfactors of the indecomposable injective modules 
which are successors of $I_{i-1}/\soc I_{i-1}$ and predecessors of $I_i$
 form a pile. 
These piles are called the {\it descend piles}
of $A$.
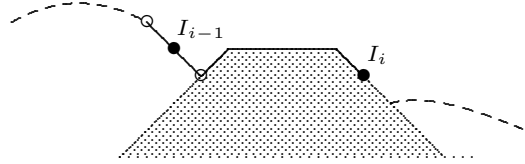
$$
{\beginpicture
    \setcoordinatesystem units <.3cm,.3cm>
\multiput{} at 0 0  13 6 /
\plot 1 5  3 3   4 4  8 4  9 3  /
\multiput{$\bullet$} at 2 4  9 3 /
\multiput{$\circ$} at 1 5  3 3 /
\put{$\ss I_{i-1}$} at 3 4.7 
\put{$\ss I_{i}$} at 9.5 3.8
\setdots <.5mm>
\plot 0 0  4 4  8 4  12 0 /
\setdots <1mm>
\plot 0 0  13 0 /
\setdashes <1mm>
\setquadratic
\plot -4 5  -1.5 5.7  1 5 /
%\plot 9 3  11 3  14 2 /
\plot 10 2  12 2  15 1 /
\setlinear
\setshadegrid span <.4mm>
\hshade 0 0 12 <z,z,,> 4 4 8   /
\endpicture}
$$
	\bigskip 
{\bf Corollary.} {\it  Let $A$ be a $d$-bound concave Nakayama algebra. 
Then the summit pile and the descent piles are $d$-piles.}
 $\s$
	\medskip
{\bf 6.6. Lemma.} {\it Let $A$ be a $d$-bound concave Nakayama algebra
and $R$ an indecomposable module. 
Let $\Gamma$ be a $d$-pile with radical $R$.
\item{\rm (a)} If $R = \rad P$, where $P$ is the first summit, then $\Gamma$ is
the summit pile.
\item{\rm (b)} If $R$ is not injective, and
$R = I/\soc I$, where $I$ is indecomposable injective and not projective,
then $\Gamma$ is a descent pile.\par}
	\medskip
Proof. (a) Let $R = \rad P$ with $P$ the first summit. The summit pile is a
$d$-pile with radical $R$. According to 6.3 (a), there is at most one $d$-pile with 
a given radical.

(b) Assume that $R$ is not injective, and 
$R = I/\soc I$, where $I$ is indecomposable injective and not projective.
The corresponding descent pile with radical $R$ is a $d$-pile. 
According to 6.3 (a), there is at most one $d$-pile with 
a given radical.
$\s$
	\bigskip
{\bf 6.7. Remark.} We have seen in 6.2: Given a pile $(\Gamma,\mu),$  
the maximum of $\mu$ will be obtained only on 
subfactors of the cliff. On the other hand, one should be aware that 
{\it any indecomposable non-projective module $M$ occurs as the cliff of a 
uniquely determined pile, namely of the summit pile of $A(\char M)$.}
	\bigskip\medskip
%%%%%%%%%%%%%%%%%%%%%%%%%%%%%%%%%%%%%%%%%%%%%%%%%%%%%%%%%%%%%%%%%% 
{\bf 7. The case $d$ odd.}
	\medskip
In this section, $d$ is always an odd natural number.
	\medskip
First, we present a general property of odd modules over Nakayama algebras which are 
higher Auslander algebras of odd global dimension.
	\medskip
Let $A$ be a Nakayama algebra. 
An indecomposable module $N$ will be said to be {\it decreasing} provided $N$ is odd 
and $\char N$ is a decreasing sequence of odd numbers.
An indecomposable  module $M$ will be said to be {\it plus-decreasing} 
provided $\rho\char M$ is a decreasing sequence of odd numbers.
	\medskip
{\bf 7.1. Proposition.} {\it Let $A$ be a Nakayama algebra which is a 
higher Auslander algebra with odd global dimension. Any odd indecomposable module 
is decreasing.}
	\medskip
Proof. We assume that  $A$ is a Nakayama algebra which is a higher Auslander algebra with 
odd global dimension $d$. 
Let us first show the following property: 
	\smallskip
(a) {\it If $M'$ is a non-zero submodule of an odd indecomposable module $M$, then $\pd M' = \pd M$.}
We use downward induction on $e = \pd M$. Since $A$ has global dimension $d$, we have $e\le d$.
Since $M$ is odd, we have $\pd M' = e' \le e.$ 
If $e = d$, then $M$ cannot be torsionless, since otherwise $IM$ is projective and
$\pd IM/M = \pd M + 1 > d.$ Thus, also $M'$ is not torsionless. According to 5.1, $\pd M' = d.$

Now assume that $e < d.$ Then, according to 5.1, $M$ is torsionless, thus $IM$ is projective.
There are non-split short exact sequences $0 \to M \to IM \to \Sigma M \to 0$ and
$0 \to M' \to IM \to \Sigma M' \to 0$. Since $IM$ is projective, $\pd \Sigma M = e+1$
and $\pd \Sigma M' = e'+1.$ Since $e'+1 \le e+1 < d,$ the modules $\Sigma M$ and $\Sigma M'$
both are again torsionless (and non-projective), thus 
there are non-split short exact sequences $0 \to \Sigma M \to I\Sigma M \to \Sigma^2 M \to 0$ and $0 \to \Sigma M' \to I\Sigma M' \to \Sigma^2 M' \to 0$. 
Since $I\Sigma M$ and $I\Sigma M'$ both are projective, $\pd \Sigma^2 M = e+2$ and
$\pd \Sigma^2 M' = e'+2$. Note that $\Sigma^2M'$ is a (non-zero) submodule of
$\Sigma^2 M$. By induction, $e'+2 = \pd \Sigma^2 M' = \pd \Sigma^2 M = e+2$ and therefore
$e' = e.$ This completes the proof of (a). 
	\smallskip
There is the following consequence: 
(b) {\it 
If $S$ is a simple module and both $S$ and $\tau S$ are odd, then $\pd \tau S \ge \pd S.$} 

Namely, If $S$ is simple and not projective, then
there exists an indecomposable module $M$ of length 2 with $\top M = S$. Now $\soc M =
\tau S$ and if both $S$ and $\tau S$ are odd, then the maximum principle shows that also
$M$ is odd, with $\pd M = \max\{\pd S,\pd \tau S\}.$ According to (a), $\pd M = \pd \tau S$,
thus $\pd \tau S = \pd M \ge \pd S.$
	\smallskip
The assertion mentioned in the Proposition is an immediate consequence: Let $M$
be an odd indecomposable module, say with $S = \top M.$
Let $\char M = (c_1,\dots,c_m)$. 
By the definition of a characteristic sequence, $c_i = \pd \tau^{m-i} S$, for
$1\le i \le m.$ Since $M$ is odd, the maxmimum principle asserts that all the
numbers $c_i$ are odd. 
According to (b), we have $\pd \tau^{m-i}S \ge \tau^{m-i-1}S$ for  $1\le i < m.$ 
But this just means that $c_i \ge c_{i+1}$ for $1\le i < m$.
$\s$
	\medskip
Let us stress that the property described in Proposition 7.1 has 
some interesting consequences (and actually is equivalent to these properties,
see our proof of 7.1), namely: (a) {\it If $M'$ is 
a non-zero submodul of an odd indecomposable module $M$, then $\pd M' = \pd M$.} And:
(b) {\it If $S$ is a simple module and both $S$ and $\tau S$ are odd, then $\pd \tau S \ge
\pd S.$}
	\medskip
{\bf Corollary.}  
{\it Let $d$ be odd. Let $A$ be a Nakayama algebra which is a 
higher Auslander algebra of global dimension $d$. 
Let $I$ be indecomposablle injective, not projective. Then $\pd I = d$ and $I$
is decreasing.}  
$\s$
	\bigskip
{\bf 7.2. Proposition.} {\it Let $d \ge c_1\ge c_2 \ge \cdots \ge c_u$ be odd
numbers. Then the algebra $C_dA(0,c_1,\dots,c_u,1)$ is a
concave higher Auslander algebra
of global dimension $d$.}  
	\medskip
Proof. Let $A = A(0,c_1,\dots,c_u,1)$ and 
$H = C_dA$. 
We show that the cliff modules of the summit pile and of the descent piles 
have projective dimension $d$ and decreasing characteristic sequences. 

We start with the summit pile.
Its radical $R$ has $\char R = (0,c_1,\dots,c_u)$ with decreasing odd numbers 
$c_1\ge c_2 \ge \cdots \ge c_u.$  
Thus, let $\Gamma$ be the pile with radical $R$, let $t = \frac12(d-c_1)$ 
and $s = (u+2)t$. Then, according to 4.5 (1), the cliff module $Y$ of $\Gamma$ has $\char Y = 
(c_1+2t,\dots,c_u+2t,1+2t)$, thus $\char Y$ is decreasing. 
Now, all the entries are odd, 
$c_1+2t = d$, and the remaining entries of $\char Y$
are bounded by $d$. Thus, we see that $Y$ is odd, thus $\pd$-controlled . Also, $\pd Y = d.$  
According to 6.6, $\Gamma$ is the summit pile of $H$.
In this way, we have shown that the cliff module of the summit pile has projective
dimension $d$ and a decreasing characteristic sequence.

Since $A$ is concave and of height $h$, there is a unique indecomposable
injective non-projective module $I_i$ of length $h-i$ for $1\le i \le h-1.$
We use induction on $i$ in order to show that $I_i$ has projective dimension $d$
and a decreasing characteristic sequence.
We have already seen that $I_1$ has projective dimension $d$ and that 
$\char I_1$ is a decreasing sequence. 
Now assume we know that $\pd I_i = d$ and that $\char I_i$ is decreasing
characteristic sequence, for some $1\le i < h-1.$ Let $\char I_i = (c_1,\dots,c_{h-i})$.
Of course, $c_1 = d.$ There are two possibilities.

First case: $I_{i+1}$ is a factor module of $I_i$, thus $I_{i+1} = I_i/\soc I_i.$
In this case $\char I_{i-1} = (c_2,\dots,c_{h-1})$ is a decreasing sequence and
$\pd I_{i-1} = c_2.$ Since $c_2$ is odd, $\pd U \le c_2$ for all submodules $U$ of
$I_{i-1}$. Since $A$ is $d$-closed and $I_{i-1}$ is not torsionless, we must have
$c_2 = d.$ 

Second, $I_{i+1}$ is not a factor module of $I_i$. Let $R_i = I_i/\soc I_i$.
Then there is the pile $\Gamma$ with radical $R_i$ and 
with $s$ summits, where $t = \frac12(d-c_1)$ and $s = ht$. According to 4.5 (2),
the cliff module of $\Gamma$ is $Y_i$ with $\char Y_i = (c_1+2t,\dots,c_{h-i}+2t)$.
Again, we use 6.6. It asserts that $\Gamma$ is the descent pile of $A$ with radical $R_i$, 
thus $Y_i = I_{i+1}$.
We see that $\char I_{i+1} = (c_1+2t,\dots,c_{h-i}+2t)$ is a decreasing sequence.
Also, $\pd I_{i+1} = c_1+2t = d.$
	\smallskip
Let $I$ be an indecomposable module which is injective and not projective.
As we have seen, $\pd I = d$ and $\char I$ is decreasing. Now, if $M$ is
an odd module with decreasig characteristic sequence, then $\pd U = \pd M$ for
any non-zero submodule $U$ of $M$. This shows that any non-zero submodule of $I$
has projective dimension $d$. According to 5.1, $H$ is a higher Auslander algebra
of global dimension $d$. 
$\s$
	\bigskip
{\bf 7.3. Lemma.} {\it Let $A$ be a Nakayama algebra. Let $M$ be an indecomposable module which is
plus-decreasing. Then $\soc M$ is even.  
If $U$ is a submodule of $M$, then  $\pd U \le \pd\soc M$.}
	\medskip
Proof. Let $\rho^-\char M = (c_1,\dots,c_m)$, then  
   $c_1 \ge \dots \ge c_m$, and all these numbers are odd. 
   Let $(\Gamma,\mu)$ be the pile with a unique
   summit and radical $M$.
   The maximum principle asserts that $\mu(\tau_\gamma^-M) = c_1$ 
   and that  $\mu(y) \le c_1$ for all vertices $y$ 
   which do not lie on the first ray, thus 
   $\mu(x) \le c_1-1$
   for all vertices $x$ on the first ray. Also, 
   $\mu(\tau_\gamma^-M) = c_1$ implies that 
   $\pd\soc M = c_1-1$. Altogether, we see: 
   If $U$ is a non-zero submodule of $M$,
   then $\pd U \le c_1-1 = \pd \soc M.$
$\s$
	\bigskip
{\bf 7.4. Proposition.}
{\it Let $d$ be odd and $m\ge 1$. Let $A$ be a concave Nakayama algebra of height $m+1$
which is $d$-bound. 
The following assertions are equivalent:
\item{\rm(1)} $A$ is a higher Auslander algebra of global dimension $d$.
\item{\rm(2)} $A$ has a  decreasing module $Y$ of length $m$
   with $\pd\top Y = 1.$
\item{\rm($2'$)} $A$ has a 
  decreasing module $Y'$ of length $m$.
\item{\rm($2''$)} $A$ has an injective 
decreasing module $Y''$ of length $m$.
\item{\rm(3)} $A$ has a projective
plus-decreasing module $R$ of length $m$.
\item{\rm($3'$)} $A$ has a 
  plus-decreasing module $R'$ of length $m$.
\item{\rm($3''$)} $A$ has a 
 plus-decreasing module $R''$ of length $m$ with $\pd\soc R'' = d-1.$
\item{\rm($4$)} $A = C_dA(0,c_1,\dots,c_{m-1},1)$ with odd numbers $c_i$ such that
  $d \ge c_1 \ge c_2 \ge \cdots \ge c_{m-1}.$\par}
	\medskip

For $A$ a concave Nakayama algebra of height $m+1$ which is a higher Auslander algebra 
of global dimension $d$, the assertions (2) to (3$''$) concern the indecomposabe modules $N$
of length $m$ such that $N$ or $\tau^- N$ is odd. Here are sketches which show the possible positions
of such modules $N = Y'$ of $N = R'$:
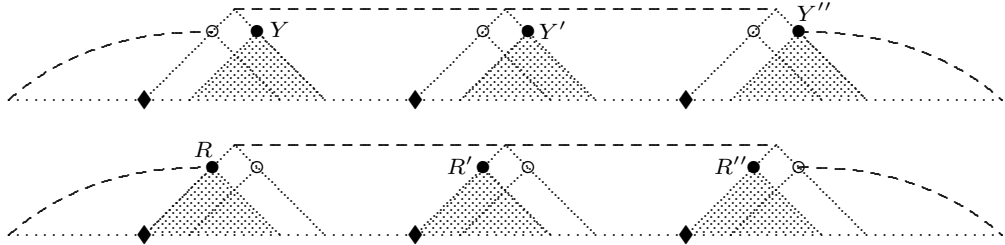
$$
{\beginpicture
    \setcoordinatesystem units <.25cm,.25cm>

%%%%%%%%%%%%%%%%%%%%%%%%%%%%%%%%%%%%%%%%%%%%%
\put{\beginpicture
\multiput{} at -6 0  36 5 /
\multiput{$\ss\blacklozenge$} at 0 0  12 0  24 0 /
\setdots <1mm>
\plot -6 0  38 0 /
\setdots <.5mm>
\plot 0 0  4 4  8 0 /
\plot 12 0  16 4  20 0 /
\plot 24 0  28 4  32 0 /

\plot 3 3  6 0 /
\plot 2 0  5 3 /

\plot 15 3  18 0 /
\plot 14 0  17 3 /

\plot 27 3  30 0 /
\plot 26 0  29 3 /

\multiput{$\bullet$} at 5 3  17 3  29 3 /
\multiput{$\circ$} at 3 3  15 3  27 3 /
\put{$\ss Y$} at 6 3.1
\put{$\ss Y'$} at 18.1 3.1
\put{$\ss Y''$} at 29.7 3.9

\setdashes <1mm>
\setquadratic
\plot -6 0  -2 2.3  3 3 /
\plot 29 3  34 2.3  38 0 /

\setlinear
\plot 4 4  28 4 /
\setdots <1mm>

\setshadegrid span <.4mm>
\hshade 0 2 8 <z,z,,> 3 5 5  /
\hshade 0 14 20 <z,z,,> 3 17 17  /
\hshade 0 26 32 <z,z,,> 3 29 29  /
\endpicture} at 0 0
%%%%%%%%%%%%%%%%%%%%%%%%%%%%%%%%%%%%%%%%%%%%%%%%
\put{\beginpicture
\multiput{} at -6 0  36 5 /
\multiput{$\ss\blacklozenge$} at 0 0  12 0  24 0 /
\setdots <1mm>
\plot -6 0  38 0 /
\setdots <.5mm>
\plot 0 0  4 4  8 0 /
\plot 12 0  16 4  20 0 /
\plot 24 0  28 4  32 0 /

\plot 3 3  6 0 /
\plot 2 0  5 3 /

\plot 15 3  18 0 /
\plot 14 0  17 3 /

\plot 27 3  30 0 /
\plot 26 0  29 3 /

\multiput{$\bullet$} at 3 3  15 3  27 3 /
\multiput{$\circ$} at 5 3  17 3  29 3 /
\put{$\ss R$} at 2.6 3.8
\put{$\ss R'$} at 14 3.1
\put{$\ss R''$} at 26 3.1

\setdashes <1mm>
\setquadratic
\plot -6 0  -2 2.3  3 3 /
\plot 29 3  34 2.3  38 0 /

\setlinear
\plot 4 4  28 4 /

\setshadegrid span <.4mm>
\hshade 0 0 6 <z,z,,> 3 3 3  /
\hshade 0 12 18 <z,z,,> 3 15 15  /
\hshade 0 24 30 <z,z,,> 3 27 27  /

\endpicture} at 0 -6

\endpicture}
$$

Proof of Proposition. 
We assume that $A$ is a concave Nakayama algebra of height $m+1$, that $A$ is $d$-closed and has
global dimension at most $d$. 
	\smallskip
(1) implies (2$''$). Here we assume in addition 
that $A$ is a higher Auslander algebra of global
dimension $d$. 
Let $Y''$ be the principal cliff module. Then $Y''$ is indecomposable injective and not
projective. Thus $Y''$ is decreasing by 7.1.
	\smallskip
(2$''$) implies (2$'$). Trivial.
	\smallskip
(2$'$) implies (2): 
Proof. Let $Y'$ be a decreasing module of length $m$.
Let $(c_1',\dots,c_m')$ be the characteristic of $Y'$. Then
$Y'$ is odd and $c_1'\ge c_2' \ge \dots \ge c_m'.$
Let $t = \frac12(c_m'-1)$, thus $t$ is an integer with $0 \le t < \frac12 c'_i$ for all
$c'_i$. Let $s = (m+1)t$. We consider the pile $\Gamma$ with $s$ summits and cliff $Y'$.
Let $Y$ be its radical. According to 7.3,
$\char Y = (c'_1-2t,\dots,c_{m}'-2t)$ and $c'_m-2t = 1$, by definition of $t$.
This shows that $Y$ is decreasing with $\pd \top Y = 1.$  
	\smallskip
(2) implies (3): Since $Y$ is odd, $Y$ is not projective. Let $R = \tau Y$.
Then $R$ is indecomposable, not injective, and $Y = \tau^-R$ is decreasing,
thus $R$ is plus-decreasing. Also, $R = \rad PY$. Since $\pd\top Y = 
\pd\top PY = 1$, we see that $R$ is projective. 
	\smallskip
(3) implies (3$'$). Trivial.
	\smallskip
(3$'$) implies (3$''$): Here we use again the shift lemma. 
   Let $R'$ be plus-decreasing of length $m$. Let $e= \pd \soc R'.$
   According to Lemma 7.3, all submodules $U$ of $R'$ have $\pd U \le e.$
   Let $t=\frac12(d-1-e)$ and $s = (m+1)t.$ 
   Let $(\Gamma,\mu)$ be the
   pile with $s = 2t$ summits and radical $R'$. Now $A$ is $d$-closed and has 
   global dimension at most $d$. Since the height of $A$ is $m+1$ and the length of
   $R'$ is $m$, we see that $\Gamma$ is part of the Auslander-reiten quiver of $A$.
   It follows that $R'' = \tau^{-s}R'$ is indecomposable with $\char R'' =
   \rho^{-s}R'.$ The shift lemma asserts that $R''$ is plus-decreasing and that
   $\pd\soc R'' = e+2s = d-1$.

	\smallskip
(3$''$) implies (2$''$): Let $R''$ be a plus-decreasing module of length $m$
   with $\pd\soc R'' = d-1$. 
Let $\rho\char R'' = (c_1,\dots,c_m)$, thus 
   $c_1 \ge \dots \ge c_m$ are odd numbers. 
Let $(\Gamma,\mu)$ be the pile with a unique
   summit and radical $R''$.
The maximum principle asserts
   that $\mu(\tau_\mu^-R'') = c_1,$ thus $c_1 = 1+\pd\soc R'' = d.$
   Since $A$ is $d$-closed, $R''$ is not injective. Since $A$ has height $m+1$ and $R''$
   has length $m$, the injective envelope of $R''$ has length $m+1$. As a consequence,
   $Y'' = \tau^-R''$ is the cliff of $\Gamma$ and $\char Y'' = \rho^-\char R'' = 
   (c_1,\dots,c_m)$ is decreasing and $\pd Y'' = d$. 
   Since the global dimension of $A$ is at most $d$, the module $Y''$ has to be injective.
	\smallskip
(2) implies (4). As in the proof that (2) implies (3), we consider the projective
cover $PY$ of $Y$. This is a summit, and since $R = \rad PY$ is projective, $PY$
is the first summit. Let $\mod B$ the predecessors
of $\top P.$ This is an ascending algebra, namely $B = A(P).$ 
Since $C$ is a concave higher Auslander algebra of global dimension $d$, 
it follows that $C = C_d(B)$. 
On the other hand, let $\char Y = (c_1,\dots,c_m).$ Since $Y$ is odd, we must have  
$\char P = (0,c_1,\dots,c_m)$. As we know, 
$d\ge c_1 \ge c_2\ge \cdots\ge c_m = 1.$
	\smallskip
(4) implies (1). This is 7.2.
$\s$
	\bigskip
{\bf 7.5. Proof of Theorems 1, 2, and 3.}
First, let $d\ge c_1\ge \cdots \ge c_u$ be odd numbers and $H = H_d(c_1,\dots,c_u) =
C_dA(0,c_1,\dots,c_u,1)$. Now $A(0,c_1,\dots,c_u,1)$ is ascending, thus
$C_dA(0,c_1,\dots,c_u,1)$ is concave. The equivalence of (4) and (1) in 7.4
yields Theorems 1 and 2. 

Let $P$ be the first summit of $H$ and $R = \rad P$. 
By construction, we have $\char P = (0,c_1,\dots,c_u,1)$ and $\char R = (0,c_1,\dots,c_u).$
If $u = 0$, then the Auslander-Reiten quiver of $H$ is just the pile of height $2$
with $d$ summits, and $\char Q = (0,d).$ Thus, we assume that $u\ge 1.$
Let $t= \frac12(d-c-1)$ and $s = (u+2)t$. 
The shift lemma asserts that $\tau^{-s}R$ has characteristic 
$(2t,c_1+2t,\dots,c_u+2t)$. It follows that
$Y = \tau^{-s-1}R$ has $\char Y = \rho^-\char Y = (c_1+2t,\dots,c_u+2t,2t+1)$.
Of course, $Y$ is the principal cliff module and its projective cover $Q$ 
(the last summit) has $\char Q = (0,c_1+2t,\dots,c_u+2t,2t+1).$
$\s$
	\bigskip
{\bf 7.6.} 
The characterizations presented in 7.4 refer to the existence of an idecomposable
module of length $h-1$, where $h$ is the height of the algebra. Actually, it
is sufficient to assume the existence of a corresponding sequence of simple modules:
	\medskip
{\bf Proposition.} {\it Let $A$ be a concave Nakayama algebra of height $h$
and let $d$ be an odd integer.
Then $A$ is a higher Auslander algebra with 
global dimension $d$ iff $A$ is $d$-closed and there is a simple module $S$
such that 
$d \ge \pd \tau^{h-2}S \ge \cdots \ge \pd \tau S \ge \pd S$
are odd numbers.}
	\medskip
Proof. 
First, assume that $A$ is a
higher Auslander algebras of odd global dimension $d$ and height $h$. 
According to 7.4, $A$ is $d$-closed and there 
is an indecomposable module $M$ with $\char M = (c_1,\dots,c_{h-1}),$
where  $d \ge c_1 \ge \cdots \ge c_{h-1}$
is a sequence of odd numbers. Let $S = \top M$. Then 
$\pd\tau^i S = c_{h-1-i}$ for $0 \le i \le h-2.$ 
These numbers are odd and bounded by $d$. 

Conversely, assume that $A$ is $d$-closed and has height $h$, with a 
simple module $S$ such that
$$
 d \ge \pd \tau^{h-1}S \ge \cdots \ge \pd \tau S \ge \pd S
$$ 
are odd numbers. According to 1.1 (4), we have $|PS| \ge h.$ Since $h$ is the height
of $A$, we have $|PS| = h.$ Let $Y = PS/\soc PS.$ Then $\char Y = 
(c_1,\dots,c_{h-1})$ with $d\ge c_1\ge \cdots \ge c_{h-1}$, 
thus 7.4 asserts that 
$A$ is a higher Auslander algebra of global dimension $d$.
$\s$
	\bigskip\bigskip

%%%%%%%%%%%%%%%%%%%%%%%%%%%%%%%%%%%%%%%%%%%%%%%%%%%%%%%%%%%%%%%%%%%%%%%%%%
{\bf 8. The case $d$ even.}
	\medskip
In this section, $d$ always will be even.
	\medskip 
{\bf 8.1.} {\bf Proposition.} {\it Let $A$ be a Nakayama algebras which is a
higher Auslander algebra with even global dimension.
\item{\rm(1)} There is no
chain $P_1 \subset P_2 \subset P_3$ of indecomposable projective modules.
\item{\rm(2)} There is no composition $I_1 \to I_2 \to I_3$ of proper epimorphisms between
indecomposable injective modules.\par}
	\medskip
Proof. We show the second assertion (2); the first one follows by duality.
We assume that $A$ is a higher Auslander algebra with global dimension $d$ being even and
that there is a composition $I_1 \to I_2 \to I_3$ of proper epimorphisms between
indecomposable injective modules. Then $I_2$ and $I_3$ are not projective. 
We may assume that the kernel $S$ of $I_2 \to I_3$ is simple. 
Then $\tau^-S$ is the socle of $I_2.$ 
Both modules $S$ and $\tau^-S$ are not torsionless, thus $\pd S = d  = \pd \tau^-S,$ 
according to 5.1. In particular, both module $S, \tau^-S$ are even. 
This is a contradiction to 1.2 (1).
$\s$
	\medskip 
{\bf 8.2.}
The indecomposable module $Y$ 
is called {\it plus-strictly-increasing} provided $\char Y = (e,c_2,\dots,c_m)$ with $e$ odd 
and $c_2 < c_3 < \dots < c_m < e,$ or, equivalently, provided $\rho \char Y$
has strictly increasing odd entries. 
The indecomposable module $R$ is called 
{\it minus-strict-increasing} provided $\char R = (c_1,\dots,c_{m-1},e)$ with $e$ even and
$e-1 < c_1 < \cdots < c_{m-1},$ er, equivalently, provided $\rho^{-1}\char R$ has odd strictly
increasing entries (since we allow that $e= 0$, the first entry of $\rho^{-1}\char R$ can be
negative). Note that a simple module is plus-strictly-increasing iff it is even iff it is
minus-strictly-increasing. 
	\medskip
{\bf 8.3. Proposition.} {\it Let $d$ be even. Let $A$ be a concave 
Nakayama algebra which is a higher Auslander algebra of global dimension $d$. 
Let $I$ be indecomposable injective, not projective and of length at least $2$.
Then $I$ is plus-strictly-fincreasing with $\pd I = d.$}
	\medskip
Proof. Note that $\soc I$ is not torsionless, therefore
5.1 asserts that $\pd \soc I = d$. Since $\soc I$ is even, $N = I/\soc I$ has to be odd. 

We use induction on $|I|$ in order to show that $\char N$ is strictly increasing. 
If $|I| = 2,$ nothing has to be shown. Thus, let $m = |I| \ge 3.$
Let $\char N = (c_2,\dots,c_m)$. Thus $\pd N = c = \max c_i.$

Let $\Gamma$ be the descent pile with radical $N$. Let $s$ be the number of summits of
$\Gamma$. 
Let us denote the simple modules
which belong to $\Gamma$ by $S_1,S_2,\dots,S_{s+m-1}$, going from left to right.
According to 5.1, we have $\pd S_i < d,$ for $1\le i \le s,$ and 
$\pd S_{s+1} = d$. 
The characteristic of $N$ asserts that $\pd S_i = c_{i+1}$ for $1\le i < m$
and we have $\pd S_m = c+1,$ since $\Omega S_m = N$. For $i>m$, 
we have $\Omega^2S_i = S_{i-m}$, thus $\pd S_i = 2+\pd S_{i-m}.$
It follows that $\pd S_i$ is even iff $m$ divides $i$. 
Since $\pd S_{s+1} = d$ is even, we must have $s+1 = (t+1)m$
for some $t\ge 0$. 
Let $Y$ be the cliff module of $\Gamma$, thus $\char Y = \rho^s(c_2,\dots,c_m)$
and $s = tm+m-1.$
According to the shift lemma, $\rho^{tm}\char N = \char N +
t(2,\dots,2) = (c'_2,\dots,c'_m)$. 
We use Lemma 4.6 and get 
$$
\align
  \char Y &=  \rho^s(c_2,\dots,c_m) = \rho^{m-1}(c'_2,\dots,c'_m) \cr
  &= (c'_m+1,c'_2+2,\cdots,c'_{m-1}+2) \cr
  &= (c_m+2t+1,c_2+2t+2,\dots,c_{m-1}+2t+2).
\endalign
$$
The first coefficient $c_m+2t+1$ is even, thus equal to $\pd Y = d$.
 Since the global
dimension of $A$ is $d$, all the coefficients are bounded by $d$, thus for $2\le i \le m,$
we have $c_i+2t+2 \le d = c_m+2t+1 < c_m+2t+2$, and therefore $c_i < c_m.$
	\smallskip

Since $\Gamma$ is a descend pile of $A$, the module $Y$ is injective and not projective.
Since $|Y| = |N|-1$, we know by induction that 
$Y/\soc Y$ is odd with strictly increasing characteristic sequence.
Now $\char (Y/\soc Y) = (c_2+2t+2,\dots,c_{m-1}+2t+2),$ thus 
$c_2 < \cdots < c_{m-1}.$
Since we know already that $c_{m-1} < c_m$, we see that $\char N$ is strictly
increasing.

Altogether, we see that $\char I = (d,c_2,\dots,c_m)$ with $c_2 < \cdots < c_m < d.$
$\s$
	\bigskip
{\bf 8.4. Lemma.} {\it Let $A$ be a Nakayama algebra of height $m+1$. 
\item{\rm(a)}
If $Y$ is an indecomposable module of length $m$ which is plus-strictly-increasing, then
$\tau^{m-1}Y$ is minus-strictly-increasing. 
\item{\rm(b)}
If $R$ is an indecomposable module of length $m$ which is
minus-strictly-increasing, then $\tau^{-m+1}R$ is plus-strictly-increasing or zero.
\item{\rm(c)} Let $Y$ be an indecomposable module which is plus-strictly-increasing. If $U$ is a non-zero
submodule of $Y$, then $\pd U = \pd Y$. If $Z$ is a subfactor of $Y$, then $\pd Z \le \pd Y.$
In particular, $Y$ is $\pd$-controlled .
\par}
	\medskip
The relationship between the modules $R$ and $Y$ considered in (a) and (b) is seen in
the following picture (it is important that $\top R = \soc Y$ is even):
$$
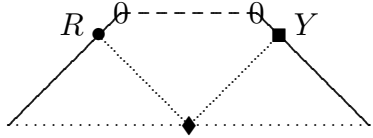
{\beginpicture
    \setcoordinatesystem units <0.25cm,.25cm>
%%%%%%%%%%%%%%%%%%%%%%%%%%%%%%%%%%%%%%%%%%%%%%%%%%
\put{\beginpicture
\multiput{} at 0 0  16 5 /
\setdots <1mm>
\plot 0 0  16 0 /
\setsolid
\plot 0 0  5 5 /
\plot 11 5  16 0 /
\setdashes <1mm>
\plot 5 5  11 5 /
\setdots <.5mm>
\plot 4 4  8 0  12 4 /
\put{$R$} at 2.8 4.5
\put{$Y$} at 13.2 4.5
\multiput{$\bullet$} at 4 4  /
\multiput{$\sssize\blacksquare$} at   12 4 /
\put{$\ss\blacklozenge$} at 8 0 
\multiput{$0$} at 5 5  11 5 /
\endpicture} at 20  0
\endpicture} 
$$

Proof of Lemma. For the proof of (a) and (b), we use Lemma 4.7.

(a) We assume that $Y$ is an indecomposable module of length $m$ which is
plus-strictly-increasing. Thus, $\char Y = (e,c_2,\dots,c_m)$ with $e > 0$ even
and $c_2 < c_3 < \dots < c_m < e$. Then 4.7 (b) asserts that
$\rho^{-m+1}\char Y = (c_3-2,\dots,c_m-2,e-1,c_2-1)$. Since $c_2$ is odd, $c_2-1$ is
even, and we have $c_2-2 < c_3-1,\dots,c_m-2 < e-1.$ It remains to be seen that
$\char \tau^{m-1}Y = \rho^{-m+1}\char Y.$ But, by assumption, $\pd Y = e > 0$, thus $Y$
itself is not projective, and for $1\le i \le m-2$ we have $\pd \tau^iY = c_{m-i}-1 > c_1 > 0$,
therefore also $\tau^iY$ is not projective.

(b) Assume that $R$ is an indecomposable module of length $m$ which is  
minus-strictly-increasing with $\char R = (c_1,\dots,c_{m-1},e).$ 
According to 4.7 (a), $\rho^{-m+1}R = (c_{m-1}+1,e+1,c_1+2,\dots,c_{m-2}+2)$. 
Since $c_{m-1}$ is odd, $c_{m-1}+1$ is even and we have 
$e+1 <c_1+2 < \cdots <c_{m-2}+2) < e+1.$ 
If $\tau^{-m+1}R$ is non-zero, then $\char \tau^{m+1}R = \rho^{-m+1}R.$  
	\smallskip
(c) Again, let $Y$ be plus-strictly-increasing. Let $\char Y = (e,c_2,\dots,c_m),$ 
thus $\rho\char Y = (c_2,\dots,c_m,e+1)$. Let $(\Gamma,\mu)$ be the memory pile
with radical $Y$ and a unique summit, let $Y'$ be the cliff of $\Gamma$. 
Since $e+1 \ge c_i$ for $2\le i \le m,$ and
these are odd numbers, the maximum principle asserts that $\mu(z) = e+1$ for all
all non-zero factor modules $z$ of $Y'$, thus $\mu(U)= e$ for all non-zero submodules
$U$ of $Y.$ Also, if $Z$ is a subfactor of $Y$, and not a submodule of $Y,$ then
$Z$ is a subfactor of $Y/\soc Y$. Since $Y/\soc Y$ is odd with projective dimension
$c_m$, it follows that $\pd Z \le c_m < e.$
$\s$
	\bigskip
{\bf 8.5. Lemma.} {\it Let $d$ be even and $u\ge 1.$
Let $c_1 < c_2 < \cdots < c_u$ be odd
numbers bounded by $d$. Let $A = C_dA(c_1,\dots,c_u,0,1)$. 
\item{\rm(a)} Let $t= \frac12(d-c_u-1).$ Then $A$ has $(u+2)t+u$ summits. 
   If $Y$ is the principal cliff of $A$, then 
    $\char Y = (d,1+2t,c_1+2t+2,\dots,c_{u-1}+2t+2).$
\item{\rm(b)} Let $I$ be
 an indecomposable module which is
injective and not projective. Then $I$ is plus-strictly-increasing with $\pd I = d$.
\item{\rm(c)} $A$ is a concave higher Auslander algebra of global dimension $d$.\par}
	\medskip
Proof. Let $A = C_dA(c_1,\dots,c_u,0,1)$.

(a) Let $P$ be the first summit, and $R = \rad P$.
Thus $\char R = (c_1,\dots,c_u,0).$  
Let $t= \frac12(d-c_u-1)$ and $s = (u+2)t+u$. Let $(\Gamma,\mu)$ be the pile with
radical $R$ and $s$ summits. Note that $R$ is minus-strictly-increasing, thus we can use
4.7 (see also 8.3). According to 4.7 (a), we have 
$x = \rho^u\char R = (c_u+1,1,c_1+2,\dots,c_u+2)$, thus $x$ is
plus-strictly-increasing with $\mu(x) = c_u+1.$ 
The shift lemma 4.4 shows that
$\rho^s\char R = \rho^{(u+2)t}x = (c_u+1+2t,1+2t,c_1+2+2t,\dots,c_u+2+2t).$
By definition of $t$, we have $c_u+1+2t = d$. 
Thus, $y = \rho^s(\char R)$ is plus-strictly-increasing and $\mu(y) = d.$ This shows
that $\Gamma$ is a $d$-pile. According to 6.6, 
we have constructed the summit pile of $A$, and $y$ is the pricipal cliff
of $A$.
	\smallskip
(b) Let $I_i$ be the indecomposable injective module of length $u+1-i$, where $0 \le i \le u$
(these are all the indecomposable injective modules which are not projective). We show
by induction on $i$ that $I_i$ is plus-strictly-increasing and $\pd I_i = d.$

The module $I_0$ is just the principl cliff-module, thus (a) shows that $I_0$ is
plus-strictly-increasing and that $\pd I_0 = d.$

Now, let $1 \le i \le u.$ Let $R_i = I_{i-1}/\soc I_{i-1}$.
We want to construct the descent pile with radical $R_i$ and cliff module $I_i.$
Since $I_{i-1}$ is plus-strictly-increasing, we have $\char R_i = (z_1,\dots,z_m)$
with odd numbers $z_1 < \cdots < z_m$ (here, $m = |R_i| = u+1-i,$ but this is
not relevant). Of course, all $z_i < d$; let $t = \frac12(d-z_m).$
We construct the pile $(\Gamma,\mu)$ with radical $R_i$ and $s = (m+1)t+m$
summits. 
According to 4.6, we have $x = \rho^m(z_1,\dots,z_m) =
(z_m+1,z_1+2,\dots,z_{m-1}+2).$ Note that $x$ is plus-strictly-increasing and $\mu(x) = z_m+1$.
According to the shift lemma 4.4, 
$y = \rho^s (\char R) = \rho^{(m+1)t+m}(\char R) = \rho^{m+1}t x$ is also
plus-strictly-increasing, and $\mu(y) = z_m+1+2t = d.$ 
This shows
that $\Gamma$ is a $d$-pile. According to 6.6, 
we have constructed the descent pile of $A$ with radical $R_i$, thus with 
cliff module $I_i = y$. As we have shown, $I_i$ is plus-strictly-increasing and $\pd I_i = d.$ 
	\smallskip
(c) 
If $M$ is indecomposable and not torsionless, then $IM$ is indecomposable, injective and
non-projective, thus equal to $I_i$ for some $0\le i \le u.$ According to (b), $I_i$
is plus-strictly-increasing with $\pd I_i = d$. According to Lemma 8.3 (c), we have $\pd M = d$, since
$M$ is a non-zero submodule of $IM$.
According to 5.1, $A$ is a higher Auslander algebra of global dimension $d$.
$\s$
	\bigskip
{\bf 8.6. Proposition.}
{\it Let $d$ be even and $m\ge 2$. 
Let $A$ be a concave Nakayama algebra of height $m+1$
which is $d$-bound. 
	\smallskip
The following assertions are equivalent:
\item{\rm(1)} $A$ is a higher Auslander algebra of global dimension $d$.
\item{\rm(2)} $A$ has a plus-strictly-increasing module $Y$ of length $m$
    with $\pd (\soc_2Y/\soc Y) = 1.$
\item{\rm($2'$)} $A$ has a 
  plus-strictly-increasing module $Y'$ of length $m$.
\item{\rm($2''$)} $A$ has a plus-strictly-increasing module $Y''$ of length $m$ is injective.
\item{\rm(3)} $A$ has a 
minus-strictly-increasing module $R$ of length $m$ which is projective.
\item{\rm($3'$)} $A$ has a 
  minus-strictly-increasing module $R'$ of length $m$.
\item{\rm($3''$)} $A$ has a 
 minus-strictly-increasing module $R''$ of length $m$ with $\pd\top R'' = d.$

\item{\rm($4$)} $A = C_dA(c_1,\dots,c_{m-1},0,1)$ with odd numbers $c_i$ such that
  $c_1 < c_2 < \cdots < c_{m-1} < d.$\par}
	\bigskip
\noindent

{\it There is also the following condition:
	\smallskip
\item{\rm(5)} $A$ has an indecomposable module $N$ of length $m$ with $\char N$
being strictly increasing.
	\smallskip
\noindent
It implies the equivalent conditions {\rm(1)} to {\rm(4)}, but the converse is not true.
\par}
	\bigskip
For $A$ a concave Nakayama algebra of height $m+1$ which is a higher Auslander algebra 
of global dimension $d$, the assertions (2) to (3$''$) concern the indecomposabe modules $N$
of length $m$ which are plus-strictly-increasing or minus-strictly-increasing. 
Here are sketches which show the possible positions
of such modules $N = Y'$ of $N = R'$:
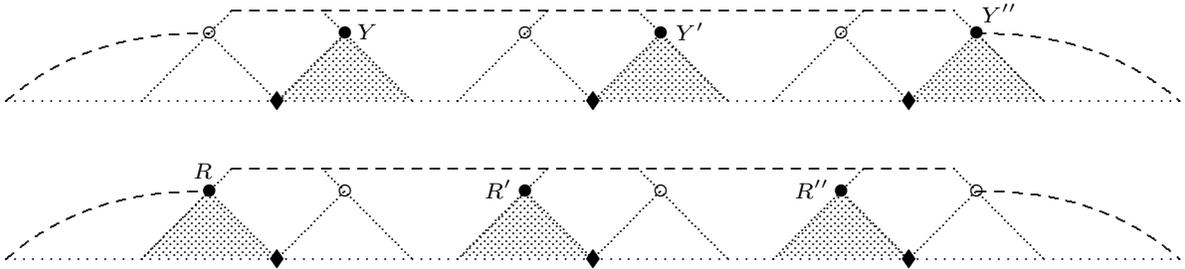
$$
{\beginpicture
    \setcoordinatesystem units <.25cm,.25cm>

%%%%%%%%%%%%%%%%%%%%%%%%%%%%%%%%%%%%%%%%%%%%%
\put{\beginpicture
\multiput{} at -6 0  46 5 /
\multiput{$\ss\blacklozenge$} at 6 0  20 0  34 0 /
\setdots <1mm>
\plot -6 0  46 0 /
\setdots <.5mm>
\plot 0 0  4 4  /
\plot 3 3  6 0  9 3 /
\plot 8 4  12 0 /

\plot 14 0  18 4 /
\plot 17 3  20 0  23 3 /
\plot 22 4  26 0 /

\plot 28 0  32 4 /
\plot 31 3  34 0  37 3 /
\plot 36 4  40 0 /

\multiput{$\bullet$} at 9 3  23 3  37 3 /
\multiput{$\circ$} at 3 3  17 3  31 3 /
\put{$\ss Y$} at 10 3.1
\put{$\ss Y'$} at 24.3 3.1
\put{$\ss Y''$} at 38 3.9

\setdashes <1mm>
\setquadratic
\plot -6 0  -2 2.3  3 3 /
\plot 37 3  42 2.3  46 0 /

\setlinear
\plot 4 4  36 4 /
\setdots <1mm>

\setshadegrid span <.4mm>
\hshade 0 6 12 <z,z,,> 3 9 9   /
\hshade 0 20 26 <z,z,,> 3 23 23  /
\hshade 0 34 40 <z,z,,> 3 37 37  /
\endpicture} at 0 0

%%%%%%%%%%%%%%%%%%%%%%%%%%%%%%%%%%%%%%%%%%%%%
\put{\beginpicture
\multiput{} at -6 0  46 5 /
\multiput{$\ss\blacklozenge$} at 6 0  20 0  34 0 /
\setdots <1mm>
\plot -6 0  46 0 /
\setdots <.5mm>
\plot 0 0  4 4  /
\plot 3 3  6 0  9 3 /
\plot 8 4  12 0 /

\plot 14 0  18 4 /
\plot 17 3  20 0  23 3 /
\plot 22 4  26 0 /

\plot 28 0  32 4 /
\plot 31 3  34 0  37 3 /
\plot 36 4  40 0 /

\multiput{$\circ$} at 9 3  23 3  37 3 /
\multiput{$\bullet$} at 3 3  17 3  31 3 /
\put{$\ss R$} at 2.7 3.9
\put{$\ss R'$} at 15.8 3.1
\put{$\ss R''$} at 29.7 3.1

\setdashes <1mm>
\setquadratic
\plot -6 0  -2 2.3  3 3 /
\plot 37 3  42 2.3  46 0 /

\setlinear
\plot 4 4  36 4 /
\setdots <1mm>

\setshadegrid span <.4mm>
\hshade 0 0 6 <z,z,,> 3 3 3   /
\hshade 0 14 20 <z,z,,> 3 17 17  /
\hshade 0 28 34 <z,z,,> 3 31 31  /
\endpicture} at 0 -7

\endpicture}
$$

Proof of Proposition. We assume that $A$ is a concave Nakayama algebra of height $m+1$ which is $d$-closed 
and has global dimension at most $d$.
	\smallskip
(1) implies (2$''$). Here, we assume in addition that $A$ has a higher Auslander algebra
of global dimension $d$. 

Let $Y''$ be the principal cliff module. It has length $m$ and is
indecomposable, injective and not projective.
Thus $Y''$ is plus-strictly-increasing, see 8.2.
	\smallskip
(2$''$) implies (2$'$). Trivial.
	\smallskip
(2$'$) implies (2). Let $Y'$ be plus-strictly-increasing of length $m$ with $\char Y' =
(z'_1,\dots,z'_m).$ Let $t = \frac12(z'_2-1)$ and $s = (m+1)t.$ We consider the
pile with cliff $Y'$ and $s$ summits. Let $Y$ be its radical 
The shift lemma asserts that $\char Y = \char Y' - 2t(1,\dots,1),$ thus
$Y$ is plus-strictly-increasing with $\char Y = (z'_1-2t,\dots,z'_m-2t)$, where
$z_2 = z'_2-2t = 1.$ 
	\smallskip
The equivalences of (2) and (3), of (2$'$) and (3$')$ as well as of (2$''$) and  (3$''$)
are given by Lemma 8.4.
	\smallskip
The equivalence of (3) and (4) is trivial.
	\smallskip
(4) implies (1): see Proposition 8.3.
	\medskip
%%%%%%%%%%%%%%%%%%%%%%%%%%%%%%%%%%%%%%%%%%%%%%%%%%%%%%%%%%%%%%%%%%%%%%%%%%
Finally, let us deal with condition (5). First, assume that there exists an
indecomposable module of length $m$ which is odd. Then $N$ is not projective and
$\pi N$ is plus-strictly-increasing, thus the condition (3$'$) is satisfied.

However there are many concave higher Auslander algebras with even global dimension and
height $m+1$, which have no odd indecomposable module of length $m$, for example
$H_4(3)$ or $H_4(1,3).$ 
$\s$
	\medskip
{\bf 8.7. Proof of Theorems 1$'$, 2$'$, and 3$'$.} Let $c_1 < \cdots < c_u \le d$
be odd numbers and $H =  C_dA(c_1,\dots,c_u,0,1).$
According to 8.5, $H$ is a concave higher Auslander algebra of global dimension $d$
and any such algebra is obtained in this way. 
	\smallskip
Let $P$ be the first summit of $H$, and $R = \rad P$. By construction, we have $\char P =
(c_1,\dots,c_u,0,1)$ and $\char R = (c_1,\dots,c_u,0)$. Let $Q$ be the last summit
and $Y = Q/\soc Q.$ We have calculated $\char Y$ in 8.4 (a). It follows that
$\pd \soc Q = d-1,$ since $\soc Q = \Omega Y$. In this way, we obtain the explicit
form of $\char Q$ as mentioned.
$\s$
	\medskip
{\bf 8.7.} 
The characterizations presented in 8.5 refer to the existence of idecomposable
modules of length $h-1$, where $h$ is the height of the algebra. Actually, it
is sufficient to assume the existence of corresponding sequences of simple modules.
	\medskip
{\bf Proposition.} 
{\it Let $d$ be even and $m\ge 2$. 
Let $A$ be a concave Nakayama algebra of height $m+1$ and global dimension $d$.
The algebra $A$ is a higher Auslander algebra iff $A$ is $d$-closed and there exists
a simple module $S$ such that the numbers $z_i = \pd \tau^iS$ with $0\le i \le m-1$
satisfy the following conditions: $z_0,\dots,z_{m-2}$ are odd, $z_{m-1}$ is even, 
and $z_{m-2} < z_{m-2}< \cdots < z_{0} < z_{m-1}.$}
	\medskip
Proof. According to 1.1 (4), we have $|PS| \ge m.$ If we would have $|PS| = m,$ 
then $\soc PS = \tau^{m-1}S$. Thus $\pd \soc PS = z_{m-1}$ and therefore $\pd (PS/\soc PS) =
z_{m-1}+1$, whereas the maximum principle asserts that $\pd (PS/\soc PS) = z_0,$
a contradiction. This shows that $|PS| = m+1$. Let $Y = PS/\soc PS$. Then $Y$
is plus-strictly-increasing (with $\char Y = (z_{m-1},\dots, z_0)$), thus condition (2$'$) of 8.5
is satisfied.
$\s$
	\bigskip\bigskip

%%%%%%%%%%%%%%%%%%%%%%%%%%%%%%%%%%%%%%%%%%%%%%%%%%%%%%%%%%%%%%%%%%%%%%%%%%
%%%%%%%%%%%%%%%%%%%%%%%%%%%%%%%%%%%%%%%%%%%%%%%%%%%%%%%%%%%%%%%%%%%%%
{\bf 9. Final remark: Parity.} 
	\medskip
The paper  shows that the 
higher Auslander algebras $A$ with global dimension $d(A)$
have a
quite extreme homological behavior, and that the extreme conditions 
satisfied by higher Auslander algebras $A$ depend on the parity of $d(A)$:
they are of completely different nature, depending on whether
$d(A)$ is even or odd.

Let $A$ be a higher Auslander algebra. Let $I$ be indecomposable injective. Then
$\pd I \in \{0,d(A)\}$. 
For $d$ even, there are many Nakayama algebras with
$\pd I \in \{0,d\}$ for all indecomposable injective modules (for example, this
happens always, if the Kupisch series is of the form $(2m,2m+1)$, but for $m \ge 2$, these
algebras have infinite global dimension), whereas for $d$ odd, this condition
implies that $A$ has finite global dimension.

Similarly, if $A$ is a higher Auslander algebra and $S$ is simple, then $S$ is torsionless
or else $\pd S = d(A).$ 
For $d$ odd, there are many Nakayama algebras such that any simple module is torsionless
or has projective dimension $d$ (for example, if
the Kupisch series is of the form $(2m-1,2m)$, then $A$ has two simple modules,
one of them is torsionless, the other has projective dimension $1$, and for $m\ge 2$,
the algebra $A$ has infinite global dimension). On the other hand, 
if $A$ is a Nakayama algebra 
such that any simple module is torsionless or has even projective dimension, then
$A$ has to be self-injective. 
	\bigskip
%%%%%%%%%%%%%%%%%%%%%%%%%%%%%%%%%%%%%%%%%%%%%%%%%%%%%%%%%%%%%%%%%%%%%
{\bf Acknowledgment.} The author thanks the referee for a very careful reading of
a first version, pointing out a huge number of inacurracies and misprints. In fact,
on the basis of these comments, the paper was completely rewritten. 
	\bigskip\bigskip
%%%%%%%%%%%%%%%%%%%%%%%%%%%%%%%%%%%%%%%%%%%%%%%%%%%%%%%%%%%%%%%%%%%%%
{\bf 10. References.} 
	\medskip
\item{[1]} A\. Chan, O\. Iyama, R\. Marczinzik. Auslander-Gorenstein algebras from Serre-formal
  algebras via replication. Adv. Math. 345 (2019), 222--262.
\item{[2]} O\. Iyama. Auslander correspondence. 
   Adv. Math. 210 (2007). 51--82.
\item{[3]} O\. Iyama. Cluster tilting for higher Auslander algebras.
   Advances in Mathematics 226 (2011), 1--61.
\item{[4]} D\. O\. Madsen. Projective dimensions and Nakayama
   algebras. Fields Institute Communications. 45.
   Amer\. Math\. Soc\., Providence, RI (2005). 247--265.
\item{[5]} D\. O\. Madsen, R\. Marczinzik, G\. Zaimi.
  On the classification of higher Auslander algebras for Nakayama algebras.
  Journal of Algebra 556 (2020), 776--805.
\item{[6]} C\. M\. Ringel.
    The finitistic dimension of a Nakayama algebra.
    Journal of Algebra (2021).
    https://doi.org/10.1016/j.jalgebra.2021.01.040
\item{[7]} E\. Sen. Nakayama Algebras which are Higher Auslander Algebras. (2020).
   \newline    arXiv:2009.03383.
\item{[8]} L. Vaso. $n$-cluster tilting subcategories of representation-directed 
   algebras. (2017) arXiv:1705.01031
	\bigskip\medskip
%%%%%%%%%%%%%%%%%%%%%%%%%%%%%%%%%%%%%%%%%%%%%%%%%%%%%%%
{\baselineskip=1pt
\rmk
C. M. Ringel\par
Fakult\"at f\"ur Mathematik, Universit\"at Bielefeld \par
POBox 100131, D-33501 Bielefeld, Germany  \par
ringel\@math.uni-bielefeld.de
\par}

\vfill\eject
{\bf Some examples:} the concave higher Auslander algebras $H = H_d(\bold c)$ of height $2,3,4$ with $d = 3$ 
and $d = 4.$ The set of subfactors of $Z(H)$ or $Z'(H)$ have been shaded.
	\medskip
First, for $d = 3$.
$$
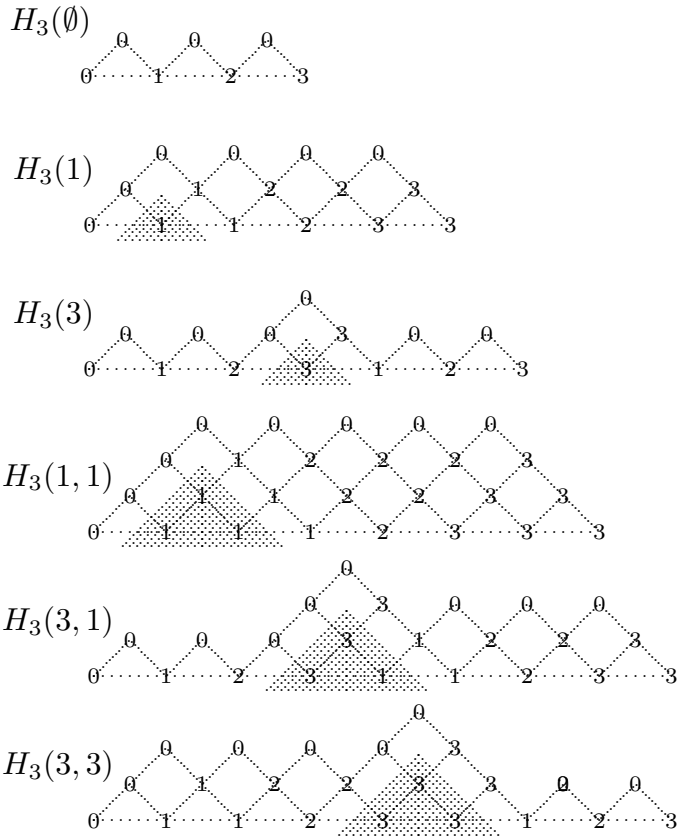
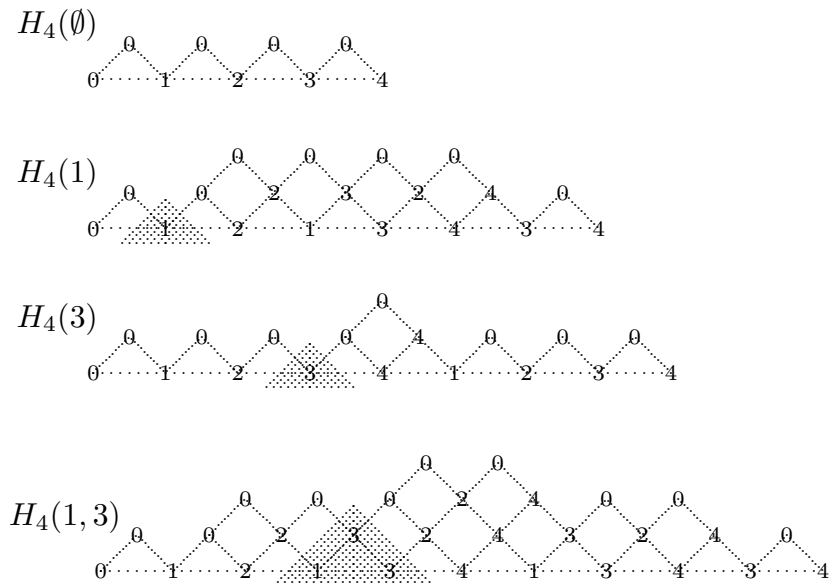
{\beginpicture
    \setcoordinatesystem units <.4cm,.4cm>
%%%%%%%%%%%%%%%%%%%%%%%%%%%%%%%%%%%%%%%%%%%%%%%%%%H_3()
\put{\beginpicture
\put{$H_3(\emptyset)$} at -1 1.5
\multiput{} at 0 0  6 1   /
\setdots <1mm>
\plot 0 0  6 0 /
\setdots <0.5mm>
\plot 0 0  1 1  2 0  3 1  4 0  5 1  6 0 /
\multiput{$\ss 0$} at 0 0  1 1  3 1  5 1 /
\multiput{$\ss 1$} at 2 0 /
\multiput{$\ss 2$} at 4 0 /
\multiput{$\ss 3$} at 6 0 /
\endpicture} at 0 0
%%%%%%%%%%%%%%%%%%%%%%%%%%%%%%%%%%%%%%%%%%%%%%%%%H_3(1)
\put{\beginpicture
\put{$H_3(1)$} at -1 1.5
\multiput{} at 0 0  12 2   /
\setdots <1mm>
\plot 0 0  10 0 /
\setdots <0.5mm>
\plot 0 0  1 1  2 0  3 1  4 0  5 1  6 0  7 1  8 0  9 1  10 0    /
\plot 1 1  2 2  3 1  4 2  5 1  6 2  7 1  8 2  9 1   /
%\plot 2 2  3 3  4 2  5 3  6 2  7 3  8 2  9 3  10 2  11 3  12 2 /
\multiput{$\ss 0$} at 0 0  1 1  2 2  4 2  6 2  8 2    /
\multiput{$\ss 1$} at 2 0  3 1  4 0     /
\multiput{$\ss 2$} at 5 1  6 0  7 1     /
\multiput{$\ss 3$} at 8 0  9 1  10 0   /
\multiput{$\ss 4$} at /
%\put{$\bigcirc$} at 2 0
\setshadegrid span <.4mm>
\vshade .5 -.5 -.5 <z,z,,> 2 -.5 1 <z,z,,> 3.5 -.5 -.5 /
\endpicture} at 3  -4
%%%%%%%%%%%%%%%%%%%%%%%%%%%%%%%%%%%%%%%%%%%%%%%%%%H_3(3)
\put{\beginpicture
\put{$H_3(3)$} at 1 1.5
\multiput{} at 0 0  16 2  /
\setdots <1mm>
\plot 2 0  14 0 /
\setdots <0.5mm>
\plot 2 0  3 1  4 0  5 1  6 0  7 1  8 0  9 1  10 0  11 1  12 0  13 1  14 0   /
\plot 7 1  8 2  9 1  /
\multiput{$\ss 0$} at  2 0  3 1  5 1  7 1  8 2  11 1  13 1   /
\multiput{$\ss 1$} at 4 0  10 0 /
\multiput{$\ss 2$} at 6 0  12 0  /
\multiput{$\ss 3$} at 8 0  9 1  14 0 /
%\put{$\bigcirc$} at 8 0
\setshadegrid span <.4mm>
\vshade 6.5 -.5 -.5 <z,z,,> 8 -.5 1 <z,z,,> 9.5 -.5 -.5 /
\endpicture} at 4 -8

%%%%%%%%%%%%%%%%%%%%%%%%%%%%%%%%%%%%%%%%%%%%%%%%%%H_3(1,1)
\put{\beginpicture
\put{$H_3(1,1)$} at -1 1.5
\multiput{} at 0 0  12 2   /
\setdots <1mm>
\plot 0 0  14 0 /
\setdots <0.5mm>
\plot 0 0  1 1  2 0  3 1  4 0  5 1  6 0  7 1  8 0  9 1  10 0  11 1  12 0  13 1  14 0 /
\plot 1 1  2 2  3 1  4 2  5 1  6 2  7 1  8 2  9 1  10 2  11 1  12 2  13 1 /
\plot 2 2  3 3  4 2  5 3  6 2  7 3  8 2  9 3  10 2  11 3  12 2 /
\multiput{$\ss 0$} at 0 0  1 1  2 2  3 3  5 3  7 3  9 3  11 3   /
\multiput{$\ss 1$} at 2 0  3 1  4 0  4 2  5 1  6 0   /
\multiput{$\ss 2$} at 6 2  7 1  8 2  8 0  9 1  10 2     /
\multiput{$\ss 3$} at 10 0  11 1  12 2  12 0  13 1  14 0   /
\multiput{$\ss 4$} at /
%\put{$\bigcirc$} at 3 1
\setshadegrid span <.4mm>
\vshade .5 -.5 -.5 <z,z,,> 3 -.5 2 <z,z,,> 5.5 -.5 -.5 /

\endpicture} at 4 -12

%%%%%%%%%%%%%%%%%%%%%%%%%%%%%%%%%%%%%%%%%%%%%%%%%%H_3(3,1)
\put{\beginpicture
\put{$H_3(3,1)$} at -1 1.5
\multiput{} at 0 0  12 2   /
\setdots <1mm>
\plot 0 0  16 0 /
\setdots <0.5mm>
\plot 0 0  1 1  2 0  3 1  4 0  5 1  6 0  7 1  8 0  9 1  10 0  11 1  12 0  13 1  14 0
  15 1  16 0 /
\plot  5 1  6 2  7 1  8 2  9 1  10 2  11 1  12 2  13 1  14 2  15 1 /
\plot  6 2  7 3  8 2  /
\multiput{$\ss 0$} at 0 0  1 1  3 1  5 1  6 2  7 3  10 2  12 2  14 2  /
\multiput{$\ss 1$} at  2 0  8 0  9 1  10 0   /
\multiput{$\ss 2$} at  4 0  11 1  12 0  13 1     /
\multiput{$\ss 3$} at  6 0  7 1  8 2  14 0  15 1  16 0 /
\multiput{$\ss 4$} at /
%\put{$\bigcirc$} at 7 1
\setshadegrid span <.4mm>
\vshade 4.5 -.5 -.5 <z,z,,> 7 -.5 2 <z,z,,> 9.5 -.5 -.5 /

\endpicture} at 5 -16

%%%%%%%%%%%%%%%%%%%%%%%%%%%%%%%%%%%%%%%%%%%%%%%%%%H_3(3,3)
\put{\beginpicture
\put{$H_3(3,3)$} at -1 1.5
\multiput{} at 0 0  12 2   /
\setdots <1mm>
\plot 0 0  16 0 /
\setdots <0.5mm>
\plot 0 0  1 1  2 0  3 1  4 0  5 1  6 0  7 1  8 0  9 1  10 0  11 1  12 0  13 1  14 0
  15 1  16 0 /
\plot 1 1  2 2  3 1  4 2  5 1  6 2  7 1  8 2  9 1  10 2  11 1  /
\plot 8 2  9 3  10 2 /
\multiput{$\ss 0$} at 0 0  1 1  2 2  4 2  6 2  8 2  9 3  13 1  15 1   /
\multiput{$\ss 1$} at 2 0  3 1  4 0  12 0  /
\multiput{$\ss 2$} at 5 1  6 0  7 1  13 1  14 0 /
\multiput{$\ss 3$} at 8 0  9 1  10 2  10 0  11 1  16 0 /
\multiput{$\ss 4$} at /
%\put{$\bigcirc$} at 9 1
\setshadegrid span <.4mm>
\vshade 6.5 -.5 -.5 <z,z,,> 9 -.5 2 <z,z,,> 11.5 -.5 -.5 /
\endpicture} at 5 -20

\put{} at 20 0
\endpicture}
$$
	\bigskip
Second, for $d = 4$
$$
{\beginpicture
    \setcoordinatesystem units <.4cm,.4cm>
%%%%%%%%%%%%%%%%%%%%%%%%%%%%%%%%%%%%%%%%%%%%%%%%%%H_4()
\put{} at 19 0

\put{\beginpicture
\put{$H_4(\emptyset)$} at -1 1.5
\multiput{} at 0 0  8 1   /
\setdots <1mm>
\plot 0 0  8 0 /
\setdots <0.5mm>
\plot 0 0  1 1  2 0  3 1  4 0  5 1  6 0  7 1  8 0  /
\multiput{$\ss 0$} at 0 0  1 1  3 1  5 1  7 1 /
\multiput{$\ss 1$} at 2 0 /
\multiput{$\ss 2$} at 4 0 /
\multiput{$\ss 3$} at 6 0 /
\multiput{$\ss 4$} at 8 0 /
\endpicture} at 0 0
%%%%%%%%%%%%%%%%%%%%%%%%%%%%%%%%%%%%%%%%%%%%%%%%%%H_4(1)
\put{\beginpicture
\put{$H_4(1)$} at -1 1.5
\multiput{} at 0 0  8 2  /
\setdots <1mm>
\plot 0 0  14 0 /
\setdots <0.5mm>
\plot 0 0  1 1  2 0  3 1  4 0  5 1  6 0  7 1  8 0  9 1  10 0  11 1  12 0  13 1  14 0  /
\plot 3 1  4 2  5 1  6 2  7 1  8 2  9 1  10 2  11 1 /
\multiput{$\ss 0$} at 0 0  1 1  3 1  4 2  6 2  8 2  10 2  13 1  /
\multiput{$\ss 1$} at 2 0  6 0 /
\multiput{$\ss 2$} at 4 0  5 1  9 1 /
\multiput{$\ss 3$} at 7 1  8 0  12 0 /
\multiput{$\ss 4$} at 10 0  11 1  14 0 /
%\put{$\bigcirc$} at 4 0
\setshadegrid span <.4mm>
\vshade .5 -.5 -.5 <z,z,,> 2 -.5 1 <z,z,,> 3.5 -.5 -.5 /
\endpicture} at 3 -4
%%%%%%%%%%%%%%%%%%%%%%%%%%%%%%%%%%%%%%%%%%%%%%%%%%H_4(3)
\put{\beginpicture
\put{$H_4(3)$} at -1 1.5
\multiput{} at 0 0  16 2  /
\setdots <1mm>
\plot 0 0  16 0 /
\setdots <0.5mm>
\plot 0 0  1 1  2 0  3 1  4 0  5 1  6 0  7 1  8 0  9 1  10 0  11 1  12 0  13 1  14 0  15 1  16 0 /
\plot 7 1  8 2  9 1  /
\multiput{$\ss 0$} at 0 0  1 1  3 1  5 1  7 1  8 2  11 1  13 1  15 1   /
\multiput{$\ss 1$} at 2 0  10 0  /
\multiput{$\ss 2$} at 4 0  12 0  /
\multiput{$\ss 3$} at 6 0  14 0 /
\multiput{$\ss 4$} at 8 0  9 1  16 0 /
%\put{$\bigcirc$} at 4 0
\setshadegrid span <.4mm>
\vshade 4.5 -.5 -.5 <z,z,,> 6 -.5 1 <z,z,,> 7.5 -.5 -.5 /
\endpicture} at 4 -8
%%%%%%%%%%%%%%%%%%%%%%%%%%%%%%%%%%%%%%%%%%%%%%%%%%H_4(1,3)
\put{\beginpicture
\put{$H_4(1,3)$} at -1 1.5
\multiput{} at 0 0  20 2   /
\setdots <1mm>
\plot 0 0  20 0 /
\setdots <0.5mm>
\plot 0 0  1 1  2 0  3 1  4 0  5 1  6 0  7 1  8 0  9 1  10 0  11 1  12 0  13 1  14 0
   15 1  16 0  17 1  18 0  19 1  20 0 /
\plot  3 1  4 2  5 1  6 2  7 1  8 2  9 1  10 2  11 1  12 2  13 1
   14 2  15 1  16 2  17 1  /
\plot  8 2  9 3  10 2  11 3  12 2 /
\multiput{$\ss 0$} at 0 0  1 1  3 1  4 2  6 2  8 2  9 3  11 3  14 2  16 2  19 1  /
\multiput{$\ss 1$} at  2 0  6 0  12 0 /
\multiput{$\ss 2$} at  4 0  5 1  9 1  10 2  15 1  /
\multiput{$\ss 3$} at  7 1  8 0  13 1  14 0  18 0  /
\multiput{$\ss 4$} at  10 0  11 1  12 2  16 0  17 1  20 0 /
%\put{$\bigcirc$} at 9 1
\setshadegrid span <.4mm>
\vshade 4.5 -.5 -.5 <z,z,,> 7 -.5 2 <z,z,,> 9.5 -.5 -.5 /
\endpicture} at 6 -13
\endpicture}
$$
	\bigskip
In general, for arbitrary height $h$, there are $h-1$ algebras
$H_3(\bold c)$ of height $h$. But there are no
algebras $H_4(\bold c)$ of height greater than $4$.

%%%%%%%%%%%%%%%%%%%%%%%%%%%%%%%%%%%%%%%%%%%%%%%%%%%%%%%%%%%%%%%%%%%%%%%%%%%%%%%
%%%%%%%%%%%%%%%%%%%%%%%%%%%%%%%%%%%%%%%%%%%%%%%%%%%%%%%%%%%%%%%%%%%%%%%%%%%%%%%
\bye